\documentclass[reqno]{amsart}
\usepackage{mathrsfs}
\usepackage{amsmath,amsfonts,amssymb,amsthm}
\usepackage{upgreek}
\usepackage[usenames,dvipsnames]{xcolor}
\usepackage{enumitem}
\usepackage{graphicx}
\usepackage{overpic}
\usepackage{cancel}
\usepackage[outdir=./]{epstopdf}
\usepackage[colorlinks=true]{hyperref} 
\usepackage{cleveref}
\usepackage{url}
\usepackage[utf8]{inputenc}
\usepackage{booktabs,multirow}
\usepackage[font=footnotesize,width=\textwidth]{caption}
\hypersetup{linkcolor=Violet,citecolor=PineGreen}
\newtheorem{teor}{Theorem}[section]
\newtheorem{lema}[teor]{Lemma}
\newtheorem{prop}[teor]{Proposition}
\newtheorem{coro}[teor]{Corollary}
\theoremstyle{definition}

\newtheorem{exa}[teor]{Example}

\newtheorem{nota}[teor]{Remark}
\newtheorem{notas}[teor]{Remarks}
\numberwithin{equation}{section}

\newcommand{\R}{\mathbb R}

\newcommand{\Z}{\mathbb{Z}}

\newcommand{\N}{\mathbb{N}}

\newcommand{\mA}{\mathcal{A}}
\newcommand{\mB}{\mathcal{B}}
\newcommand{\mC}{\mathcal{C}}

\newcommand{\mK}{\mathcal{K}}
\newcommand{\mM}{\mathcal{M}}

\newcommand{\mU}{\mathcal{U}}
\newcommand{\mV}{\mathcal{V}}
\newcommand{\ep}{\varepsilon}

\newcommand{\mI}{\mathcal{I}}
\newcommand{\mJ}{\mathcal{J}}

\newcommand{\W}{\Omega}
\newcommand{\WR}{\W\times\R}
\newcommand{\RR}{\R\times\R}
\newcommand{\w}{\omega}

\newcommand{\ma}{\mathfrak{a}}
\newcommand{\mb}{\mathfrak{b}}

\newcommand{\mr}{\mathfrak{r}}
\newcommand{\ml}{\mathfrak{l}}
\newcommand{\mm}{\mathfrak{m}}
\newcommand{\muk}{\mathfrak{u}}
\newcommand{\mf}{\mathfrak{f}}
\newcommand{\mg}{\mathfrak{g}}
\newcommand{\mh}{\mathfrak{h}}
\newcommand{\tma}{\tilde\ma}
\newcommand{\tmr}{\tilde\mr}
\newcommand{\tml}{\tilde\ml}
\newcommand{\tmm}{\tilde\mm}
\newcommand{\tmuk}{\tilde\muk}

\newcommand{\upalfa}{$\upalpha$}
\newcommand{\upomeg}{$\upomega$}

\newcommand{\G}{\Gamma}

\newcommand{\ws}{\w{\cdot}s}
\newcommand{\wt}{\w{\cdot}t}
\newcommand{\bwt}{\bar\w{\cdot}t}

\newcommand{\n}[1]{\left\|#1\right\|}

\newcommand{\lsm}{\left[\begin{smallmatrix}}
\newcommand{\rsm}{\end{smallmatrix}\right]}

\newcommand{\merg}{\mathfrak{M}_\mathrm{erg}(\W,\sigma)}

\begin{document}
\title[Critical transitions with applications in Ecology]
{Critical transitions for asymptotically concave or d-concave
nonautonomous differential equations with applications in Ecology}
\author[J. Due\~{n}as]{Jes\'{u}s Due\~{n}as}
\author[C. N\'{u}\~{n}ez]{Carmen N\'{u}\~{n}ez}
\author[R. Obaya]{Rafael Obaya}
\address{{\rm(J.~Due\~{n}as, C.~N\'{u}\~{n}ez, R.~Obaya)}. Departamento de Matem\'{a}tica Aplicada, Universidad de
Va\-lladolid, Paseo Prado de la Magdalena 3-5, 47011 Valladolid, Spain. The authors are members of
IMUVA: Instituto de Investigaci\'{o}n en Matem\'{a}ticas, Universidad de Valladolid.}
\email[J.~Due\~{n}as]{jesus.duenas@uva.es}
\email[C.~N\'{u}\~{n}ez]{carmen.nunez@uva.es}
\email[R.~Obaya]{rafael.obaya@uva.es}
\thanks{All the authors were supported by Ministerio de Ciencia, Innovaci\'{o}n y Universidades (Spain)
under project PID2021-125446NB-I00 and by Universidad de Valladolid under project PIP-TCESC-2020.
J.~Due\~{n}as was also supported by Ministerio de Universidades (Spain) under programme FPU20/01627.}
\date{}
\begin{abstract}
The occurrence of tracking or tipping situations for a transition equation $x'=f(t,x,\G(t,x))$
is analyzed under the assumptions on concavity in $x$ either of the maps giving rise
to the asymptotic equations $x'=f(t,x,\G_\pm(t,x))$ or of their derivatives with respect
to the state variable (d-concavity), but without assuming these conditions
on the transition equation itself. The approaching condition is just
$\lim_{t\to\pm\infty}(\G(t,x)-\G_\pm(t,x))=0$ uniformly on compact real sets, and so
there is no restriction to the dependence on time of the limit equations. The analysis provides
a powerful tool to analyze the occurrence of critical transitions for one-parametric
families $x'=f(t,x,\G^c_\pm(t,x))$. The new approach significatively widens the field
of application of the results, since the evolution law of the transition
equation can be essentially different from those of the limit equations.
Among these applications, some scalar population dynamics models subject
to non trivial predation and migration patterns are analyzed, both theoretically and numerically.
\par
Some key points in the proofs are: to understand the transition equation
as part of an orbit in its hull which approaches the \upalfa-limit and
\upomeg-limit sets; to observe that these sets concentrate all the ergodic measures;
and to prove that in order to describe the dynamical possibilities of the equation
it suffices that the concavity or d-concavity conditions hold for a complete measure subset of the
equations of the hull.
\end{abstract}
\renewcommand{\subjclassname}{\textup{2020} Mathematics Subject Classification}

\maketitle
\section{Introduction}
{\em Tipping points\/} or {\em critical transitions\/}
are significant nonlinear phenomena that occur in complex
systems subject to smooth changes of the external conditions. Roughly speaking, they
are sudden and often irreversible
changes in the state of the system caused by small changes in
the external input. During the last years, they have frequently appeared in the literature
as an explanation of abrupt changes in climate \cite{aajq,awvc,lenton1,schellnhuber}, ecology
\cite{altw,scheffer2,scheffer1,vanselow}, biology \cite{hill1,nene} or finances \cite{may1,yukalov},
among other scientific areas of great interest. For this reason, critical transitions have
become an important topic of multidisciplinary research.
\par
A branch of the mathematical formulation of this problem focuses on one-parametric
ordinary differential equations (ODEs). The parameter is replaced by
a map (a {\em parameter shift}) with constant asymptotic limits as $t\to\pm\infty$, and
the two {\em limit equations},
given by these constant limits of the parameter, are frequently assumed to have the same
type of global dynamics. The initial ODE is understood as a transition between the
past equation and the future equation. Sometimes, the dynamics of this transition equation
reproduces those of the asymptotic limits and approaches them as time decreases and
increases. This situation is usually referred to as {\em tracking}, and
the remaining situations as {\em tipping}. In general, tracking means the survival
over time of a local pullback attractor which represents the desirable state
of the system, and hence tipping may mean a catastrophe.
If the parameter shift, itself, depends on a parameter,
a {\em critical value\/} of the parameter, resulting in a {\em critical transition},
occurs if there is tracking to its left and tipping to its right, or viceversa.
\par
This is the approach of the reference work \cite{aspw}, as well as of
\cite{alas, kiersjones, keeffe1, sebastian1, wieczcompact}, among many other papers.
In the recent works \cite{dolo1,dno3,lno2,lno3,lnor,remo1}, the limit equations are allowed to be nonautonomous,
and the law of the (scalar) ODE is assumed either to be concave with respect to the state variable,
or to have concave derivative (also with respect to the state).
This paper delves into these approaches, which are now optimized in several senses that
considerably widen the field of application of the results.
\par
Let us briefly describe this much more general setting.
We work with a nonautonomous scalar ODE
\begin{equation}\label{eq:1intro}
 x'=g(t,x)\,,
\end{equation}
given by a regular enough map $g$ for which the hull $\W_g$, defined
as the closure of the set of time shifts $g_s(t,x):=g(t+s,x)$ in the compact-open
topology of $C(\RR,\R)$, is compact. Equation \eqref{eq:1intro} is
embedded in the family
\begin{equation}\label{eq:1hull}
 x'=\w(t,x)\,,\qquad \w\in\W_g\,,
\end{equation}
which defines a skewproduct flow on $\W_g\times\R$. The set $\W_g$ is composed by three
(possibly nondisjoint) sets: $\{g_s\,|\;s\in\R\}$ and the corresponding
\upalfa-limit and \upomeg-limit sets, $\W_g^\upalpha$ and $\W_g^\upomega$; and,
as a consequence, the ergodic measures on the hull are concentrated in
$\W_g^\upalpha\cup\W_g^\upomega$. This point is key in our approach: it turns out that
assuming conditions on concavity with respect to $x$ of $\w$ (in the so-called {\em concave case})
or of its derivative $\w_x$ (in the {\em d-concave case}) for all $\w\in\W_g^\upalpha\cup\W_g^\upomega$
suffices to determine the maximum number of hyperbolic solutions for \eqref{eq:1intro}:
two in the concave case, and three in the d-concave case. In addition, if they exist,
and if an extra coercivity assumption on all the involved equations holds
and hence the set of bounded solutions is bounded (if nonempty),
then the hyperbolic solutions yield a very fixed type of global dynamics.
These assertions are consequences of general results on skewproducts
that we prove in Sections \ref{sec:3} and \ref{sec:5}
in the concave and d-concave cases, respectively. In our opinion,
the skewproduct formalism provides the most suitable framework to understand
the occurrence of critical transitions in nonautonomous models, and this is one of the main
contributions of this paper. We point out that
no concavity condition is required on $g$ or on $g_x$, which is a significative
difference with respect to the mentioned previous approaches.
\par
In order to formulate the necessary conditions on the \upalfa-limit and \upomeg-limit sets
in the language of processes rather than in the language of skewproducts, and thus
be able to apply our results to particular examples, we assume the existence of
two regular enough functions $g_\pm$ with $\lim_{t\to\pm\infty}(g(t,x)-g_\pm(t,x))=0$ uniformly on
compact sets of $\R$, and which satisfy the required concavity and coercivity conditions. These
conditions are inherited by all the elements of the corresponding hulls $\W_{g_\pm}$, and the asymptotic
approach yields $\W_g^\upalpha=\W_{g_-}^\upalpha$ and $\W_g^\upomega=\W_{g_+}^\upomega$. So,
we are in the suitable framework.
In order to describe scenarios as rich in dynamical possibilities as possible, we
assume that the {\em limit equations\/} $x'=g_\pm(t,x)$ have the maximal number
of hyperbolic solutions. And then we describe all the
dynamical cases for \eqref{eq:1intro} (which are three)
and explain the strong connection among critical transition and nonautonomous
saddle-node bifurcation.
These are the basic contents of the beginning of Sections \ref{sec:4} (concave case)
and \ref{sec:6} (d-concave case).
\par
When trying to build realistic mathematical models, it is common to replace a
parameter that appears in a first and simple approximation to the law of evolution of a given system,
say $x'=f(t,x,\gamma)$, by a map that may depend on time, state, and new parameters.
Our approach allows us to deal with
\begin{equation}\label{eq:1introfG}
 x'=f(t,x,\G(t,x))
\end{equation}
assuming the existence of
two maps $\G_\pm$ such that $\lim_{t\to\pm\infty}(\G(t,x)-\G_\pm(t,x))=0$ uniformly on
compact sets of $\R$. Of course, we must assume conditions on $f,\G$ and $\G_\pm$
guaranteing the previous hypotheses on $g(t,x):=f(t,x,\G(t,x))$ and
$g_\pm(t,x):=f(t,x,\G_{\pm}(t,x))$. The second fundamental difference arises here.
In previous approaches, $\G$ is just a map of $t$ (and often the
unique time-dependent part of the evolution law) with constant asymptotic limits.
But, in this new formulation, also the asymptotic part $\G_\pm$ of the
past and future equations $x'=f(t,x,\G_\pm(t,x))$ may depend on $t$ and $x$.
In order to analyze the occurrence of critical transitions,
we let the transition map to depend on a parameter $c$ which moves, getting
$x'=f(t,x,\G^c(t,x))$. It is reasonable to assume that the asymptotic equations
do not depend on the parameter $c$, whose variations represent different
ways to approach the past and the future.
Subsections \ref{subsec:41} (concave case) and \ref{subsec:61} (d-concave case) describe
several scenarios of occurrence and/or absence of critical transitions, focusing
on rate-induced, phase-induced and size-induced tipping points.
\par
As said before, the lack of general requirements on
concavity of $f$ or $f_x$ with respect to $x$ combined with the
possible dependence of $\G^c_\pm$ on $t$ and $x$ significatively increases
the number of possible applications of our results. We complete the paper by
combining theoretical and numerical techniques to analyze some of these examples in
Subsections \ref{subsec:42} (in the concave case) and \ref{subsec:62} (in the d-concave case).
\par
Let us describe with some more detail the contents of the paper.
Section \ref{sec:2} summarizes several basic concepts and results used throughout the paper:
we characterize concavity and d-concavity in terms of divided differences, describe
the skewproduct construction from a scalar nonautonomous ODE, and recall some properties of
hyperbolicity and Lyapunov exponents. The remaining sections present the analysis in
the concave case (Sections \ref{sec:3} and \ref{sec:4}) and the d-concave
case (Sections \ref{sec:5} and \ref{sec:6}). Sections \ref{sec:3} and \ref{sec:5}
deal with general properties of skewproduct flows, which include those arising from
families of the type \eqref{eq:1hull}. Section \ref{sec:3} deals with the concave case.
Under some assumptions which are weaker than the strict concavity of the equations
of a set with full measure for any ergodic measure, we prove the existence of
at most two hyperbolic solutions for each one of the equations, characterize this existence
by the occurrence of two uniformly separated bounded solutions, and describe the
global dynamics in this situation. Section \ref{sec:5} contains the analogous results
for the d-concave case: now there are at most three hyperbolic solutions for each equation,
which happens if and only if there are three uniformly separated bounded solutions and
forces a certain type of global dynamics.
\par
At the beginning of Section \ref{sec:4}, with $g$ and $g_\pm$ as above described,
we assume the strict concavity on $x$ of the maps $g_\pm(t,x)$ and the existence of the
maximum number of hyperbolic solutions for $x'=g_\pm(t,x)$: two, which form an
attractor-repeller pair. Then, there are three possibilities for \eqref{eq:1intro}:
{\sc Case A}, when it also has an attractor-repeller pair which connects with that of
the past as time decreases and with that of the future as time increases
(i.e., when there is tracking); {\sc Case C}, when it has no bounded solutions; and
{\sc Case B}, when it has exactly one bounded solution, which is nonhyperbolic.
In Section \ref{sec:6}, we assume the strict concavity on $x$ of the
derivatives $(g_\pm)_x(t,x)$ and the existence of the
maximum number of hyperbolic solutions for $x'=g_\pm(t,x)$: three.
Again, the dynamical possibilities for \eqref{eq:1intro} are three:
{\sc Case A}, if it also has three hyperbolic solutions which connect with those of
the past as time decreases and with those of the future as time increases
(tracking); {\sc Case C}, if it has two hyperbolic solutions, which approach each other as time
increases; or {\sc Case B}, if it has just one hyperbolic solution. In both cases, a typical
critical transition occurs when a small variation on $g$ causes the dynamics to
move from {\sc Case A} to {\sc Case C}. We also establish nonrestrictive conditions
guaranteing the persistence under small perturbations of these two cases,
and show that a critical transition means that the (highly nonpersistent)
{\sc Case B} occurs and can be understood as a nonautonomous saddle-node bifurcation phenomenon
(see \cite{anja,dno3,lnor,nuob6}).
\par
Subsections \ref{subsec:41} and \ref{subsec:61} are centered in
equations of the type \eqref{eq:1introfG}, and hence the corresponding general
hypotheses are given for $f$, $\G^c$, and the maps $\G^c_\pm$ before mentioned.
In scenarios suitable for raising the question of the occurrence of
rate-induced, phase-induced and size-induced critical transitions, we add extra
conditions of $f$ and $\G$ ensuring the existence and/or absence of these types
of tipping points. These results, as well as the techniques used in their proofs,
are the key points to analyze some population dynamics models, which is the goal
of Subsections \ref{subsec:42} and \ref{subsec:62}: the survival of
a given population subject to emigration and predation depends on several factors,
as the speed of arrival of the predators, the quantity of them, the moment at which they arrive,
or the length of their periods of permanence in the preys' habitat.
We point out here that the purpose of the
models that we consider is not being realistic, but showing the applicability of this
new approach to the analysis of critical transitions due to
a fairly general parametric variation in a nonautonomous dynamical system.
\section{Some preliminary results}\label{sec:2}
The following subsections collect basic concepts and some general results
needed throughout the paper.  We also provide some suitable references for further
information on those results which we do not prove here.
\subsection{Concave and d-concave real functions, and divided differences}
Recall that a map $h\in C(\R,\R)$ is \emph{concave} if
$h(\alpha x_1+(1-\alpha)x_2)\geq \alpha h(x_1)+(1-\alpha)h(x_2)$
for all $x_1,x_2\in\R$ and $\alpha\in[0,1]$ (which ensures that $h'$ is nonincreasing
if $h\in C^1(\R,\R)$).
And we say that $h\in C^1(\R,\R)$ is \emph{d-concave} if $h'$ is concave.
Our next result explains that both properties can be characterized in terms of
the {\em divided differences of first\/} and {\em second order\/} of $h$, defined as
\[
 h[x_1,x_2]:=\frac{h(x_2)-h(x_1)}{x_2-x_1}\quad\text{and}\quad
 h[x_1,x_2,x_3]:=\frac{h[x_2,x_3]-h[x_1,x_2]}{x_3-x_1}\,.
\]
Both of them are invariant under any permutation of their nodes.
\begin{prop} \label{prop:concaveuna}
\begin{enumerate}[label=\rm{(\roman*)}]
\item $h\in C(\R,\R)$ is concave if and only if
$h[x_0,x_1]\geq h[x_0,x_2]$
whenever $x_1<x_2$ and $x_0\neq x_i$ for $i\in\{1,2\}$.
\item $h\in C^1(\R,\R)$ is d-concave if and only if $h[x_1,x_0,x_2]\geq h[x_1,x_0,x_3]$
whenever $x_1<x_2<x_3$ and $x_0\neq x_i$ for $i\in\{1,2,3\}$.
\end{enumerate}
\end{prop}
\begin{proof} (i) Sufficiency is proved by taking $x_0:=\alpha x_1+(1-\alpha)\,x_2$ for $\alpha\in(0,1)$.
To check necessity, we rewrite the intermediate node as a convex combination of the other two, write the parameter
(always in $(0,1)$) in terms of the nodes, and apply the definition of concavity.
(ii) See \cite[Lemma~2.1 and remark after it]{tineo1}.
\end{proof}
Let us take $x_1<x_2<x_3$ and $x_0\ne x_i$ for $i\in\{1,2,3\}$, and define
\[
\begin{split}
a(x_0,x_1,x_2)&:=h[x_0,x_1]-h[x_0,x_2]\,,\\
b(x_0,x_1,x_2,x_3)&:=h[x_1,x_0,x_2]-h[x_1,x_0,x_3]\,,
\end{split}
\]
Clearly, $\lim_{x_0\to x_i}h[x_0,x_i]= h'(x_i)$ for all $x_i\in\R$ if $h\in C^1(\R,\R)$.
Hence, there exist $\lim_{x_0\to x_i} h[x_1,x_0,x_j]$ for $i\in\{1,2,3\}$ and $j\in\{2,3\}$. We call
\begin{equation}\label{def:2ai}
\begin{split}
a_i(x_1,x_2)&:=\lim_{x_0\to x_i} a(x_0,x_1,x_2)\,,\\
b_j(x_1,x_2,x_3)&:=\lim_{x_0\to x_j} b(x_0,x_1,x_2,x_3)
\end{split}
\end{equation}
for $i\in\{1,2\}$ and $j\in\{1,2,3\}$.
Our next result establishes an equivalence
between the sign of $a_i$ (resp.~$b_i$) and the decreasing properties of
$h'$ (resp.~$h''$).
\begin{teor} \label{th:2positiveconcave}
\begin{enumerate}[label=\rm{(\roman*)}]
\item Let $h\in C^1(\R,\R)$ be concave and $x_1<x_2$. Then, for $i\in\{1,2\}$,
$a_i(x_1,x_2)\ge 0$, and $a_i(x_1,x_2)>0$ if and only of $h'(x_1)>h'(x_2)$.
\item Let $h\in C^2(\R,\R)$ be d-concave and $x_1<x_2<x_3$. Then, for
$i\in\{1,2,3\}$, $b_i(x_1,x_2,x_3)\ge 0$, and $b_i(x_1,x_2,x_3)>0$ if and only
if $h''(x_1)>h''(x_3)$.
\end{enumerate}
\end{teor}
\begin{proof}
(i) The first assertion follows from Proposition~\ref{prop:concaveuna}(i).
For the second one, we take $i=1$ and write
\begin{equation}\label{eq:2des}
 a_1(x_1,x_2)=h'(x_1)-h[x_1,x_2]=\int_0^1\big(h'(x_1)-h'(sx_1+(1-s)x_2)\big)\, ds\,.
\end{equation}
Since $h$ is $C^1$ and concave, the integrand is continuous on $s$, nonnegative
for all $s\in[0,1]$ and nonincreasing with respect to $s$. Hence, the
integral is strictly positive if and only if $h'(x_1)>h'(x_2)$.
An analogous argument proves the assertion for $i=2$.
\smallskip\par
(ii) The first assertion follows from Proposition~\ref{prop:concaveuna}(ii).
For the second one, we work in the case $i=2$. It is easy to check that
\[
\begin{split}
b_2(x_1,x_2,x_3)
&=\frac{1}{x_2-x_1}\left(h'(x_2)-\frac{x_3-x_2}{x_3-x_1}\,h[x_1,x_2]-\frac{x_2-x_1}{x_3-x_1}\,h[x_2,x_3]\right)\\
&=\frac{1}{x_2-x_1}\left(\frac{x_3-x_2}{x_3-x_1}\int_0^1\left(h'(x_2)-h'(sx_1+(1-s)x_2)\right)\, ds\right.\\
&\qquad\qquad\qquad\left.+\frac{x_2-x_1}{x_3-x_1}\int_0^1 \left(h'(x_2)-h'(sx_3+(1-s)x_2)\right)\, ds\right)\\
&=\frac{x_3-x_2}{x_3-x_1}\int_0^1\int_0^1 s\big(h''(sx_1+(1-s)x_2+ts(x_2-x_1))\\
&\qquad\qquad\qquad\qquad\quad-h''(sx_3+(1-s)x_2-ts(x_3-x_2))\big)\, dt \,ds\,.
\end{split}
\]
Since $h$ is $C^2$ and d-concave, the integrand is continuous on $s,t$ and nonnegative for all $t\in[0,1]$.
If $h''(x_1)>h''(x_3)$, then the integrand is strictly positive at $(t,s)=(0,1)$, and hence $b_2(x_1,x_2,x_3)>0$.
Conversely, if $h''(x_1)=h''(x_3)$, then $h''$ is constant on $[x_1,x_3]$, and hence the integrand is identically zero.
We proceed analogously with $b_1$ and $b_3$.
\end{proof}
\subsection{Skew-product flows}\label{subsec:skew}
Throughout the paper, the basic concepts of flows, orbits, invariant sets, ergodic measures,
\upalfa-limit
sets and \upomeg-limit sets will be used. Their well-known definitions
and some basic properties can be found, e.g., in \cite{cano}.
\par
Let $\W$ be a compact metric space
and $\sigma\colon\R\times\W\to\W$, $(t,\w)\mapsto\sigma(t,\w)=:\w{\cdot}t$
a global continuous flow on $\W$. Along the paper, $C^{0,1}(\WR,\R)$
represents the set of continuous functions $\mh\colon\WR\to\R$ for which the
derivative $\mh_x$ with respect to the second variable exists and is continuous, and
$C^{0,2}(\WR,\R)$ is the subset of $C^{0,1}(\WR,\R)$ of maps $\mh$ for which
the second derivative $\mh_{xx}$ exists and is continuous. Given $\mh\in C^{0,1}(\WR,\R)$, we consider the
family of scalar nonautonomous differential equations
\begin{equation}\label{eq:2fam}
x'=\mh(\w{\cdot}t,x)\,,\quad\w\in\W\,.
\end{equation}
For each $\w\in\W$, \eqref{eq:2fam}$_\w$ is the particular equation of the family. We use
similar notation along the paper to refer to elements of the hull or parameters.
For each $\w\in\W$ and $x\in\R$, $(\alpha_{\w,x},\beta_{\w,x})\to\R,\,t\mapsto v(t,\w,x)$ is
the maximal solution of \eqref{eq:2fam}$_\w$ with $v(0,\w,x)=x$. Throughout the paper, any
solution will be assumed to be maximal. By uniqueness of solutions,
$v(t+s,\w,x)=v(t,\ws,v(s,\w,x))$ when the right-hand term is defined.
Hence, if $\mV:=\bigcup_{(\w,x)\in\WR}((\alpha_{\w,x},\beta_{\w,x})\times\{(\w,x)\})$, then
\begin{equation}\label{def:2tau}
 \tau\colon\mV\subseteq\R\times\WR\to\WR\,,\;\;(t,\w,x)\mapsto (\wt,v(t,\w,x))
\end{equation}
defines a (possibly local) continuous flow on $\WR$, of {\em skewproduct} type.
As we will see in Subsection~\ref{subsec:hull}, families of this type appear in a natural way
when we construct the hull of a single equation.
\par
A {\em $\tau$-equilibrium} is a map $\mb\colon\W\to\R$ whose graph is  $\tau$-invariant
(i.e., with $v(t,\w,\mb(\w))=\mb(\wt)$ for all $\w\in\W$ and $t\in\R$). If it is
continuous, then its compact graph, which we represent by $\{\mb\}$, is a {\em
$\tau$-copy of the base} or {\em $\tau$-copy of $\W$}.
A $\tau$-copy of the base $\{\mb\}$ is {\em hyperbolic attractive} if there exists
$\delta>0$, $\gamma>0$ and $k\ge 1$ such that, if $|\mb(\w)-x|<\delta$ for any $\w\in\W$,
then $v(t,\w,x)$ is defined for all $t\ge 0$ and
$|\mb(\wt)-v(t,\w,x)|\leq k\,e^{-\gamma\,t}\,|\mb(\w)-x|$ for $t\ge 0$.
Changing $t\ge 0$ by $t\le 0$ provides the definition of
{\em repulsive hyperbolic} $\tau$-copy of the base. We will usually call $\tilde\mb:=\mb$ if
$\{\mb\}$ is a hyperbolic $\tau$-copy of the base, and we will often omit the prefix $\tau$.
\par
Given a bounded $\tau$-invariant set $\mB\subset\WR$ projecting onto $\W$, the maps
$\w\mapsto\inf\{x\in\R\,|\;(\w,x)\in\mB\}$ and $\w\mapsto\sup\{x\in\R\,|\;(\w,x)\in\mB\}$
define $\tau$-equilibria. We will refer to these maps as
the {\em lower\/} and {\em upper equilibria of $\mB$}. Observe that, if $\mB$ is
compact, then they are lower and upper semicontinuous, respectively,
and hence $m$-measurable for all $m\in\merg$.
\par
We complete this subsection with two definitions more.
If there exists a compact $\tau$-invariant set $\mA\subset\WR$ such that
$\lim_{t\to\infty} \text{dist}(\mC{\cdot}t,\mA)=0$ for every bounded set $\mC$,
where $\mC{\cdot}t=\{(\wt,v(t,\w,x))\,|\;(\w,x)\in\mC\}$ and
\[
 \text{dist}(\mC_1,\mC_2)=\sup_{(\w_1,x_1)\in\mC_1}\left(\inf_{(\w_2,x_2)\in\mC_2}
 \big(\mathrm{dist}_{\WR}((\w_1,x_1),(\w_2,x_2))\big)\right),
\]
then $\mA$ is the {\em global attractor for $\tau$}. And given two compact
subsets $\mK_1$ and $\mK_2$ of $\W\times\R$, we say that they are {\em ordered\/} with
$\mK_1<\mK_2$ if $x_1<x_2$ whenever there exists $\w\in\W$ such that
$(\w,x_1)\in\mK_1$ and $(\w,x_2)\in\mK_2$.
\subsection{Admissible processes and their hull extensions}\label{subsec:hull}
Let $\mU\subseteq\R^n$ be an open set.
We say that a continuous map $h\colon\R\times\mU\to\R$ is {\em admissible\/} and
write $h\in C^{0,0}(\R\times\mU,\R)$ if the restriction of $h$ to $\R\times\mJ$ is bounded
and uniformly continuous for any compact set $\mJ\subset\mU$. In most of the cases,
we will work with $\mU=\R$. We say that
$h\colon\RR\to\R$ is {\em $C^1$-admissible\/} (resp.~{\em $C^2$-admissible\/}) and write
$h\in C^{0,1}(\RR,\R)$ if there exists its derivative $h_x$ with respect to
the second variable and it is admissible (resp.~there exist $h_x$ and $h_{xx}$ and they are admissible).
\par
Given $h\in C^{0,1}(\RR,\R)$, we represent by $x_h(t,s,x)$ the maximal solution of
\begin{equation}\label{eq:2proceso}
 x'=h(t,x)
\end{equation}
with $x_h(s,s,x)=x$. By uniqueness of solutions,
$x_h(t,s,x_h(s,r,x))=x_h(t,r,x)$. Often, the map
$(t,s,x)\mapsto x_h(t,s,x)$ is called a {\em process}.
\par
We say that two solutions $b_1(t)$ and
$b_2(t)$ of \eqref{eq:2proceso} are {\em uniformly separated\/} if
they are bounded and $\inf_{t\in\R}|b_1(t)-b_2(t)|>0$. A bounded solution
$\tilde b(t)$ of \eqref{eq:2proceso} is {\em hyperbolic attractive\/}
(resp.~{\em hyperbolic repulsive}) if
there exist $k\ge 1$ and $\gamma>0$ such that
$\exp\Big(\int_s^t h_x(r,\tilde b(r))\,dr\Big)\le ke^{-\gamma(t-s)}$ whenever $t\ge s$
(resp.~$\exp\Big(\int_s^t h_x(r,\tilde b(r))\,dr\Big)\leq ke^{\gamma (t-s)}$ whenever $t\le s$);
and, in both cases, $(k,\gamma)$ is a {\em dichotomy constant pair of\/} $\tilde b$.
For the reader's convenience, we state the next fundamental result. A partial proof, strongly
based on \cite[Lecture 3]{copp}, can be found in \cite[Theorem 2.2]{dno3};
and a proof of the last assertion, strongly based on \cite[Theorem~III.2.4]{hale},
can be found in \cite[Theorem 3.2.3]{tfmduenas}
(just use the admissibility of $h_x$ instead of the existence and
boundedness of the second derivative).
We denote $\n{h}_{1,\rho}:=\sup_{(t,x)\in\R\times[-\rho,\rho]} |h(t,x)|+
\sup_{(t,x)\in\R\times[-\rho,\rho]} |h_x(t,x)|$ and $\n{b}_\infty:=\sup_{t\in\R}|b(t)|$.
\begin{teor}\label{th:2persistencia}
Let $h$ be $C^1$-admissible, let $\tilde b_h$ be an attractive (resp. repulsive) hyperbolic solution
of \eqref{eq:2proceso} with dichotomy constant pair $(k_0,\gamma_0)$, and take $\rho>\sup_{t\in\R}|\tilde b_h(t)|$.
Then, for every $\gamma\in(0,\gamma_0)$ and $\ep>0$, there exists $\delta_\ep>0$ and $\rho_\ep>0$ such that,
if $g$ is $C^1$-admissible and $\n{h-g}_{1,\rho}<\delta_\ep$, then
\begin{itemize}
\item[\rm(i)] there exists an attractive (resp. repulsive) hyperbolic solution $\tilde b_g$ of $x'=g(t,x)$
with dichotomy constant pair $(k_0,\gamma)$ which satisfies
$\big\|\tilde b_h-\tilde b_g\big\|_\infty<\ep$;
\item[\rm(ii)] if $|\tilde b_g(t_0)-x_0|\le\rho_\ep$, then $|\tilde b_g(t)-x_g(t,t_0,x_0)|
\le k_0\,e^{-\gamma(t-t_0)}|\tilde b_g(t_0)-x_0|$ for all $t\ge t_0$
(resp.~$|\tilde b_g(t)-x_g(t,t_0,x_0)|\le k_0\,e^{\gamma(t-t_0)}|\tilde b_g(t_0)-x_0|$
for all $t\le t_0$).
\end{itemize}
\end{teor}
A solution $\bar b\colon(-\infty,\beta)\to\R$ of \eqref{eq:2proceso}
is \emph{locally pullback attractive} if there exists $s_0<\beta$ and $\delta>0$
such that, if $s\leq s_0$, then $x_h(t,s,\bar b(s)\pm\delta)$ exists for $t\in[s,s_0]$,~and
\[
 \lim_{s\to-\infty}|\bar b(t)-x_h(t,s,\bar b(s)\pm\delta)|=0\qquad\text{for all }t\leq s_0\,.
\]
Analogously, a solution $\bar b\colon(\alpha,\infty)\to\R$ of \eqref{eq:2proceso}
is said to be \emph{locally pullback repulsive} if and only if there exist $s_0>\alpha$ and $\delta>0$
such that, if $s\geq s_0$ and $|x-\bar b(s)|<\delta$, then $x_h(t,s,x)$ exists for
$t\in[s_0,s]$ and
\[
 \lim_{s\to \infty}|\bar b(t)-x_h(t,s,\bar b(s)\pm\delta)|=0\qquad\text{for all }t\geq s_0\,.
\]
\par
Let us describe the already mentioned {\em hull construction}.
Given an admissible function $h\colon\RR\to\R$, we define $h{\cdot}t(s,x):=h(t+s,x)$.
The {\em hull\/} $\W_h$ of $h$ is the closure of the set
$\{h{\cdot}t\,|\;t\in\R\}$ on the set $C(\RR,\R)$ provided with the compact-open
topology. The set $\W_h$ is a compact metric space,
the time-shift map $\sigma_h\colon\R\times\W_h\to\W_h,(t,\w)\mapsto\w{\cdot}t$
defines a global continuous flow, and the map $\mh$ given by $\mh(\w,x)=\w(0,x)$
is continuous on $\W_h\times\R$. In addition, if $h$ is $C^1$-admissible
then $\W_h\subset C^{0,1}(\RR,\R)$,
and the continuous map $\mh_x(\w,x):=\w_x(0,x)$ is the derivative of
$\mh$ with respect to $x$; and, if $h$ is $C^2$-admissible
then $\W_h\subset C^{0,2}(\RR,\R)$,
and the continuous map $\mh_{xx}(\w,x):=\w_{xx}(0,x)$ is the
second derivative of $\mh$ with respect to $x$.
The proof of these properties can be found in \cite[Theorem~I.3.1]{shenyi} and
\cite[Theorem IV.3]{selltopdyn}.
Note that $(\W_h,\sigma_h)$ is a {\em transitive flow},
i.e., there exists a dense $\sigma_h$-orbit: that of the point $h\in\W_h$.
More precisely, if $\W_h^\upalpha$ and $\W_h^\upomega$ are the \upalfa-limit set
and \upomeg-limit set of the element $h\in\W_h$, then
\begin{lema}\label{lema:2hull}
$\W_h=\W_h^\upalpha\cup\{h{\cdot}t\,|\;t\in\R\}\cup\W_h^\upomega$.
\end{lema}
\begin{proof}
We can write any $\w\in\W_h$ as $\w=\lim_{n\to\infty} h{\cdot}t_n$
in the compact-open topology for a suitable sequence $(t_n)$.
If a subsequence $(t_k)$ has limit $-\infty$ or $+\infty$,
then $\w$ belongs to $\W_h^\upalpha$ or $\W_h^\upomega$. Otherwise,
there exists a subsequence $(t_k)$ with limit $t_0\in\R$, and it is easy to check that
$\w=h{\cdot}t_0$.
\end{proof}
The map $h$ is {\em recurrent} if $(\W_h,\sigma_h)$ is a {\em minimal flow}, i.e., if every
$\sigma_h$-orbit is dense in~$\W_h$.
\par
Assume that $h$ is (at least) $C^1$-admissible, and let us call $\tau_h$ the
skewproduct flow defined on $\W_h\times\R$ by the family of equations \eqref{eq:2fam}
corresponding to the constructed function $\mh$. Note that this family includes
\eqref{eq:2proceso}: it is given by the element $\w=h\in\W_h$, and, if
$\tau_h(t,\w,x)=(\wt,v_h(t,\w,x))$, then
$x_h(t,s,x)=v_h(t-s,h{\cdot}s,x)$. This is the {\em skew-product flow
induced by $h$ on its hull}. The next basic result will be used in Sections
\ref{sec:4} and \ref{sec:6}.
\begin{prop}\label{prop:2solhull}
Let $h$ be $C^1$-admissible and let $\W_h$ be its hull. If $x'=h(t,x)$
has a bounded solution $b$ (resp.~$n$ uniformly separated solutions $b_1<b_2<\ldots<b_n$),
then $x'=\w(t,x)$ has a bounded solution (resp.~$n$ uniformly separated solutions)
for all $\w\in\W_h$.
\end{prop}
\begin{proof}
We write $\w=\lim_{n\to\infty}h{\cdot}t_n$ in the compact-open topology for a
sequence $(t_n)$. Let $x_0$ be the limit of a suitable subsequence
$(b(t_k))$ of $(b(t_n))$. Then, the
solution $v_h(t,\w,x_0)$ of $x'=\w(t,x)$ (with value $x_0$ at $0$) is bounded, since
$v_h(t,\w,x_0)=\lim_{k\to\infty}b(t+t_k)$. The same
argument proves the other assertion.
\end{proof}
In what follows, we will consider both processes and skewproduct flows.
Observe that $(t,x)\mapsto \mh(\wt,x)$ is $C^1$-admissible for all
$\w\in\W$ if $\mh\in C^{0,1}(\WR,\R)$.
\begin{prop}\label{prop:2haciaatras}
\begin{itemize}
\item[\rm(i)] Let $\tilde b(t)$ be an attractive (resp.~repulsive) hyperbolic
solution of \eqref{eq:2proceso} for $h\in C^{0,1}(\RR,\R)$.
Then, $\inf_{t<t_0}|\tilde b(t)-\bar x(t)|>0$
(resp.~$\inf_{t>t_0}|\tilde b(t)-\bar x(t)|>0$) for any $t_0\in\R$ and any
solution $\bar x(t)\ne \tilde b(t)$ defined on $(-\infty, t_0]$
(resp~on $[t_0,\infty)$).
\item[\rm(ii)] Let $\tilde b_{\bar\w}(t)=v(t,\bar\w,b_0)$ be an attractive
(resp.~repulsive) hyperbolic solution of the equation corresponding to $\bar\w$
of the family \eqref{eq:2fam} given by a function $\mh\in C^{0,1}(\WR,\R)$,
and assume that the \upalfa-limit set (resp.~\upomeg-limit set) of
$(\bar\w,b_0)$ is an attractive (resp.~repulsive) hyperbolic copy of the
\upalfa-limit set $\W^\upalpha_{\bar\w}$ (resp.~\upomeg-limit set $\W^\upomega_{\bar\w}$) of $\bar\w$, say
$\{\tilde \mb\}$. Then, $\{\tilde\mb\}$ does not intersect the \upalfa-limit set
(resp.~\upomeg-limit set) of any $(\bar\w,x)$ with $x\ne b_0$ and bounded backward
semiorbit (resp.~bounded forward semiorbit).
\end{itemize}
\end{prop}
\begin{proof}
(i) We reason in the attractive case, assuming for contradiction the existence of
$(t_n)\downarrow-\infty$ such that $\lim_{n\to\infty}|\tilde b(t_n)-\bar x(t_n)|=0$.
According to the First Approximation Theorem (see \cite[Theorem~III.2.4]{hale} and
\cite[Proposition 2.1]{dno3}), the
attractive hyperbolicity of $\tilde b$ provides $k\ge 1$ and $\gamma>0$ such that,
for large enough~$n$,
\[
\begin{split}
 |\tilde b(t_0)-\bar x(t_0)|=|x_h(t_0,t_n,\tilde b(t_n))-x_h(t_0,t_n,\bar x(t_n))|
 \le k\,e^{-\gamma(t_0-t_n)}|\tilde b(t_n)-\bar x(t_n)|\,.
\end{split}
\]
The contradiction follows, since the last term tends to $0$ as $n\to\infty$.
\smallskip\par
(ii) We reason in the attractive case. Assume the existence of the \upalfa-limit set
$\mK$ of a point $(\bar\w,x)$ with $x\ne b_0$, and, for contradiction,
the existence of $(\w,\tilde\mb(\w))\in\mK$.
We write $(\w,\tilde\mb(\w))=\lim_{n\to\infty}(\bar\w{\cdot}t_n,v(t_n,\bar\w,x))$ for a
suitable sequence $(t_n)\downarrow-\infty$, assume without restriction the
existence of $\lim_{n\to\infty} v(t_n,\bar\w,b_0)$, observe that this limit is
also $\tilde\mb(\w)$, and note that this contradicts (i).
\end{proof}
\begin{prop}\label{prop:2extiende}
Let $\mh\in C^{0,1}(\WR,\R)$, and let $\tilde\mb\colon\W\to\R$ determine an attractive
(resp.~repulsive) copy of the base for \eqref{eq:2fam}. For any $\w\in\W$, the function $\tilde b_\w$
defined by $\tilde b_\w(t):=\tilde\mb(\wt)$ is an attractive (resp.~repulsive)
hyperbolic solution of \eqref{eq:2fam}$_\w$.
\end{prop}
\begin{proof}
Let us reason in the attractive case, fixing $\w\in\W$. Let us define $\w^*(t,x):=\mh(\wt,x)$ and
$v$ as in \eqref{def:2tau}.
Then, the solution $x_\w(t,s,x)$ of $x'=\w^*(t,x)$ (i.e., of \eqref{eq:2fam}$_\w$),
coincides with $v(t-s,\ws,x)$, and $\tilde b_\w(t)=v(t-s,\ws,\tilde b_\w(s))$.
Hence, the hyperbolicity of $\tilde\mb$ ensures that, if
$|x-\tilde b_\w(s)|\le\delta$ for an $s\in\R$, then $x_\w(t,s,x)$ exists
for all $t\ge s$ and it satisfies $|x_\w(t,s,x)-\tilde b_\w(t)|\le k e^{-\gamma(t-s)}|x-\tilde b_\w(s)|$.
Therefore, $(\partial/\partial x)x_\w(t,s,x)|_{x=\tilde b_\w(s)}
=\lim_{\ep\to 0} \big(x_\w(t,s,\tilde b_\w(s)+\ep)-x_\w(t,s,\tilde b_\w(s))\big)/\ep
\le k e^{-\gamma(t-s)}$. This derivative
solves the variational equation $z'=(\w^*)_x(t,\tilde b_\w(t))\,z$
and has value 1 at $t=s$, from where the assertion follows.
\end{proof}
\subsection{Lyapunov exponents}\label{subsec:Lyap}
Let $\mK\subset\WR$ be a $\tau$-invariant compact set projecting onto $\W$,
and let $\mathfrak{M}_\mathrm{erg}(\mK,\tau)$ be the
(nonempty) set of the $\tau$-ergodic measures on $\mK$.
The \emph{Lyapunov exponent of $\mK$ with respect to
$\nu\in\mathfrak{M}_\mathrm{erg}(\mK,\tau)$} is
$\gamma_{\mh_x}(\mK,\nu):=\int_\mK \mh_x(\w,x)\, d\nu$.
According to \cite[Theorem 4.1]{furstenberg1} and \cite[Theorem~1.8.4]{arnol},
for any $\nu\in\mathfrak{M}_\mathrm{erg}(\mK,\tau)$ projecting onto $m\in \merg$,
there exists an $m$-measurable $\tau$-equilibrium $\mb\colon\W\to\R$
such that
\begin{equation}\label{eq:21lyapunovexponents}
 \gamma_{\mh_x}(\mK,\nu)=\int_\W \mh_x(\w,\mb(\w))\, dm\,.
\end{equation}
Conversely, if $m\in\merg$
and $\mb\colon\W\to\R$ is a $m$-measurable $\tau$-equilibrium with graph in $\mK$, then
$\int_\W \mh_x(\w,\mb(\w))\,dm$ defines the Lyapunov exponent on $\mK$ with respect to
a suitable $\tau$-ergodic measure on $\mK$.
\begin{teor}\label{th:2copia}
Let $\mK\subset\WR$ be a $\tau$-invariant compact set projecting onto $\W$.
Assume that its upper and lower equilibria
coincide (at least) on a point of each minimal subset $\mM\subseteq\W$.
Then, all the Lyapunov exponents of $\mK$
are strictly negative (resp.~positive) if and only if $\mK$ is
an attractive (resp. repulsive) hyperbolic copy of the base.
\par
In addition, if either $\mK$ (and hence $\W$) is minimal or its upper and lower equilibria
coincide on a $\tau$-invariant subset $\W_0\subseteq\W$
with $m(\W_0)=1$ for all $m\in\merg$, then the condition on its upper and lower equilibria holds.
\end{teor}
\begin{proof}
We reason in the attractive case.
Let $\ml$ and $\muk$ be the lower and upper equilibria of $\mK$. We take
$(\w,\ml(\w))\in\mK$, a point $\w_0$ in a minimal subset of the \upalfa-limit set of $\w_0$
with $\ml(\w_0)=\muk(\w_0)$, and a sequence $(t_n)\downarrow-\infty$ with
$(\w_0,\ml(\w_0))=\lim_{n\to\infty}\tau(t_n,\w,\ml(\w))$ and such that there exists
$(\w_0,x_0):=\lim_{n\to\infty}\tau(t_n,\w,\muk(\w))$. Then, $\ml(\w_0)\le x_0\le \muk(\w_0)=\ml(\w_0)$;
i.e., $x_0=\ml(\w_0)$.
This property allows us to repeat the proof \cite[Proposition~2.8]{cano} in the
attractive case:
just replace the points $(\w_1,x_1)$ and $(\w_1,x_2)$ of that proof by $(\w,\ml(\w))$
and $(\w,\muk(\w))$. For the repulsive case, we work with the \upomeg-limit set.
\par
The last assertion is proved in \cite[Subsection 2.4]{cano} in the minimal case.
In the other one, it follows from the existence of a measure $m\in\merg$ concentrated
on each minimal set.
\end{proof}
\section{The concave and nonlinear case}\label{sec:3}
The main purpose in this section is to extend previous results on
families of equations on which hypotheses on concavity were globally assumed
(as in \cite{alob3,nuos4,lno2,lnor,dolo1}) to
a significantly less restrictive setting, on which the concavity
hypotheses are assumed just in measure.
\par
Let $(\W,\sigma)$ be a global continuous real flow on a compact metric space,
and let us consider the family of scalar ordinary differential equations
\begin{equation}\label{eq:3ini}
 x'=\mh(\wt,x)
\end{equation}
for $\w\in\W$, where $\mh\colon\WR\to\R$ satisfies (all or part of) the next conditions:
\begin{list}{}{\leftmargin 17pt}
\item[\hypertarget{c1}{{\bf c1}}] $\mh\in C^{0,1}(\WR,\R)$,
\item[\hypertarget{c2}{{\bf c2}}] $\limsup_{x\to\pm\infty}\mh(\w,x)<0$ uniformly on $\W$,
\item[\hypertarget{c3}{{\bf c3}}] $m(\{\w\in\W\,|\;x\mapsto\mh(\w,x) \text{ is concave}\})=1$ for all $m\in\merg$,
\item[\hypertarget{c4}{{\bf c4}}] $m(\{\w\in\W\,|\; x\mapsto\mh_x(\w,x)$ is strictly decreasing on $\mJ\})>0$
for all compact interval $\mJ\subset\R$ and all $m\in\merg$.
\end{list}
Recall that $\merg$ is the (nonempty) set of $\sigma$-ergodic measures on $\W$.
\begin{nota}\label{rm:3Ctambien}
Let $\W_0\subset\W$ be a nonempty compact $\sigma$-invariant subset. Then, any
$m_0\in{\mathfrak{M}_\mathrm{erg}(\W_0,\sigma)}$ can be extended to $m\in\merg$
by $m(\mU)=m_0(\mU\cap\W_0)$. So, if $\mh$ satisfies
\hyperlink{c1}{\bf cj} for {\bf j} $\in\{1,2,3,4\}$, also the restriction
$\mh\colon\W_0\times\R\to\R$ satisfies \hyperlink{c1}{\bf cj}.
\end{nota}
\par
Let $\tau$ be the skewproduct flow defined by \eqref{def:2tau}, with $\tau(t,\w,x)=(\wt,v(t,\w,x))$.
Before working under the coercive property \hyperlink{c2}{\bf c2},
we want to explain some consequences of the concavity assumptions \hyperlink{c3}{\bf c3}
and \hyperlink{c4}{\bf c4}, fundamental in what follows. They are based on the
next theorem, which establishes conditions under which two different ordered
$m$-measurable $\tau$-equilibria give rise to two Lyapunov exponents
\eqref{eq:21lyapunovexponents} with different signs.
Its proof is based on that of \cite[Theorem~4.1]{dno1}.
\begin{teor} \label{th:3C2medidas}
Let $\mh\colon\WR\to\R$ satisfy \hyperlink{c1}{\bf c1}, let us fix $m\in\merg$, and
let $\mb_1,\mb_2\colon\W\to\R$ be bounded $m$-measurable $\tau$-equilibria
with $\mb_1(\w)<\mb_2(\w)$ for $m$-a.e.~$\w\in\W$.
Assume that $m(\{\w\in\W\,|\;x\mapsto \mh(\w,x) \text{ is concave}\})=1$ and
$m(\{\w\in\W\,|\;\mh_{x}(\w,\mb_1(\w))>\mh_{x}(\w,\mb_2(\w))\})>0$. Then,
\[
 \int_\W \mh_x(\w,\mb_1(\w))\, dm>0\quad\text{and}\quad \int_\W \mh_x(\w,\mb_2(\w))\, dm<0\,.
\]
In particular, there are at most two bounded $m$-measurable $\tau$-equilibria which
are strictly ordered $m$-a.e.
\end{teor}
\begin{proof}
We call $\W^c:=\{\w\in\W\,|\;x\mapsto \mh(\w,x)\text{ is concave}\}$, which satisfies
$m(\W^c)=1$, and
$\W_0:=\{\w\in\W\,|\;\mb_1(\w)<\mb_2(\w)\}$, which is $\sigma$-invariant
(since $\mb_i(\wt)=v(t,\w,\mb_i(\w))$ for $i\in\{1,2\}$) and with $m(\W_0)=1$.
For each $\w\in\W^c$, we represent by $a_i(\w,x_1,x_2)$ the expression
$a_i(x_1,x_2)$ of \eqref{def:2ai} associated to the concave map
$x\mapsto \mh(\w,x)$ and observe that $(x_1,x_2)\mapsto a_i(\w,x_1,x_2)$ is continuous on
$\R^2$ for every $\w\in\W^c$: see \eqref{eq:2des}.
For $i\in\{1,2\}$, we define $a_i^*\colon\W\to\R$ by
$a_i^*(\w):=a_i(\w,\mb_1(\w),\mb_2(\w))$
if $\w\in\W^c\cap\W_0$ and $a_i^*(\w):=0$ otherwise,
and observe that $a_i^*$ is $m$-measurable and that $a_i^*\geq0$
(see Theorem~\ref{th:2positiveconcave}(i)). Let us take $i=1$ and write
\[
 \mh_x(\wt,\mb_1(\wt))-a_1^*(\wt)=\frac{\mh(\wt,\mb_2(\wt))
 -\mh(\w,\mb_1(\wt))}{\mb_2(\wt)-\mb_1(\wt)}=
 \frac{\mb_2'(\wt)-\mb_1'(\wt)}{\mb_2(\wt)-\mb_1(\wt)}
\]
for $\w\in\W_0$, where $\mb_i'(\wt)$ is the derivative of $t\mapsto\mb_i(\wt)$. This yields
\begin{equation}\label{eq:21intheorem23}
 \frac{1}{t}\int_0^t \mh_x(\w{\cdot}s,\mb_1(\w{\cdot}s))\,ds=\frac{1}{t}\int_0^t a_1^*(\w{\cdot}s)\, ds
 +\frac{1}{t}\log\left(\frac{\mb_2(\w{\cdot}t)-\mb_1(\w{\cdot}t))}{\mb_2(\w)-\mb_1(\w)}\right)\,.
\end{equation}
Lusin's Theorem provides a compact subset $\Delta\subset\W^c$ with
$m(\Delta)>0$ such that $\mb_1|_{\Delta},\mb_2|_{\Delta}\colon\Delta\to\R$ are continuous.
Since $\mh_x(\cdot,\mb_1(\cdot))$ is bounded and $a_1^*(\cdot)$ is nonnegative,
Birkhoff's Ergodic Theorem (see \cite[Theorem~1 in Section 1.2]{sinai1}
and \cite[Proposition~1.4]{johnson1}) ensures the existence of a $\sigma$-invariant subset
$\W_0^*\subseteq\W_0$ with $m(\W_0^*)=1$ such that, for every $\w\in\W_0^*$,
\[
\begin{split}
 \lim_{t\to\infty}\frac{1}{t}\int_0^{t} \mh_x\big(\w{\cdot}s,\mb_1(\w{\cdot}s)\big)\, ds
 &=\int_\W \mh_x\big(\w,\mb_1(\w)\big)\, dm\in\R\,,\\
 \lim_{t\to\infty} \frac{1}{t}\int_0^{t}a_1^*(\w{\cdot}s)\, ds&=\int_\W a_1^*(\w)\, dm\in [0,\infty]\,,\\
 \lim_{t\to\infty}\frac{1}{t}\int_0^t\chi_\Delta(\w{\cdot}s)\; ds&=m(\Delta)>0\,.
\end{split}
\]
In particular, if $\w\in\W_0^*$, there exists a sequence $(t_n)\uparrow\infty$
such that $\w{\cdot}t_n\in\Delta$. Hence, the sequence
$\{\log((\mb_2(\w{\cdot}t_n)-\mb_1(\w{\cdot}t_n))/(\mb_2(\w)-\mb_1(\w))\}_{n\in\N}$
is bounded. We write \eqref{eq:21intheorem23} for $t=t_n$
and take limit as $n\to\infty$ to get $\int_\W \mh_x(\w,\mb_1(\w))\, dm=\int_\W a_1^*(\w)\, dm$.
So, $\int_\W \mh_x(\w,\mb_1(\w))\, dm>0$ follows from $\int_\W a_1^*(\w)\, dm>0$. To
prove this last inequality, we deduce from Theorem~\ref{th:2positiveconcave}(i) that
$a_1^*(\w)>0$ if and only if $\mh_x(\w,\mb_1(\w))>\mh_x(\w,\mb_2(\w))$, and hence use
the last hypothesis on $\mh_x$ to get $m(\{\w\in\W^c\colon\, a_1^*(\w)>0\})>0$.
An analogous argument proves that $\int_\W \mh_x(\w,\mb_2(\w))\, dm<0$. The last assertion is
an easy consequence of the previous ones.
\end{proof}
\begin{teor}\label{th:3C2copias}
Let $\mh$ satisfy \hyperlink{c1}{\bf c1}, \hyperlink{c3}{\bf c3} and \hyperlink{c4}{\bf c4}.
Then, there exist two disjoint and ordered $\tau$-invariant compact sets $\mK_1<\mK_2$
projecting onto $\W$ if and only if there exist two different hyperbolic copies
of the base $\{\tmr\}$ and $\{\tma\}$ with $\tmr<\tma$. In this case, $\mK_1=\{\tmr\}$ and
it is repulsive; $\mK_2=\{\tma\}$ and it is attractive; and
$\mB:=\{(\w,x)\in\WR\,|\;\tmr(\w)\le x\le\tma(\w)\}$ is the set of globally
bounded orbits.
In particular, there are at most two disjoint and ordered $\tau$-invariant compact sets
projecting onto $\W$.
\end{teor}
\begin{proof}
Sufficiency is obvious.
To check necessity, we observe that $\mh$ satisfies the conditions of
Theorem \ref{th:3C2medidas} for any $m\in\merg$ and any pair of bounded ordered $m$-measurable equilibria.
This result ensures that all the Lyapunov exponents of $\mK_1$ are positive and
all the Lyapunov exponents of $\mK_2$ are negative.
The last assertion in Theorem \ref{th:3C2medidas} ensures that the upper and
lower equilibria of $\mK_i$ coincide on a $\tau$-invariant set with $m(\W_0)=1$
for all $m\in\merg$ for $i\in\{1,2\}$, and hence Theorem \ref{th:2copia} ensures that
$\mK_1$ is an attractive hyperbolic copy of $\W$ and
$\mK_2$ is a repulsive hyperbolic copy of $\W$.
Observe that this fact precludes the existence of more
than two disjoint and ordered $\tau$-invariant compact sets projecting onto $\W$.
\par
Let us write $\mK_1=\{\tmr\}$ and $\mK_2=\{\tma\}$ for continuous maps
$\tmr,\tma\colon\W\to\R$. Clearly,
$\bigcup_{\w\in\W}(\{\w\}\times[\tmr(\w),\tma(\w)])\subseteq\mB$. To prove the
converse inclusion, we assume for contradiction the existence of $(\w_0,x_0)$ with
$x_0>\tma(\w_0)$ and with globally defined and bounded $\tau$-orbit. Then,
the \upalfa-limit set $\mK$ of this orbit exists and is a $\tau$-invariant compact set
projecting onto a compact set $\W_\mK\subseteq\W$. Since $\{\tma\}$ is attractive,
Proposition \ref{prop:2haciaatras}(ii) restricted to $\W_\mK$ (see Remark~\ref{rm:3Ctambien})
ensures that $\mK>\mK_2|_{\W_\mK}>\mK_1|_{\W_\mK}$,
which contradicts the last assertion of the previous paragraph.
A similar argument working with \upomeg-limit sets shows that $x_0\ge\tmr(\w_0)$
for all $(\w_0,x_0)\in\mB$.
\end{proof}
We say that {\em there exists an attractor-repeller pair $(\tma,\tmr)$
of copies of the base\/} (or {\em of\/} $\W$)
{\em for \eqref{eq:3ini}}, or that {\em $(\tma,\tmr)$ is an
attractor-repeller pair of copies of the base for \eqref{eq:3ini}},
if $\{\tma\}$ is an attractive hyperbolic copy of $\W$ and $\{\tmr\}$ is a repulsive
hyperbolic copy of $\W$. So, Theorem \ref{th:3C2copias} characterizes its existence
under conditions \hyperlink{c1}{\bf c1}, \hyperlink{c3}{\bf c3} and \hyperlink{c4}{\bf c4},
in which case, in addition, $\tmr<\tma$.
\begin{nota}\label{rm:3basta}
Theorem~\ref{th:3C2copias} shows that, if there exists an attractor-repeller pair of copies of the base,
then the set $\mB$ of bounded $\tau$-orbits is nonempty, bounded, and the pair is given by its upper and lower
equilibria. Assume now that we previously know that $\mB$ is nonempty and bounded, with
$\mB\subseteq\W\times\mJ_\mB$ for a compact interval $\mJ_\mB\subset\R$.
Then all the conclusions of Theorem~\ref{th:3C2copias} apply if, for all $m\in\merg$,
$m(\{\w\in\W\,|\;x\mapsto \mh(\w,x) \text{ is concave}\})=1$ and
$m(\{\w\in\W\,|\; x\mapsto \mh_x(\w,x)$ is strictly decreasing on $\mJ_\mB\})>0$.
\end{nota}
\par
Let us now derive consequences of the coercivity property.
\begin{prop}\label{prop:3Ccoer}
Let $\mh$ satisfy \hyperlink{c1}{\bf c1} and \hyperlink{c2}{\bf c2}, and take
$\delta>0$ and $m_1<m_2$ with $\mh(\w,x)\le-\delta$ for all $\w\in\W$ if $x\notin(m_1,m_2)$.
Then,
\begin{itemize}
\item[\rm(i)] $\liminf_{t\to(\alpha_{\w,x})^+}v(t,\w,x)>m_1$ and
$\limsup_{t\to(\beta_{\w,x})^-}v(t,\w,x)<m_2$ for any solution $v(t,\w,x)$: any
solution remains lower bounded as time decreases and upper bounded as time increases.
\item[\rm(ii)] If $v(t,\w,x)$ is bounded, then $v(t,\w,x)\in[m_1,m_2]$ for all $t\in\R$: the set
\[
 \mB:=\Big\{(\w,x)\,|\;\sup_{t\in\R}|v(t,\w,x)|<\infty\Big\}
\]
is either empty or contained in $\W\times[m_1,m_2]$.
\item[\rm(iii)] If $\mB$ is nonempty, then the projection $\W^b$ of $\mB$ onto $\W$ is a $\sigma$-invariant compact set.
\item[\rm(iv)] For each $\w\in\W^b$, let us write $\mB_\w:=\{x\,|\;(\w,x)\in\mB\}
=[\mr(\w),\ma(\w)]$. Then, the maps $\mr,\ma\colon\W^b\to[m_1,m_2]$ are lower and upper semicontinuous
equilibria for the restriction of $\tau$ to $\W^b\times\R$.
\item[\rm(v)] If, for a point $\w\in\W$, there exists a bounded $C^1$ function $b\colon\R\to\R$
such that $b'(t)\le \mh(\wt,b(t))$ for all $t\in\R$, then
$\w\in\W^b$, and $\mr(\wt)\le b(t)\le \ma(\wt)$ for all $t\in\R$.
If $b'(t)<\mh(\wt,b(t))$ for all $t\in\R$, then
$\mr(\wt)<b(t)<\ma(\wt)$ for all $t\in\R$.
\item[\rm(vi)] If $\w\in\W^b$, then $v(t,\w,x)$ is bounded from below if and only if $x\ge \mr(\w)$, and
from above if and only if $x\le \ma(\w)$.
\item[\rm(vii)] Assume that $\mh$ satisfies also \hyperlink{c3}{\bf c3} and \hyperlink{c4}{\bf c4},
and that $(\tma,\tmr):=(\ma,\mr)$ is an attractor-repeller pair of copies of the base.
Then, $\lim_{t\to\infty}(v(t,\w,x)-\tma(\wt))=0$ if and only if $x>\tmr(\w)$,
$\lim_{t\to-\infty}(v(t,\w,x)-\tmr(\wt))=0$ if and only if $x<\tma(\w)$,
and $t\mapsto\tmr(\wt),\tma(\wt)$ define the two unique hyperbolic solutions of \eqref{eq:3ini}$_\w$.
\end{itemize}
\end{prop}
\begin{proof} The existence of $m_1$ and $m_2$ is ensured by property \hyperlink{c2}{\bf c2}.
The proofs of (i) and (ii) are classical exercises on ODEs. It is easy to check that $\mB$ is closed,
and hence compact, and clearly it is $\tau$-invariant. Assertions (iii) and (iv) follow from here.
The properties stated in (v) follow from (i) and standard comparison arguments: see e.g. the
proof of \cite[Theorem 3.1(v)]{lnor}.
Easy contradiction arguments using (i) prove (vi).
\par
To check the first assertion in (vii), we first take $x>\tmr(\w)$, assume for contradiction
that the  \upomeg-limit set of $(\w,x)$ is not contained in $\{\tma\}$, deduce from Theorem
\ref{th:3C2copias} (restricted to the projection of the \upomeg-limit set: see Remark \ref{rm:3basta})
that it intersects $\{\tmr\}$, and observe that this contradicts Proposition \ref{prop:2haciaatras}(ii).
Conversely, if $\lim_{t\to\infty}(v(t,\w,x)-\tma(\wt))=0$,
then (vi) and the $\tau$-invariance of $\{\tmr\}$ ensure that $x>\tmr(\w)$. The same arguments prove
the second assertion in (vii), and the last one
follows from Propositions \ref{prop:2extiende} and \ref{prop:2haciaatras}.
\end{proof}
Note that the previous property (iv) ensures that $\mr$ and $\ma$ are $m$-measurable equilibria
for all $m\in\mathfrak{M}_\mathrm{erg}(\W^b,\sigma)$. We will use this property when we apply
Theorem \ref{th:3C2medidas} to these equilibria.
\par
The last result in this section characterizes the existence of an attractor-repeller pair of
copies of the base in terms of the existence of two uniformly separated hyperbolic solutions of a
given equation when the base is constructed as the hull of that equation: see Subsection \ref{subsec:hull}.
\begin{teor}\label{th:3Cequiv}
Let $\mh\colon\WR\to\R$ satisfy \hyperlink{c1}{\bf c1}, \hyperlink{c2}{\bf c2},
\hyperlink{c3}{\bf c3} and \hyperlink{c4}{\bf c4}. Let us fix $\bar\w\in\W$.
Then, the following assertions are equivalent:
\begin{itemize}
\item[\rm(a)] Equation \eqref{eq:3ini}$_{\bar\w}$ has two hyperbolic solutions.
\item[\rm(b)] Equation \eqref{eq:3ini}$_{\bar\w}$ has two uniformly separated
hyperbolic solutions.
\item[\rm(c)] Equation \eqref{eq:3ini}$_{\bar\w}$ has two uniformly separated
bounded solutions.
\item[\rm(d)] There exists an attractor-repeller pair $(\tma,\tmr)$ of copies of the base
for the restriction of the family \eqref{eq:3ini} to the closure $\W_{\bar\w}$ of
$\{\bwt\,|\;t\in\R\}$.
\end{itemize}
In this case, $t\mapsto\tilde a(t):=\tma(\bwt)$ and $t\mapsto\tilde r(t):=\tmr(\bwt)$ are the
two unique uniformly separated solutions of \eqref{eq:3ini}$_{\bar\w}$, they are hyperbolic,
and there are no more hyperbolic solutions.
In addition, if $x_{\bar\w}(t,s,x)$ is the solution of \eqref{eq:3ini}$_{\bar\w}$ with
$x_{\bar\w}(s,s,x)=x$, then: it is bounded if and only if $x\in[\tilde r(s),\tilde a(s)]$,
$\lim_{t\to\infty}|x_{\bar\w}(t,s,x)-\tilde a(t)|=0$ if and only if $x>\tilde r(s)$, and
$\lim_{t\to-\infty}|x_{\bar\w}(t,s,x)-\tilde r(t)|=0$ if and only if $x<\tilde a(s)$.
\end{teor}
\begin{proof}
The statements after the equivalences follow from (d) and
Proposition \ref{prop:3Ccoer}(vi) and (vii).
We will check (b)$\,\Rightarrow\,$(c)$\,\Rightarrow\,$(d)$\,\Rightarrow\,$(a)$\,\Rightarrow\,$(b).
Recall that the hypotheses on $\mh$ are also valid for its restriction to $\W_{\bar\w}\times\R$:
see Remark \ref{rm:3Ctambien}.
\smallskip\par
(b)$\,\Rightarrow\,$(c)$\,\Rightarrow\,$(d). Obviously, (b) implies (c).
Now we assume (c) and observe that it ensures that the lower and upper bounded solutions,
$r(t)$ and $a(t)$, are uniformly separated. We call $\delta:=\inf_{t\in\R}(a(t)-r(t))>0$.
Let $\mK_a$ be the closure of the $\tau$-orbit of $(\bar\w,a(0))$. It projects on $\W_{\bar\w}$,
and hence $\W_{\bar\w}\subset\W^b$: there exists $\mr(\w)$ and $\ma(\w)$
for all $\w\in\W_{\bar\w}$.
Let us check that $x_0\ge\mr(\w_0)+\delta$ for all $(\w_0,x_0)\in\mK_a$.
We write $(\w_0,x_0)=\lim_{n\to\infty}(\bar\w{\cdot}t_n,a(t_n))$
and assume without restriction the existence of $(\w_0,x^0):=
\lim_{n\to\infty}(\bar\w{\cdot}t_n,r(t_n))$, which belongs to the (closed) set $\mB$.
If $x_0<\mr(\w_0)+\delta$, then $x^0\le x_0-\delta<\mr(\w_0)+\delta-\delta=\mr(\w_0)$, impossible.
Let us consider the restriction $\bar\tau$ of $\tau$ to $\W_{\bar\w}\times\R$.
Since any $\bar\tau$-equilibria with graph in $\mK_a$ is strictly above $\mr$,
Theorem~\ref{th:3C2medidas} shows that all the Lyapunov exponents of $\mK_a$ are
strictly negative, and that its upper and lower equilibria of $\mK_a$
coincide on a $\sigma$-invariant set
$\W_0$ with $m_0(\W_0)=1$ for all $m_0\in\mathfrak{M}_\mathrm{erg}(\W_{\bar\w},\sigma)$.
Hence, Theorem \ref{th:2copia}, ensures that $\mK_a$ is an attractive
hyperbolic copy of $\W_{\bar\w}$. This fact and the
previous property ensure that $\mK_a$ is strictly above the closure $\mK_r$ of
$\{(\w,\mr(\w))\,|\;\w\in\W_{\bar\w}\}$. Hence, Theorem \ref{th:3C2copias} ensures
that $\mK_r$ is a repulsive hyperbolic copy of $\W_{\bar\w}$:
(d) holds.
\smallskip\par
(d)$\,\Rightarrow\,$(a)$\,\Rightarrow\,$(b). If (d) holds, then
$t\mapsto\tma(\bwt)$ and $t\mapsto\tmr(\bwt)$ are two
hyperbolic solutions of \eqref{eq:3ini}$_{\bar\w}$
(see Proposition \ref{prop:2extiende}), which ensures (a).
Let us assume (a), and let $\tilde x_1<\tilde x_2$
be the two hyperbolic solutions of \eqref{eq:3ini}$_{\bar\w}$.
Let us first check that $\tilde x_1$ is
repulsive, assuming for contradiction that it is attractive. We call
$a(t):=\ma(\bwt)$, with $\ma$
given by Proposition \ref{prop:3Ccoer}(iv). Proposition \ref{prop:2haciaatras}(i)
ensures that $\delta:=\inf_{t\le 0}(a(t)-\tilde x_1(t))>0$.
Let $\mM_1$ and $\mM_a$ be the \upalfa-limit sets of $(\bar\w,\tilde x_1(0))$
and $(\bar\w,a(0))$, which project on the \upalfa-limit
set $\W_{\bar\w}^\upalpha\subseteq\W_{\bar\w}$ of $\bar\w$. Repeating the argument
of the previous paragraph, we check that
$x_0\le\ma(\w_0)-\delta$ whenever $(\w_0,x_0)\in\mM_1$,
and deduce that $\mM_1$ is a repulsive copy of $\W_{\bar\w}^\upalpha$, which
contradicts Proposition \ref{prop:2extiende}.
\par
Hence, $\tilde x_1$ is repulsive. Proposition \ref{prop:2haciaatras}(i)
yields $\delta>0$ such that $\inf_{t\ge 0}(a(t)-\tilde x_1(t))>\delta$. Let
$\bar\mM_1$ be the \upomeg-limit set of $(\bar\w,\tilde x_1(0))$,
which projects on the \upomeg-limit set $\W_{\bar\w}^\upomega\subseteq\W_{\bar\w}$ of $\bar\w$.
Repeating again the arguments used to prove (c)$\,\Rightarrow\,$(d), we check that
$x_0\le \ma(\w_0)-\delta$ whenever $(\w_0,x_0)\in\bar\mM_1$; and we deduce that
$\bar\mM_1$ is a repulsive copy of $\W_{\bar\w}^\upomega$. Hence,
$\bar\mM_1$ does not intersect the \upomeg-limit set $\bar\mM_2$ of $(\bar\w,\tilde x_2(0))$:
see Proposition \ref{prop:2haciaatras}(ii).
So, we have $\bar\mM_1<\bar\mM_2$. Theorem \ref{th:4DC3copias} applied to $\W_{\bar\w}^\upomega\times\R$
ensures that $\bar\mM_2$ is an attractive hyperbolic copy of $\W_{\bar\w}^\upomega$,
which according to Proposition \ref{prop:2extiende} is only possible
if $\tilde x_2$ is attractive. Proposition \ref{prop:2haciaatras}(i)
ensures that the two solutions are uniformly separated. So, (b) holds.
\end{proof}
We will refer to the situation described by the equivalences of Theorem
\ref{th:3Cequiv} as the {\em existence of an attractor-repeller pair of
solutions of\/} \eqref{eq:3ini}$_{\bar\w}$.
\section{Asymptotically concave transition equations}\label{sec:4}
Let $g\colon\RR\to\R$ be a $C^1$-admissible function.
The hull construction described
in Subsection \ref{subsec:hull} allows us to understand the $\sigma$-orbit
of $g$, $\{g{\cdot}t\,|\;t\in\R\}$, which is dense in the hull $\W_g$,
as a connection between its
\upalfa-limit set $\W_g^\upalpha$ and its \upomeg-limit set $\W_g^\upomega$.
In fact, the hull $\W_g$ is the union of these three sets: see Lemma \ref{lema:2hull}.
Our goal in this section is to describe the dynamical possibilities of
an ``asymptotically concave" equation
\begin{equation}\label{eq:4Cgtran}
 x'=g(t,x)
\end{equation}
under conditions which ensure that the families of equations
defined over $\W_g^\upalpha$ (\upalfa-family) and
$\W_g^\upomega$ (\upomeg-family) satisfy the regularity, coercivity
and strict concavity properties \hyperlink{c1}{\bf c1}-\hyperlink{c4}{\bf c4},
as well as the existence of attractor-repeller pairs
of copies of the base for the \upalfa-family and the \upomeg-family.
Since the structures of these sets represent
the past and future of $g$, we are understanding \eqref{eq:4Cgtran}
as a transition between the \upalfa-limit and \upomeg-limit families.
\par
Proposition \ref{prop:2solhull} precludes the existence of
uniformly separated solutions of
\eqref{eq:4Cgtran} unless all the equations of the \upalfa-family and
the \upomeg-family have uniformly separated solutions. Hence,
to consider a transition scenario with interesting dynamical
possibilities, it is reasonable to assume the
existence of attractor-repeller pairs of solutions for all the
(concave) equations of the limit families: see Theorem \ref{th:3Cequiv}.
We will achieve these properties by assuming the existence of strictly
concave (in $x$) maps $g_-$ and $g_+$ such that $g$ and $g_-$ (resp.~$g$ and $g_+$) form
an asymptotic pair as $t\to-\infty$ (resp.~as $t\to\infty$) in the common hull of
$g$ and $g_-$ (resp.~$g$ and $g_+$). (The common hull of two admissible maps
$h_1$ and $h_2$ is the compact metric space defined as the closure of
$\{h_i{\cdot}t\,|\;i=1,2,\,t\in\R\}$ in the compact-open topology,
and $\w_1$ and $\w_2$ form an asymptotic pair as $t\to\pm\infty$ if
the distance from $\w_1{\cdot}t$ to $\w_2{\cdot}t$ tends to 0.) The required
existence of these maps does not imply their uniqueness, but
Lemma \ref{lema:4Cunion} below also shows that $\W_g^\upalpha$ and
$\W_g^\upomega$ respectively coincide with $\W_{g_-}^\upalpha$ and $\W_{g_+}^\upomega$,
which is a key point in our analysis.
\par
So, we fix $g$ and assume the existence of $g_-$ and $g_+$ such that:
\begin{list}{}{\leftmargin 23pt}
\item[\hypertarget{gc1}{\bf gc1}] $g,g_-,g_+\in C^{0,1}(\RR,\R)$.
\item[\hypertarget{gc2}{\bf gc2}] $\lim_{t\to\pm\infty}(g(t,x)-g_\pm(t,x))=0$
uniformly on each compact subset $\mJ\subset\R$.
\item[\hypertarget{gc3}{\bf gc3}] $\limsup_{x\to\pm\infty} h(t,x)<0$
uniformly on $\R$ for $h=g,g_-,g_+$.
\item[\hypertarget{gc4}{\bf gc4}] $\inf_{t\in\R}\big((g_\pm)_x(t,x_1)-(g_\pm)_x(t,x_2)\big)>0$
whenever $x_1<x_2$.
\item[\hypertarget{gc5}{\bf gc5}] Each one of the equations
\begin{equation}\label{eq:4Cglim}
 x'=g_-(t,x) \quad\text{and}\quad x'=g_+(t,x)
\end{equation}
has two hyperbolic solutions, $\tilde r_{g_-}<\tilde a_{g_-}$ and $\tilde r_{g_+}<\tilde a_{g_+}$.
\end{list}
As Lemma \ref{lema:4Chiposc} will prove, conditions
\hyperlink{gc1}{\bf gc1}-\hyperlink{gc3}{\bf gc4} provide a setting satisfying
the hypotheses of Section \ref{sec:3}.
\begin{notas}\label{rm:4Chipos}
1.~Slightly abusing language,
we will say that ``$g$ satisfies conditions \hyperlink{gc1}{\bf gc1}-\hyperlink{gc5}{\bf gc5}"
if there exist $g_-$ and $g_+$ such that all the listed conditions are satisfied.
\par
2.~To simplify the language, we will refer to \eqref{eq:4Cgtran}
as a {\em transition equation\/} between the {\em past equation\/} and the {\em future equation},
which are the first one and the second one in \eqref{eq:4Cglim}. That the use of
these words is accurate is partly justified by the previously mentioned
equalities $\W_{g_-}^\upalpha=\W_g^\upalpha$ and $\W_{g_+}^\upomega=\W_g^\upomega$,
which mean that the hyperbolic structures of the equations
\eqref{eq:4Cglim} condition that of \eqref{eq:4Cgtran} and viceversa;
and it will be better justified by the main results of this section.
But observe that the future of the dynamics of the nonautonomous
equation $x'=g_-(t,x)$ is not necessarily related to its past (since
$\W_{g_-}^\upalpha$ can be different $\W_{g_-}^\upomega$), and hence it
can be not related to the dynamics of $x'=g(t,x)$. And the
same happens with the past dynamics of $x'=g_+(t,x)$ and $x'=g(t,x)$.
\end{notas}
\par
We will classify the dynamical scenarios for
\eqref{eq:4Cgtran} and relate them to those of \eqref{eq:4Cglim}
under the above conditions, which include coercivity of all the involved equations
(\hyperlink{gc3}{\bf gc3}) but strict concavity in $x$ only of the limit
ones (\hyperlink{gc4}{\bf gc4}).
As said in the Introduction, several fundamental differences arise with respect to
previous approaches, which, on the one hand, are restricted to maps $g(t,x):=f(t,x,\G(t))$ and $g_\pm(t,x):=
f(t,x,\gamma_\pm)$, where $\gamma_\pm:=\lim_{t\to\pm\infty}\G(t)$ are assumed to exist and be real;
and, on the other hand, are analyzed under much more exigent concavity
hypothesis. So, we will extend part of the results of \cite{lnor},
formulated for $x'=-(x-\G(t))^2+p(t)$ and with constant asymptotic limits of $\G$,
to a much more general setting. The problem is also analyzed in \cite{lno2}
for $x'=h(t,x-\G(t))$ assuming (less restrictive) Carath\'{e}odory conditions
on $h$ and properties concerning its concavity with respect to the second variable and the
asymptotic limits of $\G$ which are much stronger than those assumed here. Thus, the current
formulation of our results considerably broadens their possibilities of application,
as we will see in Subsection \ref{subsec:42}.
\par
Theorem \ref{th:4Ccasos} provides the above mentioned classification of the
dynamical possibilities for
\eqref{eq:4Cgtran} when all the previous conditions hold. Its proof is strongly based on
Theorems \ref{th:3Cequiv} and \ref{th:4Cdossoluciones}, and this last one also requires some previous
work. Our first two results, fundamental for the subsequent application of Theorem \ref{th:3Cequiv},
refer to the hull extensions (see Subsection \ref{subsec:hull}).
Recall that we represent by $x_h(t,s,x)$ the maximal solution of $x'=h(t,x)$ which
satisfies $x_h(s,s,x)=x$: we will use this notation for $h$ equal to $g$, $g_-$, $g_+$,
and some other auxiliary admissible functions. In all these cases, the set $\W_h$ is the hull of
$h$, and $\W_h^\upalpha$ and $\W_h^\upomega$ are the \upalfa-limit set and \upomeg-limit set
of the element $h\in\W_h$. Recall that $h{\cdot}t(s,x):=h(t+s,x)$, and that
$\mh(\w,x):=\w(0,x)$ if $\w\in\W_h$. We represent by $\mg,\,\mg_-$ and $\mg_+$
the extensions to the corresponding hulls of $g,\,g_-$ and $g_+$. These auxiliary results do not need all the conditions
\hyperlink{gc1}{\bf gc1}-\hyperlink{gc5}{\bf gc5}: we will specify the required ones.
\begin{lema}\label{lema:4Cunion}
Let $g$ and $g_\pm$ satisfy \hyperlink{gc1}{\bf gc1} and \hyperlink{gc2}{\bf gc2}.
Then, $\W_g^\upalpha=\W_{g_-}^\upalpha$ and $\W_g^\upomega=\W_{g_+}^\upomega$. Hence,
$\W_g=\W_{g_-}^\upalpha\cup\{g{\cdot}t\,|\;t\in\R\}\cup\W_{g_+}^\upomega$.
\end{lema}
\begin{proof}
Given a sequence $(t_n)$ with limit $\pm\infty$, it is easy to check that
$\w(t,x)=\lim_{n\to\infty}g(t+t_n,x)$ uniformly on the
compact subsets of $\RR$ if and only if $\w(t,x)=\lim_{n\to\infty}g_\pm(t+t_n,x)$ uniformly on the
compact subsets of $\RR$. This proves the first equalities, which combined with Lemma \ref{lema:2hull} prove the last one.
\end{proof}
\begin{lema}\label{lema:4Chiposc}
If $h\in C^{0,1}(\RR,\R)$ then $\mh$ satisfies \hyperlink{c1}{\bf c1} on $\W_h$.
If $h\in C^{0,1}(\RR,\R)$ and $\limsup_{x\to\pm\infty} h(t,x)<0$
uniformly on $\R$, then $\mh$ satisfies \hyperlink{c2}{\bf c2} on $\W_h$.
And, if \hyperlink{gc1}{\bf gc1}, \hyperlink{gc2}{\bf gc2} and \hyperlink{gc4}{\bf gc4} hold, then
$\mg$ and $\mg_\pm$ satisfy \hyperlink{c3}{\bf c3} and \hyperlink{c4}{\bf c4} on $\W_g$ and $\W_{g_\pm}$,
respectively.
\end{lema}
\begin{proof}
As explained in Subsection \ref{subsec:hull}, \hyperlink{c1}{\bf c1} for $\mh$ follows from
the $C^1$-admissibility of $h$. If, in addition, $\limsup_{x\to\pm\infty} h(t,x)<0$ uniformly on $\R$,
then there exists $\delta>0$ and $\rho_\delta>0$ such
that $h(t,x)\le-\delta$ if $|x|\ge\rho_\delta$ and $t\in\R$. Since any $\w\in\W_h$ satisfies
$\w(0,x)=\lim_{n\to\infty}h(t_n,x)$ for a sequence $(t_n)$, we have
$\mh(\w,x)=\w(0,x)\le-\delta$ if $|x|\ge\rho_\delta$: \hyperlink{c2}{\bf c2} holds on $\W_h$.
\par
Now, we assume that \hyperlink{gc1}{\bf gc1}, \hyperlink{gc2}{\bf gc2} and \hyperlink{gc4}{\bf gc4} hold.
To prove the last assertion, it is enough to reason with $h=g$, since $g_-$ and $g_+$
satisfy the conditions assumed on $g$. To check that
$\mg$ satisfies \hyperlink{c3}{\bf c3} and \hyperlink{c4}{\bf c4} on $\W_g$,
we use the ideas of the proof of \cite[Proposition 3.16]{dno3}.
Lemma \ref{lema:2hull} ensures that
$\W_g=\W_g^\upalpha\cup\{g{\cdot}t\,|\;t\in\R\}\cup\W_g^\upomega$.
In particular, given $m\in\mathfrak{M}_\mathrm{erg}(\W_g,\sigma_g)$,
either $m(\W_g^\upalpha)=1$ or $m(\W_g^\upomega)=1$ (or both):
this is trivial if $g$ is independent of $t$ or $t$-periodic (since $\W_g=\W_g^\upalpha=\W_g^\upomega$);
and, in the remaining cases,
$\{g{\cdot}t\,|\;t\in\R\}=\bigcup_{n\in\Z} \sigma_n(\{g{\cdot}t\,|\;t\in[0,1)\})$ (where $\sigma_n(\w)=\w{\cdot}n$)
is a nonfinite union of disjoint sets. Therefore, $m(\sigma_n(\{g{\cdot}t\,|\;t\in[0,1)\}))=0$
for all $n\in\N$, since this measure is independent of $n$. Hence, it suffices to
check that $x\mapsto\mg_x(\w,x)$ is strictly decreasing on $\R$ for all $\w\in\W_g^\upalpha\cup\W_g^\upomega$:
this ensures that $m(\{\w\in\W_g\,|\;x\mapsto \mg_x(\w,x)$ is strictly decreasing on $\R\})=1$
for all $m\in\mathfrak{M}_\mathrm{erg}(\W_g,\sigma_g)$,
which is stronger than \hyperlink{c3}{\bf c3} and \hyperlink{c4}{\bf c4}.
We reason for $\w\in\W_g^\upomega$. According to Lemma \ref{lema:4Cunion},
$\w=\lim_{n\to\infty}g_+{\cdot}t_n$ (in the compact-open topology)
for a sequence $(t_n)$ with limit $\infty$. Then,
$\w_x$ is the limit of any subsequence of $((g_+)_x{\cdot}t_n)$
which uniformly converges on the compact subsets of $\RR$, and hence
$\w_x=\lim_{n\to\infty}(g_+)_x{\cdot}t_n$ uniformly on the compact subsets of $\RR$. We take $x_1<x_2$, and apply \hyperlink{gc4}{\bf gc4} to get
\[
\begin{split}
 \mg_x(\w,x_1)-\mg_x(\w,x_2)&=\w_x(0,x_1)-\w_x(0,x_2)\\
 &=\lim_{n\to\infty}((g_+)_x(t_n,x_1)-(g_+)_x(t_n,x_2))
 >0\,,
\end{split}
\]
which completes the proof.
\end{proof}
\begin{nota}\label{rm:4Creformular}
Lemma \ref{lema:4Chiposc} shows that
$\mg_-$ satisfies \hyperlink{c1}{\bf c1}, \hyperlink{c2}{\bf c2},
\hyperlink{c3}{\bf c3} and \hyperlink{c4}{\bf c4} if $g_-$
satisfies the conditions assumed on it on \hyperlink{gc1}{\bf gc1}, \hyperlink{gc3}{\bf gc3} and
\hyperlink{gc4}{\bf gc4}. Hence, in this case, and
according to Theorem \ref{th:3Cequiv}, the property corresponding to $g_-$ in
condition \hyperlink{gc5}{\bf gc5} can be reformulated as:
``the equation $x'=g_-(t,x)$ has an attractor-repeller pair
of solutions $(\tilde a_{g_-},\tilde r_{g_-})$", which determines its
corresponding global dynamics: see Theorem \ref{th:3Cequiv}.
The same applies to $g_+$. We will use these facts without further reference.
\end{nota}
The next result allows us to apply Theorem \ref{th:2persistencia} in the
proof of Theorem \ref{th:4Cdossoluciones}.
\begin{lema}\label{lem:4Cderivadas}
If $g$ and $g_\pm$ satisfy \hyperlink{gc1}{\bf gc1} and \hyperlink{gc2}{\bf gc2}, then
$\lim_{t\to\pm\infty}(g_x(t,x)-(g_\pm)_x(t,x))=0$ uniformly on each compact subset $\mJ\subset\R$.
\end{lema}
\begin{proof}
Let us reason for $g_-$, taking $(t_n)\downarrow-\infty$ and a compact subset $\mJ\subset\R$.
Since $h_-:=g-g_-$ is $C^1$-admissible, every subsequence
$(t_m)$ has a subsequence $(t_k)$ such that there exists
$d_-(x):=\lim_{k\to\infty}(h_-)_x(t_k,x)$ and is uniform on $\mJ$.
We assume for contradiction that $d_-\not\equiv 0$ and, without
restriction, that $d_-(x)\ge\ep>0$ for $x\in[x_1,x_2]\subseteq\mJ$. Then,
$0=\lim_{k\to\infty}(h_-(t_k,x_2)-h_-(t_k,x_1))=\lim_{k\to\infty}
\int_{x_1}^{x_2}(h_-)_x(t_k,s)\,ds=
\int_{x_1}^{x_2}d_-(s)\,ds\ge\ep(x_2-x_1)$, which is impossible.
\end{proof}
\begin{teor}\label{th:4Cdossoluciones}
Assume that $g$ satisfies \hyperlink{gc1}{\bf gc1}-\hyperlink{gc5}{\bf gc5}, and let
$(\tilde a_{g_\pm},\tilde r_{g_{\pm}})$ be the attractor-repeller pairs of solutions of
$x'=g_\pm(t,x)$ given by \hyperlink{gc5}{\bf gc5}. Then,
\begin{itemize}
\item[\rm(i)] there exists a unique solution $a_g$ of \eqref{eq:4Cgtran} defined
at least on a negative half-line and characterized by
``$x>a_g(s)$ if and only if $x_g(t,s,x)$ is unbounded from above as time decreases".
It satisfies $\lim_{t\to-\infty}(a_g(t)-\tilde a_{g_-}(t))=0$ and, if there exists
$a_g(s)$, then $x<a_g(s)$ if and only if $\lim_{t\to-\infty}|x_g(t,s,x)-\tilde r_{g_-}(t)|=0$.
In addition, $a_g$ is locally pullback attractive.
\item[\rm(ii)] There exists a unique solution $r_g$ of \eqref{eq:4Cgtran} defined
at least on a positive half-line and characterized by
``$x<r_g(s)$ if and only if $x_g(t,s,x)$ is unbounded from below as time increases".
It satisfies $\lim_{t\to\infty}(r_g(t)-\tilde r_{g_+}(t))=0$ and, if there exists $r_g(s)$,
then $x>r_g(s)$ if and only if $\lim_{t\to\infty}|x_g(t,s,x)-\tilde a_{g_+}(t)|=0$.
In addition, $r_g$ is locally pullback repulsive.
\item[\rm(iii)] There exists a bounded solution $b\colon\R\to\R$ of \eqref{eq:4Cgtran}
if and only if $r_g$ and $a_g$ are globally defined, in which case $r_g\le b\le a_g$.
\item[\rm(iv)] If $a_g$ and $r_g$ are bounded
and different, then $(\tilde a_g,\tilde r_g):=(a_g,r_g)$ is an attractor-repeller pair
of solutions for \eqref{eq:4Cgtran}.
In particular, this is the situation if: there exists $s\in\R$ such that $a_g(s)$ and
$r_g(s)$ exist and satisfy $r_g(s)<a_g(s)$; or if $a_g\ne r_g$ and one of them is
bounded.
\item[\rm(v)] If \eqref{eq:4Cgtran} has no hyperbolic solutions, then
it has at most one bounded solution $a_g=r_g$.
\end{itemize}
\end{teor}
\begin{proof}
If \hyperlink{gc1}{\bf gc1} and \hyperlink{gc3}{\bf gc3} hold, then
there exist $m>0$ and $\ep>0$ such that $g(t,\pm x)\le-\ep$ and $g_\pm(t,\pm x)\le-\ep$
for all $t\in\R$ if $x\ge m$. So, we can repeat the proofs of points
(i) and (ii) of \cite[Theorem 3.1]{lnor}, in order to proof that, if there are solutions which
remain bounded as time decreases (resp.~as time increases), then there exists the
(possibly local) map $a_g\colon(-\infty,s_a)\to(-\infty,m]$, with $-\infty<s_a\le\infty$
(resp.~$r_g\colon(s_r,\infty)\to[-m,\infty)$ with and $-\infty\le s_r<\infty$),
characterized as in the first assertion of (i) (resp.~of (ii)).
The same arguments show that $m$ is a bound for the absolute value
of any bounded solution of $x'=g_\pm(t,x)$. Hence, $-m\le\tilde r_{g_\pm}(t)\le
\tilde a_{g_\pm}(t)\le m$ for all $t\in\R$.
We also point out that the characterizations of $a_g$ and $r_g$ combined
with the existence of an upper bound for $a_g$ and a lower bound for $r_g$ prove (iii).
\par
Now, we proceed as in the proof of \cite[Theorem 3.4]{lno2}. We detail it,
since the scenario here is much more general, and some technical differences arise.
\par
Let us take $\ep>0$. Since \hyperlink{gc1}{\bf gc1} holds,
Theorem \ref{th:2persistencia} provides $\delta_-=\delta_-(\ep)>0$
such that, if $f$ is $C^1$-admissible and
$\n{g_--f}_{1,m}<\delta_-$,
then $x'=f(t,x)$ has an attractor-repeller pair
$(\tilde a_f,\tilde r_f)$ with $\|\tilde a_{g_-}-\tilde a_f\big\|_\infty\le\ep$ and
$\|\tilde r_{g_-}-\tilde r_f\big\|_\infty\le\ep$. It also ensures the existence
of a common dichotomy pair for all these functions $f$. We
fix the same for both hyperbolic solutions: $(k_\ep,\beta_\ep)$.
\par
We choose $t^-=t^-(\ep)<0$ such that $|g(t,x)-g_-(t,x)|<\delta_-/2$ and
$|g_x(t,x)-(g_-)_x(t,x)|<\delta_-/2$ if $t\le t^-$ and $|x|\le m$ (see Lemma \ref{lem:4Cderivadas}),
and define $f_-(t,x)$ as $g(t,x)$ if $t<t^-$ and as $g_-(t,x)-g_-(t^-,x)+g(t^-,x)$
otherwise. It is easy to check that $f_-$ is $C^1$-admissible and $\n{g_--f_-}_{1,m}\le\delta_-$,
and so $x'=f_-(t,x)$ has an attractor-repeller pair $(\tilde a_{f_-},\tilde r_{f_-})$, with
$\big\|\tilde a_{f_-}-\tilde a_{g_-}\big\|_\infty\le\ep$ and
$\big\|\tilde r_{f_-}-\tilde r_{g_-}\big\|_\infty\le\ep$.
\par
Let us now define $\hat a_{f_-}$ as the solution of $x'=g(t,x)$
with $\hat a_{f_-}(t^-)=\tilde a_{f_-}(t^-)$. We will check that
$\hat a_{f_-}=a_g$. Since $\hat a_{f_-}(t)=\tilde a_{f_-}(t)$
for $t\le t^-$, it remains bounded as $t$ decreases, which, as seen before,
ensures that $a_g$ exists and
that $\hat a_{f_-}\le a_g$. To prove that $\hat a_{f_-}\ge a_g$,
we take $x>\hat a_{f_-}(t^-)$ in order to check that $x_g(t,t^-,x)$
is unbounded as time decreases: Lemma \ref{lema:4Chiposc} guarantees that
the map $\mf_-$ defined on the hull of $f_-$ satisfies the
hypotheses of Proposition \ref{prop:3Ccoer}; hence, this result
ensures that the solution $x_{f_-}(t,t^-,x)$ of $x'=f_-(t,x)$
is unbounded as time decreases; and the assertion follows from here and from
$x_g(t,t^-,x)=x_{f_-}(t,t^-,x)$ for $t\le t^-$.
\par
Note that we have proved that $\lim_{t\to-\infty}(a_g(t)-\tilde a_{g_-}(t))=0$.
On the other hand, if $x<a_g(s)$, then there exists $t_0<t^-=t^-(\ep)$
such that $x_g(t_0,s,x)<a_g(t_0)=\tilde a_{f_-}(t_0)$.
Since $x_g(t,s,x)=x_g(t,t_0,x_g(t_0,s,x))$ solves
$x'=f_-(t,x)$ for $t\le t_0$, we conclude from Theorem~\ref{th:3Cequiv} that
$\lim_{t\to-\infty}|x_g(t,s,x)-\tilde r_{f_-}(t)|=0$.
Therefore, $|x_g(t,s,x)-\tilde r_{g_-}(t)|<2\,\ep$ if $t\le t^-(\ep)$,
which ensures that $\lim_{t\to-\infty}|x_g(t,s,x)-\tilde r_{g_-}(t)|=0$.
To check that $a_g$ is locally pullback attractive (and so complete the proof of (i)),
we observe that the attractive hyperbolicity of $\tilde a_{f_-}$ provides $\delta>0$, $k\ge 1$ and $\gamma>0$
such that $|a_g(t)-x_g(t,s,a_g(s)\pm\delta)|=|\tilde a_{f_-}(t)-x_{f_-}(t,s,\tilde a_{f_-}(s)\pm\delta)|
\le k\,\delta \,e^{-\gamma(t-s)}$ if $s\le t^-$ and $t\in [s,t^-]$
(see, e.g., \cite[Proposition 2.1]{dno3}).
\par
Analogous arguments prove (ii). If $a_g$ and $r_g$
are bounded and different, then (iii) yields $r_g<a_g$, and hence
the limiting properties established in (i) and (ii) prove
that they are uniformly separated. Therefore, Theorem \ref{th:3Cequiv}
proves that they form an attractor-repeller pair. The last assertions in (iv)
follow from (i), (ii) and (iii), and the bounds $-m\le r_g$ and $a_g\le m$.
Finally, (v) follows easily from (iv).
\end{proof}
\begin{teor}\label{th:4Ccasos}
Assume that $g$ satisfies \hyperlink{gc1}{\bf gc1}-\hyperlink{gc5}{\bf gc5}, let
$(\tilde a_{g_\pm},\tilde r_{g_{\pm}})$ be
the attractor-repeller pairs of solutions of $x'=g_\pm(t,x)$
given by \hyperlink{gc5}{\bf gc5}, and let $a_g$ and $r_g$ be the solutions of
\eqref{eq:4Cgtran} provided by Theorem {\rm \ref{th:4Cdossoluciones}}.
Then, the dynamics of the transition equation \eqref{eq:4Cgtran} fits in
one of the following dynamical scenarios:
\begin{itemize}[leftmargin=12pt]
\item \hypertarget{CCA}{{\sc Case A}}: there exists an attractor-repeller pair of solutions $(\tilde a_g,\tilde r_g)$,
with $\tilde a_g:=a_g$ and $\tilde r_g:=r_g$.
In this case, $\lim_{t\to\pm\infty}(\tilde r_g(t)-\tilde r_{g_\pm}(t))=0$
and $\lim_{t\to\pm\infty}(\tilde a_g(t)-\tilde a_{g_\pm}(t))=0$.
\item \hypertarget{CCB}{{\sc Case B}}: there are bounded solutions but no hyperbolic ones. In this case,
$r_g=a_g$ is the unique bounded solution, and it is locally pullback attractive and repulsive.
\item \hypertarget{CCC}{{\sc Case C}}: there are no bounded solutions.
\end{itemize}
\end{teor}
\begin{proof}
Assume that \hyperlink{CCC}{\sc Case C} does not hold; i.e., that there exists a
bounded solution of \eqref{eq:4Cgtran}. Theorem \ref{th:4Cdossoluciones}(iii)
ensures that $r_g$ and $a_g$ are bounded. If they are different, point (iv) of Theorem
\ref{th:4Cdossoluciones} ensures that they form an attractor-repeller pair of solutions,
and points (i) and (ii) yield the asymptotic behaviour described in \hyperlink{CCA}{\sc Case A}.
If, on the contrary, $a_g=r_g$, then Theorem \ref{th:4Cdossoluciones}(iii) ensures that $a_g=r_g$ is the
unique bounded solution. As stated in Theorem \ref{th:2persistencia}, its hyperbolicity
would ensure its exponential asymptotic stability as time either increases or decreases,
which contradicts either point (ii) or (i) of Theorem \ref{th:4Cdossoluciones}.
Hence, \hyperlink{CCB}{\sc Case B} holds.
\end{proof}
Let us analyze part of the information provided by Theorems \ref{th:4Cdossoluciones} and
\ref{th:4Ccasos}. In all the cases, the locally pullback attractive
solution $a_g$ of the transition equation ``connects" with the
attractive hyperbolic solution of the past equation as time decreases.
The differences arise with its behavior in the future: in \hyperlink{CCA}{\sc Case A},
usually referred to as ({\em end-point})
{\em tracking}, $\tilde a_g:=a_g$ also connects with the
attractive hyperbolic solution of the future as time increases, while in
\hyperlink{CCC}{\sc Case C}, of {\em tipping}, $a_g$ is unbounded, and hence
the connection is lost. In the extremely
unstable \hyperlink{CCB}{\sc Case B}, $a_g$ is still bounded but it connects with the repulsive
hyperbolic solution of the future. The interested reader can in find
\cite[Figures 1-6]{lnor} some drawings depicting the dynamical behavior in
each one of these three cases. (There is a typo there: the
graphs of \hyperlink{CCA}{\sc Cases A} and \hyperlink{CCC}{C} are interchanged).
\par
The next result provides a useful comparison criterion ensuring \hyperlink{CCA}{\sc Case A}.

\begin{prop} \label{prop:4CCaseA}
Assume that $g$ satisfies \hyperlink{gc1}{\bf gc1}-\hyperlink{gc5}{\bf gc5}.
If there exists a continuous map
$h\colon\RR\to\R$ such that $h(t,x)\le g(t,x)$ for all $(t,x)\in\RR$ and
\begin{itemize}[leftmargin=10pt]
\item[-] either $x'=h(t,x)$ has two different bounded solutions,
\item[-] or $x'=h(t,x)$ has a bounded solution which does not solve $x'=g(t,x)$,
\end{itemize}
then \eqref{eq:4Cgtran} is in \hyperlink{CCA}{\sc Case A}.
\end{prop}
\begin{proof}
In both cases, Proposition \ref{prop:3Ccoer}(v) ensures that
\eqref{eq:4Cgtran} has two bounded solutions, so that
Theorem \ref{th:4Ccasos} proves the assertion.
\end{proof}
As explained in the Introduction, a {\em critical transition} (or {\em tipping point})
occurs when a small variation on the external input of the equation causes a dramatic
variation on the dynamics. We will focus on critical transitions associated
to one-parametric families of equations which occur when the dynamics
moves from \hyperlink{CCA}{\sc Case A} to \hyperlink{CCC}{\sc Case C} of Theorem \ref{th:4Ccasos} as the
parameter crosses a particular {\em critical value}. Theorem \ref{th:4Csaddle-node}
shows that, if the parametric variation
is smooth enough, this transition means \hyperlink{CCB}{\sc Case B} for the critical value,
and that these
critical transitions can be understood as nonautonomous saddle-node bifurcations: they
occur as a consequence of the collision of an attractive hyperbolic solution with a repulsive
one as $c$ varies. It also shows the
persistence of \hyperlink{CCA}{\sc Cases A} and \hyperlink{CCC}{C}.
\begin{teor} \label{th:4Csaddle-node}
Let $\mC\subseteq\R$ be an open interval, and let $\bar g\colon\R\times\R\times\mC\to\R$
be a map such that $g^c(t,x):=\bar g(t,x,c)$ satisfies
\hyperlink{gc1}{\bf gc1}-\hyperlink{gc5}{\bf gc5} for all $c\in\mC$.
Let $\bar g_x$ be the partial derivative with respect to the second variable
and assume that $\bar g$ and $\bar g_x$ are admissible on $\R\times\R\times\mC$. Assume also that
$\limsup_{x\to\pm\infty} \bar g(t,x,c)<0$ uniformly on $\R\times\mJ$ for any compact interval
$\mJ\subset\mC$.
\begin{itemize}
\item[\rm(i)] Assume that there exist $c_1,c_2$ in $\mC$ with
$c_1<c_2$ such that the dynamics of $x'=g^c(t,x)$ is in \hyperlink{CCA}{\sc Case A} for $c=c_1$
and not for $c=c_2$. If $c_0:=\inf\{c>c_1\,|\;\text{\hyperlink{CCA}{\sc Case A} does not hold\/}\}$,
then $c_0>c_1$. Let $(\tilde a_{g^c},\tilde r_{g^c})$ be the attractor-repeller
pair for $c\in[c_1,c_0)$. Then, the dynamics of $x'=g^{c_0}(t,x)$ is in
\hyperlink{CCB}{\sc Case B}, and $\lim_{c\to c_0^-}(\tilde a_{g^c}(t)-\tilde r_{g^c}(t))=0$
for all $t\in\R$. The result is analogous if $c_1>c_2$.
\item[\rm(ii)] Assume that there exist $c_3,c_4$ in $\mC$ with
$c_3<c_4$ such that the dynamics of $x'=g^c(t,x)$ is in \hyperlink{CCC}{\sc Case C} for $c=c_3$
and not for $c=c_4$. If $c_0:=\inf\{c>c_3\,|\;\text{\hyperlink{CCC}{\sc Case C} does not hold}\}$,
then $c_0>c_3$, and the dynamics of $x'=g^{c_0}(t,x)$ is in
\hyperlink{CCB}{\sc Case B}. The result is analogous if $c_3>c_4$.
\end{itemize}
\end{teor}
\begin{proof}
(i) The admissibility hypotheses ensure that,
for $c\in\mC$, $\rho>0$ and $\delta>0$ fixed, there exists $\ep_0>0$ such that
\[
 \sup_{(t,x)\in\R\times[-\rho,\rho]}|g^c(t,x)-g^{c+\ep}(t,x)|
 +\sup_{(t,x)\in\R\times[-\rho,\rho]}|g^c_x(t,x)-g^{c+\ep}_x(t,x)|<\delta
\]
if $|\ep|\le\ep_0$. Hence, Theorems \ref{th:4Ccasos} and \ref{th:2persistencia} guarantee
the persistence of \hyperlink{CCA}{\sc Case A} under small variations of $c$, which in turn
ensures that $c_0>c_1$ and that $x'=g^{c_0}(t,x)$ is not in \hyperlink{CCA}{\sc Case A}.
(Note that the last condition on coercivity is not yet required.)
On the other hand, the hypothesis on $\limsup_{x\to\pm\infty}\bar g(t,x,c)$
(which is stronger than ``\hyperlink{gc3}{\bf gc3} for all $c$")
ensures the existence of a constant $m>0$ and $\delta>0$ such that
$g^c(t,x)\le-\delta$ if $t\in\R$, $|x|>m$, and $c\in\mC$ is
close enough to $c_0$. This fact allows us to reason as in the proof of
Proposition \ref{prop:3Ccoer}(ii) in order to check that the lower and upper
bounded solutions ($\tilde r_{g^c}$ and $\tilde a_{g^c}$) of
$x'=g^c(t,x)$ are lower bounded by $-m$ and
upper bounded by $m$ for $c\in\mC$ close enough to $c_0$. Hence,
for a common sequence $(c_n)\uparrow c_0$, there exist
$\bar r_0:=\lim_{n\to\infty}\tilde r_{g^{c_n}}(0)$ and
$\bar a_0:=\lim_{n\to\infty}\tilde a_{g^{c_n}}(0)$. It is easy to deduce that
the solutions of $x'=g^{c_0}(t,x)$ with values $\bar r_0$ and $\bar a_0$
at $t=0$ are bounded. Hence we are in \hyperlink{CCB}{\sc Case B}, and both solutions
coincide. It is easy to check that this unique bounded solution, $b^{c_0}$,
satisfies $b^{c_0}(t)=\lim_{n\to\infty}\tilde a_{g^{c_n}}(t)=
\lim_{n\to\infty}\tilde r_{g^{c_n}}(t)$
for all $t\in\R$, which combined with the uniqueness of $b^{c_0}$ guarantees that
$\lim_{c\to c_0^-}(\tilde a_{g^c}(t)-\tilde r_{g^c}(t))=0$ for all $t\in\R$, as asserted.
It is clear that the argument can be repeated if $c_1>c_2$.
\smallskip\par
(ii) We assume for contradiction the existence of $(c_n)\downarrow c_3$ such that
there exists a bounded solution $b^{c_n}$ for all $n$, and reason as before
to conclude that $\bar b_3:=\lim_{n\to\infty} b^{c_n}(0)$ is finite and provides
the value at $0$ of a bounded solution for $c_3$, impossible. This shows
the persistence of \hyperlink{CCC}{\sc Case C} under small variations in $c$, which in turn
ensures that $c_0>c_3$ and that $x'=g^{c_0}(t,x)$ is not in \hyperlink{CCC}{\sc Case C}.
The previously proved persistence of \hyperlink{CCA}{\sc Case A} proves the last assertion.
And the argument can be repeated if $c_3>c_4$.
\end{proof}
We complete this part with a consequence of Proposition \ref{prop:4CCaseA}
which ensures the existence of at most a unique
tipping point for certain parametric families:
\begin{coro} \label{coro:4Cunpunto}
Let $\mC\subseteq\R$ be an open interval and let
$\{g^c\,|\;c\in\mC\}$ be a family of functions satisfying
\hyperlink{gc1}{\bf gc1}-\hyperlink{gc5}{\bf gc5}.
Assume that there exists $c_0\in\mC$ such that the dynamics of $x'=g^{c_0}(t,x)$
is in \hyperlink{CCB}{\sc Case B}, and such that,
for all $c_-,c_+\in\mC$ with $c_-<c_0<c_+$:
$g^{c_-}(t,x)\le g^{c_0}(t,x)\le g^{c_+}(t,x)$
for all $(t,x)\in\RR$; and there exist $t_{c_-}$ and $t_{c_+}$ such that the first
and second inequality
are strict for $t=t_{c_-}$ and $t=t_{c_+}$ (respectively) and all $x\in\R$.
Then, $x'=g^c(t,x)$ is in \hyperlink{CCC}{\sc Case C} for $c\in\mC$ with $c<c_0$
and in \hyperlink{CCA}{\sc Case A} for $c\in\mC$ with $c>c_0$.
\end{coro}
\begin{proof}
The hypotheses ensure that, if $c_+>c_0$, any bounded solution of $x'=g^{c_0}(t,x)$
does not solve $x'=g^{c_+}(t,x)$, and hence
Proposition \ref{prop:4CCaseA} shows that $x'=g^c(t,x)$ is in
\hyperlink{CCA}{\sc Case A} for $c>c_0$ in $\mC$.
Analogously, if $c_-<c_0$, any bounded solution of $x'=g^{c_-}(t,x)$
does not solve $x'=g^{c_0}(t,x)$, and hence
Proposition \ref{prop:4CCaseA} also shows that $x'=g^{c_0}(t,x)$
would be in \hyperlink{CCA}{\sc Case A} (which is not true, by hypothesis)
if $x'=g^c(t,x)$ were in \hyperlink{CCA}{\sc Cases A} or \hyperlink{CCB}{B}
for $c<c_0$ in $\mC$.
\end{proof}
\subsection{Some scenarios of critical transitions in the concave case}\label{subsec:41}
Let $\mI\subseteq\R$ be an open interval, and let the functions $f\colon\RR\times\mI\to\R$ and
$\G,\G_-,\G_+\colon\RR\to\R$ satisfy
\begin{list}{}{\leftmargin 20pt}
\item[\hypertarget{fc1}{\bf fc1}] there exist the derivatives $f_x$
and $f_\gamma$, and $f$, $f_x$ and $f_\gamma$ are admissible on $\R\times\R\times\mI$.
\item[\hypertarget{fc2}{\bf fc2}] $\G,\G_-$ and $\G_+$ take values in $[a,b]\subset\mI$,
are $C^1$-admissible,
and $\lim_{t\to\pm\infty}(\G(t,x)-\G_\pm(t,x))=0$ uniformly on each compact subset $\mJ\subset\R$.
\item[\hypertarget{fc3}{\bf fc3}] $\limsup_{x\to\pm\infty} f(t,x,\gamma)<0$
uniformly in $(t,\gamma)\in\R\times\mJ$ for all compact interval $\mJ\subset\mI$.
\item[\hypertarget{fc4}{\bf fc4}]
$\inf_{t\in\R}\big((\partial/\partial x)f(t,x,\G_\pm(t,x))|_{x=x_1}-
(\partial/\partial x)f(t,x,\G_\pm(t,x))|_{x=x_2}\big)>0$ whenever $x_1<x_2$.
\item[\hypertarget{fc5}{\bf fc5}] Each equation $x'=f(t,x,\G_\pm(t,x))$ has
two hyperbolic solutions $\tilde r_{\G_\pm}<\tilde a_{\G_{\pm}}$.
\end{list}
\par
Observe that condition \hyperlink{fc2}{\bf fc2} allows us to understand the equations
\begin{equation}\label{eq:4Clim}
 x'=f(t,x,\G_-(t,x)) \quad\text{and}\quad x'=f(t,x,\G_+(t,x))
\end{equation}
as the ``past" and ``future" of
\begin{equation}\label{eq:4Ctran}
 x'=f(t,x,\G(t,x))\,.
\end{equation}
\par
We will say that the pair $(f,\G)$ satisfies
\hyperlink{fc1}{\bf fc1}-\hyperlink{fc5}{\bf fc5} whenever there exist
$\G_\pm$ such all the listed properties hold. Note that, in this case,
also the pairs $(f,\G_\pm)$ satisfy \hyperlink{fc1}{\bf fc1}-\hyperlink{fc5}{\bf fc5}.
We omit the almost immediate proof of the next result, which shows that the previous ones
apply to the current setting.
\begin{prop}\label{prop:4Chiposcomp}
Assume that $(f,\G)$ satisfies \hyperlink{fc1}{\bf fc1}-\hyperlink{fc5}{\bf fc5}.
Then, the maps $g,\,g_-$ and $g_+$ given by $g(t,x):=f(t,x,\G(t,x))$, $g_-(t,x):=f(t,x,\G_-(t,x))$
and $g_+(t,x):=f(t,x,\G_+(t,x))$ satisfy the conditions
\hyperlink{gc1}{\bf gc1}-\hyperlink{gc5}{\bf gc5}. Therefore, the dynamical
possibilities for \eqref{eq:4Ctran} are those described in Theorem {\rm\ref{th:4Ccasos}}.
\end{prop}
\begin{notas}\label{rm:4Cpar}
1.~It is easy to check that the proof of Proposition \ref{prop:4Chiposcomp} can be repeated
in the next cases: if we remove the boundedness of $\G$ and $\G_\pm$
from condition \hyperlink{fc2}{\bf fc2} but assume that $\mI=\R$ and that
the limit in \hyperlink{fc3}{\bf fc3}
is uniform in $(t,\gamma)\in\RR$; and if we remove the assumptions on the derivative
$f_\gamma$ of \hyperlink{fc1}{\bf fc1} but assume that $\G$, and hence $\G_\pm$,
depend only on $t$.
Hence, the conclusions of Theorems \ref{th:4Cdossoluciones}
and \ref{th:4Ccasos} also hold under these conditions.
\par
2.~As explained in Remark \ref{rm:4Creformular}, Proposition \ref{prop:4Chiposcomp}
applied to the pairs $(f,\G_\pm)$
allows us to reformulate condition \hyperlink{fc5}{\bf fc5} as:
``each equation $x'=f(t,x,\G_{\pm}(t,x))$ has an attractor-repeller pair
of solutions", which determines its global dynamics.
\end{notas}
\par
In this subsection (as in Subsection \ref{subsec:61}), we analyze
some mechanisms of occurrence (or lack) of tipping points for
transition equations \eqref{eq:4Ctran} due to small parametric
variations in the transition function: we work with
one-parametric families
\begin{equation}\label{eq:4conc}
 x'=f(t,x,\G^c(t,x))\,.
\end{equation}
Let us mention three of the large variety of physical mechanisms that
may cause critical transitions:
\begin{itemize}[leftmargin=*]
\item[-] \emph{Rate-induced critical transitions:}
if $\G^c(t,x)=\G(ct,x)$ for a fixed $\G$ and any $c>0$, then
the parameter $c>0$ determines the speed of the transition $\G^c$.
In order to have a past and a future independent of the rate, we require $\G_-$
and $\G_+$ to be independent of $t$. So, a larger $c$ means a significant distance
from $\G(ct,x)$ to $\G_\pm(x)$ during a shorter period.
\item[-] \emph{Phase-induced critical transitions:}
if $\G^c(t,x)=\G(c+t,x)$, then the parameter $c\in\R$ represents the initial phase of the
transition function. As before, we assume $\G_-$ and $\G_+$ independent of $t$.
\item[-] \emph{Size-induced critical transitions:} with $\G_-\equiv 0$ and
$\G^c(t,x):=c\,\G(t,x)$,
different values of $c>0$ mean different sizes of the transition function
which ``takes" $x'=f(t,x,0)$ to $x'=f(t,x,c\,\G_+(t,x))$.
\end{itemize}
\par
The next result establishes conditions on a parametric family of maps $\{\G^c\}$ and $f$
which are enough to guarantee the
persistence of \hyperlink{CCA}{\sc Cases A} and \hyperlink{CCC}{\sc C}, and to show that
the occurrence of a critical transition means the occurrence of \hyperlink{CCB}{\sc Case B}
and can be understood as a nonautonomous saddle-node bifurcation: see Theorem \ref{th:4Csaddle-node}.
We omit the (easy) proof.
\begin{teor} \label{th:4Csaddle-node-gamma}
Let $\mC\subseteq\R$ be an open interval, and let the maps
$\{\G^c\,|\;c\in\mC\}$ be a family of functions such that
all the pairs $(f,\G^c)$ satisfy \hyperlink{fc1}{\bf fc1}-\hyperlink{fc5}{\bf fc5}
and such that $\RR\times\mC\to\R,\,(t,x,c)\mapsto \G^c(t,x)$ is admissible.
Assume also that, for any $c\in\mC$, there exists $\delta_c>0$ such that
$\sup_{(t,x)\in\RR,\,|\ep|\le\delta_c}|\G^{c+\ep}(t,x)|<\infty$.
Then, the map $\bar g(t,x,c):=f(t,x,\G^c(t,x))$ satisfies
all the hypotheses of Theorem {\rm\ref{th:4Csaddle-node}}.
\end{teor}
\begin{nota}\label{rm:4basta}
Note that if, in the considered case of rate and phase variation,
with $\G_\pm$ independent of $t$, all the pairs $(f,\G^c)$ satisfy
\hyperlink{fc1}{\bf fc1}-\hyperlink{fc5}{\bf fc5} if $(f,\G)$ does, with the same maps $\G_\pm$.
The same occurs in the size-variation case if we also assume $\G_+\equiv 0$. In addition, in
the three considered cases, the admissibility and
boundedness hypotheses of Theorem \ref{th:4Csaddle-node-gamma} also hold.
\end{nota}
In the rest of this subsection, we will describe conditions ensuring the lack of
rate-induced and phase-induced critical transitions (in Theorem \ref{th:4Cnotran}),
as well as the occurrence of size-induced critical transitions (in Theorem \ref{th:4Csize}).
These scenarios assume monotonicity of $f$ with respect to $\gamma$.
Theorem \ref{th:4Cnotran} establishes conditions on $f(t,x,\gamma_0)$
ensuring that $[\gamma_0,\infty)$ or $(-\infty,\gamma_0]$ is a {\em safety halfline\/}: if
$\G$ takes values in it, then neither rate-induced tipping nor
phase-induced tipping occurs.
\begin{teor} \label{th:4Cnotran}
Assume that $(f,\G)$ satisfy \hyperlink{fc1}{\bf fc1}-\hyperlink{fc5}{\bf fc5}. Assume also that
$\gamma\mapsto f(t,x,\gamma)$ is nondecreasing (resp.~nonincreasing)
for all $(t,x)\in\RR$ and that one of the following situations holds:
\begin{itemize}[leftmargin=20pt]
\item[\rm(1)] there exists a constant $\gamma_0\le\G$
(resp.~$\gamma_0\ge\G$) such that $x'=f(t,x,\gamma_0)$
has either two different bounded solutions or
a bounded solution which does not solve \eqref{eq:4Ctran};
\item[\rm(2)] there exists a continuous map $\Delta\colon\RR\to\R$
with $\Delta\le\G$ (resp.~$\Delta\ge\G$)
such that $x'=f(t,x,\Delta(t,x))$ has either two different bounded solutions, or
a bounded solution which does not solve \eqref{eq:4Ctran};
\end{itemize}
Then, \eqref{eq:4Ctran} is in \hyperlink{CCA}{\sc Case A}.
\par
In particular, if {\rm (1)} or {\rm (2)} holds, and if we assume in addition that $\G_\pm$
do not depend on $t$, then the equations $x'=f(t,x,\G(ct,x))$ and
$x'=f(t,x,\G(t+c,x))$ are in \hyperlink{CCA}{\sc Case A} for all $c>0$ and $c\in\R$,
respectively: there are neither rate-induced nor phase-induced critical transitions.
\end{teor}
\begin{proof}
All the assertions follow easily from Propositions \ref{prop:4Chiposcomp} and \ref{prop:4CCaseA},
and from Remark \ref{rm:4basta}.
\end{proof}
\begin{nota}
In the increasing (resp.~decreasing) scenario of Theorem \ref{th:4Cnotran},
condition \hyperlink{fc5}{\bf fc5} ensures (2) if either $\G_-(t,x)\le\G(t,x)$ or $\G_+(t,x)\le\G(t,x)$ (resp.
either $\G_-(t,x)\ge\G(t,x)$ or $\G_+(t,x)\ge\G(t,x)$) for all $(t,x)\in\RR$: it
suffices to take $\Delta=\G_-$ or $\Delta=\G_+$.
\end{nota}
Our next result provides another scenario of lack of critical transitions
or occurrence of exactly one. Now, the variation of the parameter precludes the transition
map to remain always in the safety half-line, and hence the dynamics
cannot be always (if ever) in \hyperlink{CCA}{\sc Case A}.
The interval $\mI$ of variation of the third argument of $f$ must be $\R$.
\begin{teor} \label{th:4Csize}
Assume that $\mI=\R$. Let $\G\colon\RR\to\R$ and $\G_0\colon\RR\to[0,\infty)$
be globally bounded and $C^1$-admissible, and such
that the pair $(f,\G+d\,\G_0)$ satisfies \hyperlink{fc1}{\bf fc1}-\hyperlink{fc5}{\bf fc5}
for all $d\in\R$. Assume that $\G_0(t_0,x)>0$ for all $x\in\R$ and a $t_0\in\R$. Assume also that
$\gamma\mapsto f(t,x,\gamma)$ is strictly increasing on $\R$ for all $(t,x)\in\RR$,
with $\lim_{\gamma\to-\infty} f(t,x,\gamma)=-\infty$ uniformly on compact sets of $\RR$.
Then, either
\begin{equation}\label{eq:4CGG_0}
 x'=f(t,x,\G(t,x)+d\,\G_0(t,x))
\end{equation}
is in \hyperlink{CCC}{\sc Case C} for all $d\in\R$, or there exists $d_0$ such that
\eqref{eq:4CGG_0} is in \hyperlink{CCA}{\sc Case A} for $d>d_0$,
in \hyperlink{CCB}{\sc Case B} for $d=d_0$ and in \hyperlink{CCC}{\sc Case C} for $d<d_0$.
\end{teor}
\begin{proof}
Let us assume for contradiction that \eqref{eq:4CGG_0}$_d$ is in
\hyperlink{CCA}{\sc Case A} for all $d\in\R$, and fix $\bar d\in\R$.
We take $\delta>0$ and $m_1,m_2\in\R$ such that
$f(t,x,\G(t,x)+\bar d\,\G_0(t,x))\le-\delta$ for all $t\in\R$ if $x\notin(m_1,m_2)$.
Since $f$ is nondecreasing in $\gamma$ and $\G_0\ge 0$,
the map $d\mapsto f(t,x,\G(t,x)+d\,\G_0(t,x))$ is nondecreasing, and hence, for all $d\le\bar d$,
$f(t,x,\G(t,x)+d\,\G_0(t,x))\le-\delta$ for all $t\in\R$ if $x\notin(m_1,m_2)$.
As in the proof of Proposition~\ref{prop:3Ccoer}(ii),
we check that $m_1\le \tilde a_d\le m_2$ if $d\le\bar d$, where $\tilde a_d$ is the upper bounded
solution of \eqref{eq:4CGG_0}$_d$.
We look for $t_1<t_0<t_2$ and $k>0$ such that
$\G_0(t,x)>k$ if $t\in[t_1,t_2]$ and $x\in[m_1,m_2]$, and call
$k_d:=\sup_{(t,x)\in[t_1,t_2]\times[m_1,m_2]}f(t,x,\G(t,x)+d\,\G_0(t,x))$
for $d\le\bar d$. Then,
$ k_d\le\sup_{(t,x)\in[t_1,t_2]\times[m_1,m_2]}f(t,x,\G(t,x)+dk)$
if $d\le\min(0,\bar d)$, which combined with the hypothesis on $\lim_{\gamma\to-\infty} f(t,x,\gamma)$
ensures that $\lim_{d\to-\infty} k_d=-\infty$.
Hence, $m_1-m_2\le\tilde a_d(t_2)-\tilde a_d(t_1)\le (t_2-t_1)\,k_d$ for all $d\le\bar d$,
which is impossible.
\par
Now we assume that \eqref{eq:4CGG_0}$_d$ is not in \hyperlink{CCA}{\sc Case C} for all
$d\in\R$. Theorem \ref{th:4Csaddle-node}
ensures the existence of $d_0$ such that
\eqref{eq:4CGG_0}$_{d_0}$ is in \hyperlink{CCB}{\sc Case B}, and Corollary \ref{coro:4Cunpunto},
all whose hypotheses are satisfied, completes the proof.
\end{proof}
By reviewing the proof of Theorem \ref{th:4Csize}, we observe that we have in fact proved the next result,
which considerably weakens the conditions of the previous one (but with
a statement quite harder to read, so we keep both of them).
\begin{teor} \label{th:4Csize2}
Assume that $\mI=\R$. Let $\G\colon\RR\to\R$ and $\G_0\colon\RR\to[0,\infty)$
be globally bounded and $C^1$-admissible, and such
that the pair $(f,\G+d\,\G_0)$ satisfies \hyperlink{fc1}{\bf fc1}-\hyperlink{fc5}{\bf fc5}
for all $d\in\R$. Assume that there exists
$\bar d\in\R$ such that
\begin{equation}\label{eq:4CGG_02}
x'=f(t,x,\G(t,x)+d\,\G_0(t,x))
\end{equation}
is in \hyperlink{CCA}{\sc Cases A} or \hyperlink{CCB}{\sc B} for $d=\bar d$.
Let $\delta>0$ and $m_1,m_2\in\R$ satisfy
$f(t,x,\G(t,x)+\bar d\,\G_0(t,x))\le-\delta$ for all $t\in\R$ if $x\notin(m_1,m_2)$.
Assume that there exists $t_0$
such that $\G_0(t_0,x)>0$ for all $x\in[m_1,m_2]$,
that $\gamma\mapsto f(t,x,\gamma)$ is nondecreasing for all $(t,x)\in\RR$ and strictly increasing
for $(t,x)\in\R\times[m_1,m_2]$, with
$\lim_{\gamma\to-\infty} f(t,x,\gamma)=-\infty$ uniformly on compact sets of $\R\times[m_1,m_2]$.
Then, there exists $d_0\le \bar d$ such that
\eqref{eq:4CGG_02} is in \hyperlink{CCA}{\sc Case A} for $d>d_0$,
in \hyperlink{CCB}{\sc Case B} for $d=d_0$ and in \hyperlink{CCC}{\sc Case C} for $d<d_0$.
\end{teor}
\begin{nota}\label{rm:4Ccasoparticular}
Frequently, the limit maps providing condition \hyperlink{fc2}{\bf fc2}
for all $d$ are $\G_\pm+d\,\G_{0,\pm}$ for $C^2$-admissible maps
$\G_\pm$ and $\G_{0,\pm}\ge 0$.
If so, condition \hyperlink{fc5}{\bf fc5} for all $x'=f(t,x,\G_\pm(t,x)+d\,\G_{0,\pm}(t,x))$
is only possible if, for any $t_0\in\R$, each map $x\mapsto\G_{0,\pm}(t_0,x)$ vanishes at least
for an $x^\pm_{t_0}\in[m_1,m_2]$: otherwise, Theorem \ref{th:4Csize2} precludes
hyperbolic solutions if $-d$ is large enough. Also often, $\G_\pm\equiv\G$ and $\G_{0,-}\equiv 0$,
and so Theorems \ref{th:4Csize} and \ref{th:4Csize2} study the occurrence of
size-induced critical transitions: just define $f^*(t,x,d\,\G_0(t,x)):=f(t,x,\G(t,x)+d\,\G_0(t,x))$.
\end{nota}
\subsection{Numerical simulations in asymptotically concave models}\label{subsec:42}
In this section, we consider a single-species population whose intrinsic dynamics is
governed by a nonautonomous quadratic equation (see \cite{dolo1})
subject to two external phenomena: migration and predation.
The inclusion of time-dependent intrinsic parameters and functions in the model
naturally arises when considering the influence of external factors which vary over
time, such as climatic or meteorological factors (see e.g.~\cite{renshaw1}).
In this way, the evolution of the population size is governed by
\begin{equation}\label{eq:4Cpredgeneral}
 x'=r(t)\,x\,\left(1-\frac{x}{K(t)}\right)+\Delta(t,x)\,,
\end{equation}
where we assume $r$ and $K$ to be positively bounded from below quasiperiodic functions,
and $\Delta$ to be $C^1$-admissible: so, if we define $h(t,x,\delta):=r(t)\,x\,(1-x/K(t))+\delta$,
then $h$ satisfies \hyperlink{fc1}{\bf fc1} and \hyperlink{fc3}{\bf fc3}.
In addition, we assume $(h,\Delta)$ to satisfy conditions \hyperlink{fc2}{\bf fc2},
\hyperlink{fc4}{\bf fc4} and \hyperlink{fc5}{\bf fc5}
for certain maps $\Delta_\pm$ which do not play a role for the moment.
The function $r$ represents the intrinsic growth rate of the species,
that is, the growth rate per individual in an ideal situation of unlimited resources;
$K$ does no longer represent, as in the autonomous case, the maximal population
allowed by the resources if $\Delta\equiv0$ (unless it is constant), but there exists
positive hyperbolic attractive solution which describes this target population;
and the function $\Delta$ quantifies
the net contribution per unit of time of both external phenomena: migration and predation.
In this framework, the locally pullback attractive solution $a_h$ of \eqref{eq:4Cpredgeneral} provided by
Theorem~\ref{th:4Cdossoluciones} (which can be applied, as Proposition~\ref{prop:4Chiposcomp}
ensures) describes the evolution of the target population during the transition.
The dynamical possibilities described by Theorem~\ref{th:4Ccasos} are:
the stable \hyperlink{CCA}{\sc Case A}, that means the survival of
the population $a_h$, which approaches the upper bounded
solution of the future equation as time increases;
the stable \hyperlink{CCC}{\sc Case C}, which means the
extinction of that population due to predation and migration;
and \hyperlink{CCB}{\sc Case B}, which is the highly unstable situation which
separates the other two.
\par
Before choosing a particular function $\Delta$, we apply Theorem~\ref{th:4Cnotran}(1)
to prove the existence of a safety halfline: a threshold $\delta_0<0$ such that, if
$\inf_{(t,x)\in\R\times\R}\Delta(t,x)\ge\delta_0$ and $\Delta\not\equiv\delta_0$,
then \eqref{eq:4Cpredgeneral} is in \hyperlink{CCA}{\sc Case A}.
That is, if the combined effect of predation and migration
never eliminates more than
$|\delta_0|$ individuals per unit of time, then the target population persists.
To apply Theorem~\ref{th:4Cnotran}, we use the information provided by
\cite[Theorem 3.6]{lnor} (which remains unchanged except for the bound for
$\lambda^*$ when multiplying the leading term by $r$),
according to which there exists $\delta_0\in\R$ such that
$x'=r(t)\,x\,(1-x/K(t))+\delta$ is in \hyperlink{CCA}{\sc Case A} if
$\delta>\delta_0$, in \hyperlink{CCB}{\sc Case B} if $\delta=\delta_0$,
and in \hyperlink{CCC}{\sc Case C} if $\delta<\delta_0$. In addition, since
we can take $\ep>0$ such that $0<r(t)\,\ep\,(1-\ep/K(t))$ for all $t\in\R$,
Proposition~\ref{prop:3Ccoer}(v) ensures that $x'=r(t)\,x\,(1-x/K(t))$ has two distinct
bounded solutions, and hence it is in \hyperlink{CCA}{\sc Case A}; that is, $\delta_0<0$.
Observe that this choice of $\delta_0$ is optimal for the application of Theorem~\ref{th:4Cnotran}(1),
since it is the smallest value of $\delta$ for which $x'=r(t)\,x\,(1-x/K(t))+\delta$ has
a bounded solution.
\par
In the following examples, we will give explicit expressions to the function $\Delta$,
depending on several bifurcation parameters, and we will find critical transitions: changes from
\hyperlink{CCA}{\sc Cases A} to \hyperlink{CCC}{\sc Case C} through \hyperlink{CCB}{\sc Case B}
as one of those parameters changes.
We will prove the uniqueness of almost all those critical transitions and find
numerical evidence of nonuniqueness of the remaining one.
\begin{exa}\label{ex:4Cejemplo}
We assume that the predation in \eqref{eq:4Cpredgeneral} can be suitably modeled by
a H\"{o}lling type III functional response term $-\gamma\, x^2/(b(t)+x^2)$ (see \cite{dolo1,dno3}),
where $\gamma$ represents the predator density and $b$ depends on the average time between
two attacks of a predator (which is related to the food processing time); and
we also assume that the net migration rate per unit of time $\phi(t)$ is negative
for all $t\in\R$. Therefore, the model is
\begin{equation}\label{eq:4Cpred}
 x'=r(t)\,x\left(1-\frac{x}{K(t)}\right)+\phi(t)-\gamma\,\frac{x^2}{b(t)+x^2}\,,
\end{equation}
where $-\phi$ and $b$ are positively bounded from below and quasiperiodic (as $r$ and $K$).
\par
Let us introduce a specific time-variation on the predator density $\gamma$:
we assume that the habitat is initially free of predators, that a certain time a group
of predators arrives in, and that all of them leave away after some time.
\begin{figure}
     \centering
         \includegraphics[width=\textwidth]{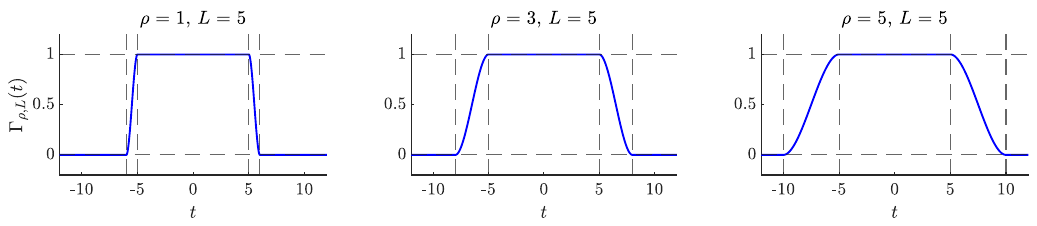}
         \caption{The $C^1$ transition map $\G_{\rho,L}$ for $L=5$ and several values of $\rho>0$.
         This map is defined as the unique $C^1$ cubic spline
         which takes value $1$ on $[-L,L]$ and 0 outside $[-L-\rho,L+\rho]$:
         if $Q(y):=2y^3-3y^2+1$, then $\G_{\rho,L}(t):=Q(-(t+L)/\rho)$ for $t\in[-L-\rho,-L]$ and
         $\G_{\rho,L}(t):=Q((t-L)/\rho)$ for $t\in[L,L+\rho]$.
         This map is increasing on $[-L-\rho,-L]$ and decreasing on $[L,L+\rho]$,
         and hence $\G_{\rho,L}(\cdot)$ is nondecreasing with respect to $L$ and with respect to~$\rho$.}
        \label{fig:Gammafunction}
\end{figure}
A wide range of causes can give rise to this transitory phenomenon in the case of predatory
birds: adverse winds, storms, orientation errors, changing attractiveness of the
breeding colony, etc. (see \cite{rappole}). A somehow related real-life example,
with foxes as predators,
can be found in the work \cite{oahgsar}, which describes the colonization of Punta de la Banya by
the Audouin's gull: an increasing population began to
severely decline from a certain time due to the arrival of foxes, whom later were
removed; and then the gull population began to increase again.
\par
To model the effect of this type of phenomena in a simple way, we
use a multiple of a $C^1$ approximation to the characteristic function of $[-L,L]$. More precisely,
we substitute the predator density parameter $\gamma$ in \eqref{eq:4Cpred} by the
compactly supported four-parametric transition function $t\mapsto d\,\G_{\rho,L}(t-p)$,
where $d\ge 0$ and $\G_{\rho,L}$ is the unique
$C^1$ cubic spline which takes the value 1 on $[-L,L]$
and 0 outside $[-L-\rho,L+\rho]$.
Its asymptotic limits are 0 for any choices of $\rho$ and $L$, and so,
the past and future equations coincide: the predation term disappears.
Figure~\ref{fig:Gammafunction} depicts $\Gamma_{\rho,L}$ for $L=5$ and some values of $\rho$.
Altogether, we get the four-parameter model
\begin{equation}\label{eq:4Cpredcompleta}
 x'=r(t)\,x\left(1-\frac{x}{K(t)}\right)+\phi(t)-d\,\G_{\rho,L}(t-p)\:\frac{x^2}{b(t)+x^2}\,,
\end{equation}
where $d\ge 0$ is proportional to the size of the group of predators, $\rho>0$ is inversely
related to the average speed at which the predators arrive and leave, $2L\ge 0$ is the period of time
during which the action of the predators is most intense,
and $p\in\R$ fixes the arrival and departure times $p-L-\rho$ and $p+L+\rho$ of the predators.

\begin{figure}
     \centering
         \includegraphics[width=\textwidth]{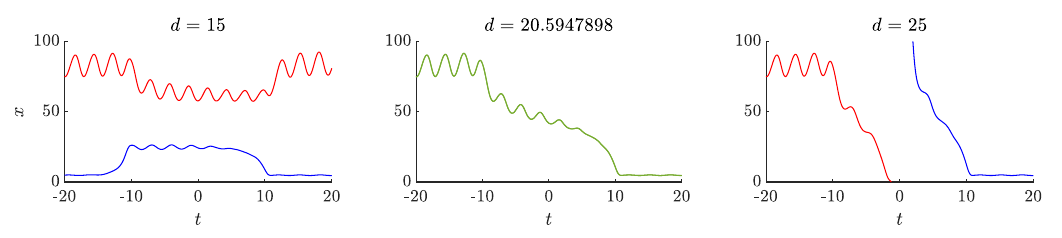}
         \caption{Numerical depiction of the existence of a unique size-tipping point for \eqref{eq:4Cpredcompleta}$_d$
         for $\rho=1$, $L=10$ and $p=0$.
         The central panel depicts the dynamics just before the tipping point $d(1,10,0)$ (see Table~\ref{table:3}):
         the two hyperbolic solutions are so close within the representation window that any of them is a good
         approximation (green) to the unique (nonhyperbolic) bounded solution of \protect\hyperlink{CCB}{\sc Case B}.
         The left panel depicts \protect\hyperlink{CCA}{\sc Case A} (persistence), which is the dynamics for any $d<d(1,10,0)$: the attractive
         hyperbolic solution which represents the behavior of the population in red, and the repulsive one in blue.
         The right panel depicts \protect\hyperlink{CCC}{\sc Case C} (extinction), which is the dynamics for any $d>d(1,10,0)$: the
         locally pullback attractive solution which represents the behavior of the population in red, and the locally pullback repulsive in blue.
         (They are given by Theorem \ref{th:4Cdossoluciones}.)}
         \label{fig:dsingletippingpoint}
\end{figure}
Let us define $f(t,x,\gamma):=r(t)\,x\big(1-x/K(t)\big)+\phi(t)-\gamma\, x^2/(b(t)+x^2)$.
It is easy to check that $(f,d\,\G_{\rho,L}(\cdot-p))$ satisfies
\hyperlink{fc1}{{\rm\textbf{fc1}}}-\hyperlink{fc4}{{\rm\textbf{fc4}}}
with $\Gamma_\pm\equiv 0$ independently of the choice of the parameters.
We emphasize that $f(t,x,d\,\G_{\rho,L}(t-p))$ is not a concave function if
$d>4\,\max_{t\in\R}r(t)\,b(t)/K(t)$, what is easy to check: we are dealing with
an asymptotically concave equation which is not concave.
In addition, we assume that $r$, $K$ and $\phi$ are chosen in such a way that
$x'=f(t,x,0)$ has two (bounded) hyperbolic solutions: \hyperlink{fc5}{\bf fc5} is also fulfilled.
Observe also that we can take $m_1<m_2$ as in
Proposition \ref{prop:3Ccoer} with $m_1>0$: we take $m_1\in (0,\inf_{t\in\R}(-\phi(t)/r(t)))$
and check that $f(t,x,d\,\G_{\rho,L}(t-p))\le f(t,x,0)\leq r(t)\,m_1+\phi(t)<-\delta$ for certain
$\delta>0$ if $x\le m_1$. Hence, all the bounded solutions of \eqref{eq:4Cpredcompleta} are strictly positive
(and hence they are biologically meaningful) for all $d\ge 0$, $\rho>0$, $L\ge 0$ and $p\in\R$.
Besides, $\gamma\mapsto f(t,x,\gamma)$ is nonincreasing for all $(t,x)\in\RR$ and strictly decreasing
for $(t,x)\in\R\times(0,\infty)\supset\R\times[m_1,m_2]$, and $\lim_{\gamma\to\infty} f(t,x,\gamma)=-\infty$
uniformly on $\R\times[m_1,m_2]$. Recall also that
we have chosen the coefficients to ensure that
$x'=f(t,x,0)$ is in \hyperlink{CCA}{\sc Case A}.
So, if we fix $(\rho,L,p)\in(0,\infty)\times[0,\infty)\times\R$ and define
$g(t,x,\gamma):=f(t,x,-\gamma)$, then the pairs
$(g,d\,\G_{\rho,L}(\cdot-p))$ satisfy all the hypotheses of Theorem \ref{th:4Csize2}
(with $\G(t):=0$, $\G_0(t):=\G_{\rho,L}(t-p)$ and $\bar d=0$).
Therefore, there exists a unique tipping point $d(\rho,L,p)>0$:
for $0\le d<d(\rho,L,p)$ (as for $d<0$), the dynamics of \eqref{eq:4Cpredcompleta}$_d$
fits \hyperlink{CCA}{\sc Case A},
and it fits \hyperlink{CCB}{\sc Cases B} and \hyperlink{CCC}{C} for $d=d(\rho,L,p)$ and
for $d>d(\rho,L,p)$, respectively.
In addition, $d(\rho,L,p)$ varies continuously with respect to each parameter,
as Theorem \ref{th:4Csaddle-node} shows.
\par
This critical transition is depicted in Figure \ref{fig:dsingletippingpoint}
for $\rho=1$, $L=10$ and $p=0$, with the next choices:
$r(t):=1+0.2\sin^2(t)$, $K(t):=90+20\sin(\sqrt{5}\,t)$, $\phi(t):=-5$, and $b(t):=20+\cos (t)$.
Numerical evidences show that the corresponding equation $x'=f(t,x,0)$ has two
(strictly positive) hyperbolic solutions, as our analysis requires.
For these choices, the right size of \eqref{eq:4Cpredcompleta}
is not concave in $x$ if $d>1.44$.
\par
In Table~\ref{table:3}, we numerically approximate the
unique bifurcation points $d(1,L,p)$ for some pairs $(L,p)$ and the previous choices.
The bifurcation points have been approximated through bisection methods.
\begin{table}[h!]
\centering
\begin{tabular}{|c | c | c|c|}
 \hline &&&\\[-2.2ex]
 $d(1,L,p)$  & $p=0$ & $p=2$ & $p=5$\\[0.4ex]
\hline &&&\\[-2.2ex]
$L=\mskip4.5mu 1\mskip4.5mu$   &40.2455300 &42.2034404&41.9617506\\[0.4ex]
$L=\mskip4.5mu 5\mskip4.5mu$   &23.0532048 &22.9017928&22.8172667\\[0.4ex]
$L=10$  &20.5947898 &20.5342198&20.4768856\\[0.4ex]
$L=15$  &19.9425668 &19.9151819&19.8875532\\[0.4ex]
$L=20$  &19.6805426 &19.6731947&19.6649049\\[0.4ex]
 \hline
\end{tabular}
\caption{Numerical approximations up to seven decimal places to the bifurcation
point $d(1,L,p)$ of \eqref{eq:4Cpredcompleta}$_d$.
The displayed number is a value of $d$ for which
\eqref{eq:4Cpredcompleta}$_d$ is in \protect\hyperlink{CCA}{\sc Case A} and
\eqref{eq:4Cpredcompleta}$_{d+1e-7}$ is in \protect\hyperlink{CCC}{\sc Case C}.
The Matlab2023a \texttt{ode45} algorithm has been used
with \texttt{AbsTol} and \texttt{RelTol} equal to \texttt{1e-12}.
The final integration has been carried out over the
interval $[-\texttt{1e4},\texttt{1e4}]$.}
\label{table:3}
\end{table} \begin{figure}
     \centering
         \includegraphics[width=\textwidth]{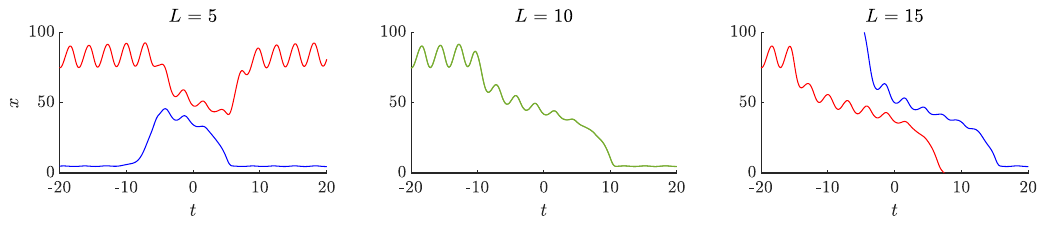}
         \caption{Numerical depiction of a time-of-action-tipping point for \eqref{eq:4Cpredcompleta}$_L$
         for $d=20.5342198$, $\rho=1$ and $p=0$.
         The central panel approximates the dynamics at the tipping point $d(1,10,0)$ (see Table~\ref{table:3}
         and Figure~\ref{fig:dsingletippingpoint}): \protect\hyperlink{CCB}{\sc Case B}.
         The left panel depicts \protect\hyperlink{CCA}{\sc Case A} (persistence) and the
         right panel depicts \protect\hyperlink{CCC}{\sc Case C} (extinction).
         Recall that the predators act during $t\in[-L-\rho,L+\rho]$.}
        \label{fig:dsingletippingpointL}
\end{figure}

Let us briefly extract some conclusions from the data of Table~\ref{table:3}.
In this three-parametric model (we have fixed $\rho=1$), we may find changes from \hyperlink{CCA}{\sc Case A}
to \hyperlink{CCC}{\sc Case C} by varying any of the three parameters $d$, $L$ or $p$ (and $\rho$, although
we work with a fixed value for simplicity).
In fact, since $\G_{\rho,L}(\cdot)$ is nondecreasing and nonconstant with respect to $L$ (and also
with respect to $\rho$), Corollary \ref{coro:4Cunpunto} (applied to the maps $g(t,x,d\,\G_{\rho,L}(t-p))$)
shows that there exists at most a critical value $L_0$ (or $\rho_0$) if the rest
of the parameters are fixed, and that \hyperlink{CCA}{\sc Case A} holds to its left.
Figure \ref{fig:dsingletippingpointL} depicts the occurrence of this critical transition
as $L$ increases for a fixed value of $d$ (which is our previous approximation to $d(1,L,0)$).
(In fact, it can be proved that this tipping point indeed exists if $d$ is large enough and $\rho$
is small enough, but we omit this nontrivial analysis.)
On the other hand, the lack of monotononicity
in the first row of Table~\ref{table:3} shows the possible lack of uniqueness of
the tipping point as $p$ varies: for $L=1$ and $d=42$, we are
in \hyperlink{CCC}{\sc Case C} for $p=0$ and $p=5$ (in those cases $d(1,1,p)<42$), and in
\hyperlink{CCA}{\sc Case A} for $p=2$ (because $d(1,1,2)>42$): there are at least two critical
transitions, i.e., two \hyperlink{CCB}{\sc Cases B}. This fact is depicted in Figure~\ref{fig:doubletippingpoint}.
This lack of uniqueness of phase-induced critical transitions was already observed in
\cite[Section 6]{lno3}.
\par
Summing up, we find a unique size-induced critical transition as $d$ varies
(see Figure \ref{fig:dsingletippingpoint}),
a unique time-of-action-induced critical transition as $L$ varies (see Figure \ref{fig:dsingletippingpointL})
if $d$ is large enough, and possibly several phase-induced critical transitions as $p$ varies
(see Figure~\ref{fig:doubletippingpoint}), also for a large enough $d$.

\begin{figure}
     \centering
         \includegraphics[width=\textwidth]{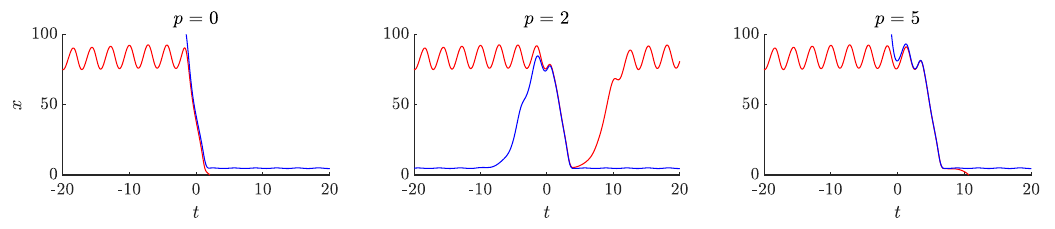}
         \caption{Numerical depiction of the possible existence of two phase-induced tipping points:
         for $L=1$ and $d=42$, we find \protect\hyperlink{CCC}{\sc Case C} (extinction) for $p=0$ and $p=5$, and
         \protect\hyperlink{CCA}{\sc Case A} (persistence) for $p=2$.
         See Figure~\ref{fig:dsingletippingpoint} for the meaning of the different elements.}
        \label{fig:doubletippingpoint}
\end{figure}
\par
Let us finally analyze the occurrence of rate-induced critical transitions, fixing
$d$, $L$, $\rho$ and taking $p=0$. By writing the expression of $\G_{\rho,L}$, we
check that $\G_{\rho,L}(ct)=\G_{\rho/c,L/c}(t)$ for all $\rho,L,c>0$ and $t\in\R$.
This formula relates the study of rate-induced critical transitions to one induced by
a simultaneous change of time-of-action $L$ and the time spent in arrival $\rho$.
It also shows that $c\mapsto \G_{\rho,L}(ct)$ is a nonincreasing and nonconstant map.
Therefore, Corollary~\ref{coro:4Cunpunto} shows the uniqueness of the tipping point in case of existence,
and that \hyperlink{CCA}{\sc Case A} holds to its right: we are in cases of rate-induced tracking.
(And, as in the case of variation of $L$, the existence
of this unique tipping value $c_0$ for each fixed $\rho>0,\,L>0,\,p\in\R$ and large enough $d$
can be indeed proved.)
\end{exa}
\begin{exa}
The analysis in Example~\ref{ex:4Cejemplo} has been carried out considering an
$x$-independent transition function $t\mapsto d\, \G_{\rho,L}(t-p)$.
Now we add another parameter shift $t\mapsto\Lambda(t-p)$ in the migration part
of the driving law. The results presented in this paper can still be applied
to an $x$-dependent transition function $(t,x)\mapsto\Delta_{d,\rho,L,p}(t,x)$
(recalling the formulation of \eqref{eq:4Cpredgeneral}) that encompasses the
above two parameter shifts. The presented theory allows to understand separately
a stationary part and a transient part of the law that generates the dynamics.
\par
So, we consider a slightly more complicated emigration phenomena: due to the arrival
to the habitat of the predator species, the attractiveness of the habitat is reduced,
and, as a consequence, the emigration increases.
To model this, we consider that $r,\,K,\,b$ and $\phi$ satisfy the conditions assumed in
Example \ref{ex:4Cejemplo}, and replace the quasiperiodic emigration function
$\phi$ by an asymptotically quasiperiodic emigration function
$\phi(t)+(\psi(t)-\phi(t))\,\Lambda_{\rho,L}(t-p)$, where $\psi$ is quasiperiodic
with $\psi(t)\leq\phi(t)$ for all $t\in\R$, and $\Lambda_{\rho,L}(t)$ is the unique $C^1$
cubic piecewise polynomial which takes values $0$ for all $x\leq -L-\rho$ and $1$ for all
$x\geq -L$. So, we work with
\begin{equation}\label{eq:4CpredcompletaCambioEmigracion}
 x'=r(t)\,x\left(1-\frac{x}{K(t)}\right)+\phi(t)+(\psi(t)-
 \phi(t))\,\Lambda_{\rho,L}(t-p)-\frac{d\,\G_{\rho,L}(t-p)\, x^2}{b(t)+x^2}\,.
\end{equation}
If, as before, $h(t,x,\delta):=r(t)\,x\,(1-x/K(t))+\delta$, and
$\Delta_{d,\rho,L,p}(t,x):=\phi(t)+(\psi(t)-\phi(t))\,
\Lambda_{\rho,L}(t-p)-d\,\Gamma_{\rho,L}(t-p)\, x^2/(b(t)+x^2)$, then
$(h,\Delta_{d,\rho,L,p})$ satisfies \hyperlink{fc1}{\bf fc1}-\hyperlink{fc4}{\bf fc4}
with $\Delta_-:=\phi$ and $\Delta_+:=\psi$. We also assume that $x'=h(t,x,\psi(t))$
has two hyperbolic solutions and deduce from Proposition \ref{prop:3Ccoer}(v)
that this property also holds for $x'=h(t,x,\phi(t))$: \hyperlink{fc5}{\bf fc5} is also
fulfilled. Proposition \ref{prop:3Ccoer}(v) also shows that
$x'=h(t,x,\phi(t)+(\psi(t)-\phi(t))\,\Lambda_{\rho,L}(t-p))$ has two hyperbolic solutions
for all $L\ge 0$ and $\rho>0$. This fact allows us to repeat the arguments of
Example \ref{ex:4Cejemplo} in order to prove the existence and uniqueness of the
bifurcation point $d(\rho,L,p)>0$, giving rise to a unique size-induced critical transition.

\begin{table}[h!]
\centering
\begin{tabular}{|c | c | c|c|}
 \hline &&&\\[-2.2ex]
 $d(1,L,p)$  & $p=0$ & $p=2$ & $p=5$\\[0.4ex]
\hline &&&\\[-2.2ex]
$L=\mskip4.5mu 1\mskip4.5mu$  &34.1938684 &36.9449750 &35.7506039\\[0.4ex]
$L=\mskip4.5mu 5\mskip4.5mu$  &18.8506812 &18.6486286 &18.6059557\\[0.4ex]
$L=10$ &16.4930418 &16.4318568 &16.3869395\\[0.4ex]
$L=15$ &15.8700118 &15.8460071 &15.8202938\\[0.4ex]
$L=20$ &15.6203137 &15.6150934 &15.6065018\\[0.4ex]
 \hline
\end{tabular}
\caption{Numerical approximations up to seven decimal places to the bifurcation
point $d(1,L,p)$ of \eqref{eq:4CpredcompletaCambioEmigracion}$_d$.
The displayed number is a value of $d$ for which
\eqref{eq:4CpredcompletaCambioEmigracion}$_d$ is in \protect\hyperlink{CCA}{\sc Case A} and
\eqref{eq:4CpredcompletaCambioEmigracion}$_{d+1e-7}$ is in \protect\hyperlink{CCC}{\sc Case C}.
The Matlab2023a \texttt{ode45} algorithm has been used
with \texttt{AbsTol} and \texttt{RelTol} equal to \texttt{1e-12}.
The final integration has been carried out over the
interval $[-\texttt{1e4},\texttt{1e4}]$.}
\label{table:3_V2}
\end{table}
We repeat the choices of $r$, $K$, $\phi$ and $b$
of the first example, and take $\psi(t):=-9-\cos t$. Again, numerical evidences show
that $x'=h(t,x,\psi(t))$ has two (positive) hyperbolic solutions, as required.
Table~\ref{table:3_V2} shows numerical approximations of the
unique bifurcation point $d(1,L,p)$ for some pairs $(L,p)$.
The arguments of Example \ref{ex:4Cejemplo} also work to prove the existence
of a unique bifurcation value $L_0>0$ if $d$ is large enough, $\rho>0$ is small
enough, and $p\in\R$.
Note that, for $L=1$, we find the phenomenon of multiple phase bifurcation points
already mentioned in Table~\ref{table:3} and shown in Figure~\ref{fig:doubletippingpoint}.
As expected, the bifurcation points $d(1,L,p)$ for \eqref{eq:4CpredcompletaCambioEmigracion}
in Table~\ref{table:3_V2} are lower than those for \eqref{eq:4Cpredcompleta} in Table~\ref{table:3}:
an increased emigration means that fewer predators are needed to cause the population extinction.
We point out that, in this example, the asymptotic limits of the transition
map are time-dependent.
\end{exa}
\section{The d-concave and nonquadratic case}\label{sec:5}
The results of this section partly extend those of \cite{dno1} to a less
restrictive setting. Let $(\W,\sigma)$ be defined as in Section \ref{sec:3}.
Now, we work with the family
\begin{equation}\label{eq:5ini}
 x'=\mh(\wt,x)
\end{equation}
and with the flow $\tau$ defined by \eqref{def:2tau}, $\tau(t,\w,x)=(\wt,v(t,\w,x))$,
assuming that $\mh\colon\WR\to\R$ satisfies (all or part of) the next conditions:
\begin{list}{}{\leftmargin 17pt}
\item[\hypertarget{d1}{{\bf d1}}] $\mh\in C^{0,2}(\WR,\R)$,
\item[\hypertarget{d2}{{\bf d2}}] $\limsup_{x\to\pm\infty} (\pm\mh(\w,x))<0$ uniformly on $\W$,
\item[\hypertarget{d3}{{\bf d3}}] $m(\{\w\in\W\,|\;x\mapsto \mh_x(\w,x) \text{ is concave}\})=1$
for all $m\in\merg$,
\item[\hypertarget{d4}{{\bf d4}}] $m(\{\w\in\W\,|\; x\mapsto \mh_{xx}(\w,x)$ is strictly decreasing
on $\mJ\})>0$ for all compact interval $\mJ\subset\R$ and all $m\in\merg$.
\end{list}
We will refer to the concavity of the derivative as {\em d-concavity\/}.
\begin{nota}\label{rm:5DCtambien}
As explained in Remark \ref{rm:3Ctambien}, if $\mh$ satisfies
\hyperlink{d1}{\bf dj} for {\bf j} $\in\{1,2,3,4\}$ and $\W_0\subset\W$ is
a nonempty compact $\sigma$-invariant subset, then also the restriction
$\mh\colon\W_0\times\R\to\R$ satisfies \hyperlink{d1}{\bf dj}.
\end{nota}
\par
The following results establish properties analogous to those of
Theorems~\ref{th:3C2medidas} and \ref{th:3C2copias}, now for the new hypotheses on d-concavity.
\begin{teor} \label{th:5DC3medidas}
Let $\mh\colon\WR\to\R$ satisfy \hyperlink{d1}{\bf d1}, let us fix $m\in\merg$,
and let $\mb_1,\mb_2,\mb_3\colon\W\to\R$ be
bounded $m$-measurable $\tau$-equilibria with $\mb_1(\w)<\mb_2(\w)<\mb_3(\w)$ for $m$-a.e.~$\w\in\W$.
Assume that $m(\{\w\in\W\,|\;x\mapsto \mh_x(\w,x) \text{ is concave}\})=1$
and $m(\{\w\in\W\,|\;\mh_{xx}(\w,\mb_1(\w))>\mh_{xx}(\w,\mb_3(\w))\})>0$. Then,
\[
\int_\W \mh_x(\w,\mb_2(\w))\, dm>0\quad\text{and}\quad \int_\W \mh_x(\w,\mb_i(\w))\, dm<0 \text{ \;for\: $i=1,3$}\,.
\]
In particular, there are at most three
bounded $m$-measurable $\tau$-equilibria which are strictly ordered $m$-a.e.
\end{teor}
\begin{proof}
We call $\W^c:=\{\w\in\W\,|\;x\mapsto \mh_x(\w,x)\text{ is concave}\}$, which satisfies
$m(\W^c)>0$, and $\W_0:=\{\w\in\W\,|\;\mb_1(\w)<\mb_2(\w)<\mb_3(\w)\}$, which is invariant and with $m(\W_0)=1$.
For each $\w\in\W^c$, we represent by $b_i(\w,x_1,x_2,x_3)$ the expression $b_i(x_1,x_2,x_3)$ of \eqref{def:2ai}
associated to the d-concave map
$x\mapsto \mh(\w,x)$. For $i\in\{1,2,3\}$, we define $b_i^*\colon\W\to\R$ by $b_i^*(\w):=
b_i(\w,\mb_1(\w),\mb_2(\w),\mb_3(\w))$ for $\w\in\W_0\cap\W^c$ and
$b_i^*(\w):=0$ if $\w\notin\W_0\cap\W^c$. The hypothesis on $\mh_{xx}$
and Theorem~\ref{th:2positiveconcave}(ii) (see also its proof)
ensure that $b_i^*$ is $m$-measurable, with $b_i^*\ge 0$ and
$m(\{\w\in\W\,|\; b_i^*(\w)>0\})>0$. These are the key properties to repeat
the arguments of \cite[Theorem~4.1]{dno1} in order to check that, if
$v_1:=1/(\mb_2-\mb_1)-1/(\mb_3-\mb_1)$ (which satisfies $v_1(\w)>0$
for $\w\in\W_0$), then $\int_\W \mh_x(\w,\mb_1(\w))\,dm=
-\int_\W(b_1^*(\w)/v_1(\w))\,dm<0$. The argument is similar for $\mb_2$ and $\mb_3$,
and the last assertion follows from the previous ones and
a simple contradiction argument.
\end{proof}
\begin{teor}\label{th:4DC3copias}
Let $\mh$ satisfy \hyperlink{d1}{\bf d1}, \hyperlink{d3}{\bf d3} and \hyperlink{d4}{\bf d4}. Then,
there exist three disjoint and ordered $\tau$-invariant compact sets $\mK_1<\mK_2<\mK_3$
projecting onto $\W$ if and only if there exist three hyperbolic copies
of the base $\{\tml\}$, $\{\tmm\}$ and $\{\tmuk\}$ with $\tml<\tmm<\tmuk$.
In this case, $\mK_1=\{\tml\}$ and $\mK_3=\{\tmuk\}$ and they are attractive;
$\mK_2=\{\tmm\}$ and it is repulsive; and
$\mB:=\{(\w,x)\in\WR\,|\;\tml(\w)\le x\le\tmuk(\w)\}$ is the set of globally bounded orbits.
In particular, there are at most three disjoint and ordered $\tau$-invariant compact sets
projecting onto $\W$.
\end{teor}
\begin{proof}
Observe that $\mh$ satisfies the conditions of
Theorem \ref{th:5DC3medidas} for any $m\in\merg$ and any three bounded
ordered $m$-measurable equilibria. This fact allows us to use the arguments
of the proof of Theorem~\ref{th:3C2copias} to check all the assertions.
\end{proof}
\begin{nota}
An analogue of Remark \ref{rm:3basta} applies to Theorem \ref{th:4DC3copias}.
\end{nota}
\par
Now, the coercivity property \hyperlink{d2}{\bf d2} comes into play.
Recall that a $\tau$-invariant compact set $\mA\subset\WR$ is the {\em global attractor of\/}
$\tau$ if it attracts every bounded set $\mC\subset\WR$;
that is, if $\tau_t(\mC)$ is defined for any $t\ge 0$ and
$\lim_{t\to\infty} \text{dist}(\tau_t(\mC),\mA)=0$, where
$\text{dist}(\mC_1,\mC_2)$ is the Hausdorff semidistance from $\mC_1$ to $\mC_2$ and
$\tau_t(\mC):=\{\tau(t,\w,x)\,|\;(\w,x)\in\mC\}$.
\begin{prop}\label{prop:5Dcoer}
Let $\mh\in C^{0,1}(\WR,\R)$ satisfy \hyperlink{d2}{\bf d2}, and take $\delta>0$ and
$m_1,m_2\in\R$ with $\mh(\w,x)\ge \delta$ if $x\le m_1$ and
$\mh(\w,x)\le -\delta$ if $x\ge m_2$ for all $\w\in\W$. Then,
\begin{itemize}
\item[\rm (i)] $v(t,\w,x)$ exists for $(t,\w,x)\in[0,\infty)\times\WR$, and
$m_1\le\liminf_{t\to\infty}v(t,\w,x)\le\limsup_{t\to\infty}v(t,\w,x)\le m_2$:
any forward $\tau$-semiorbit is bounded.
\item[\rm(ii)] There exists the global attractor for $\tau$, it is of the form
\begin{equation}\label{eq:5DdefA}
 \mA=\bigcup_{\w\in\W}\big(\{\w\}\times[\ml(\w),\muk(\w)]\big)\,,
\end{equation}
it is the union of all the globally defined and bounded $\tau$-orbits,
and it is contained in $\W\times[m_1,m_2]$.
\item[\rm(iii)] The maps $\ml$ and $\muk$ are, respectively, lower and
    upper semicontinuous $\tau$-equilibria.
\item[\rm (iv)] If, for a point $\w\in\W$, there exists a bounded $C^1$ function
    $b\colon\R\to\R$ such that $b'(t)\le \mh(\wt,b(t))$ (resp.~$b'(t)\ge
    \mh(\wt,b(t))$) for all $t\in\R$, then $b(t)\le\muk(\wt)$ (resp.~$b(t)\ge\ml(\wt)$)
    for all $t\in\R$. If $b'(t)<\mh(\wt,b(t))$ (resp. $b'(t)>\mh(\wt,b(t))$) for all $t\in\R$,
    then $b(t)<\muk(\wt)$ (resp.~$\ml(\wt)<b(t)$) for all $t\in\R$.
\item[\rm(v)] $v(t,\w,x)$ is bounded from below if and only if $x\ge\ml(\w)$, and
    from above if and only if $x\le\muk(\w)$.
\item[\rm (vi)] Assume that $\mh$ satisfies also \hyperlink{d1}{\bf d1}, \hyperlink{d3}{\bf d3} and \hyperlink{d4}{\bf d4},
    and that $\{\tml\}$, $\{\tmm\}$ and $\{\tmuk\}$ are three hyperbolic copies of the base with
    $\tml<\tmm<\tmuk$. Then, $\mA=\bigcup_{\w\in\W}
    (\{\w\}\times [\tml(\w),\tmuk(\w)])$; $\{\tml\}$ and
    $\{\tmuk\}$ are attractive and $\{\mm\}$ is repulsive;
    $\lim_{t\to\infty}(v(t,\w,x)-\tmuk(\wt))=0$ if and only if $x>\tmm(\w)$;
    $\lim_{t\to\infty}(v(t,\w,x)-\tml(\wt))=0$ if and only if $x<\tmm(\w)$;
    and $\lim_{t\to-\infty}(v(t,\w,x)-\tmm(\wt))=0$ if and only if $x\in(\tml(\w),
    \tmuk(\w))$.
\end{itemize}
\end{prop}
\begin{proof} The existence of $m_1$ and $m_2$ is ensured by \hyperlink{d2}{\bf d2}.
The properties stated in (i) are a nice exercise on ODEs.
To prove (ii), we take $n_1<m_1$ and $n_2>m_2$ and check that $v(t,\w,n_i)\in[m_1,m_2]$
for all $\w\in\W$ and $i=1,2$ if $t\ge(1/\delta)\max(m_1-n_1,n_2-m_2)$.
We deduce from this fact that
$\lim_{t\to\infty}\text{dist}(\tau_t(\mC),\W\times[m_1,m_2])=0$
for every bounded set $\mC\subset\WR$; i.e.,
$\W\times[m_1,m_2]$ is a compact absorbing set. This property and \cite[Theorem 2.2]{chks}
prove the existence of the global attractor $\mA\subseteq\W\times[m_1,m_2]$, and
\cite[Theorem 1.7]{carvalho1} ensures the last assertion in (ii).
\par
Assertion (iii) is a consequence of the compactness of $\mA$; and the properties in (iv)
follow from (i) and standard comparison results: see, e.g., the proof of \cite[Theorem 5.1(iii)]{dno1}.
The assertions in (v) follow from (i) and (ii). In the conditions of (vi),
Theorem \ref{th:4DC3copias} shows that
$\mA=\bigcup_{\w\in\W}(\{\w\}\times [\tml(\w),\tmuk(\w)])$, and that $\{\tml\}$ and $\{\tmuk\}$
are attractive and $\{\tmm\}$ is repulsive. The remaining assertions are proved with the
arguments used to check Proposition \ref{prop:3Ccoer}(vii),
working with the \upomeg-limit sets in the cases of $x>\tmm(\w)$ and $x<\tmm(\w)$, and
with the \upalfa-limit set for $\tml(\w)<x<\tmuk(\w)$.
\end{proof}
It follows from the previous property (iii) that $\ml$ and $\muk$ are $m$-measurable
equilibria for any $m\in\merg$, which we will use without further reference.
In the line of Theorem \ref{th:3Cequiv}, our next result establishes equivalences
regarding the existence of three uniformly separated hyperbolic solutions of a
given equation in terms of the existence of three ordered hyperbolic copies of the
corresponding hull.
\begin{teor}\label{th:5DCthreecopies}
Let $\mh\colon\WR\to\R$ satisfy \hyperlink{d1}{\bf d1}, \hyperlink{d2}{\bf d2},
\hyperlink{d3}{\bf d3} and \hyperlink{d4}{\bf d4}. Let us fix $\bar\w\in\W$.
Then, the following assertions are equivalent:
\begin{itemize}
\item[\rm(a)] Equation \eqref{eq:5ini}$_{\bar\w}$ has three hyperbolic solutions.
\item[\rm(b)] Equation \eqref{eq:5ini}$_{\bar\w}$ has three uniformly separated
hyperbolic solutions.
\item[\rm(c)] Equation \eqref{eq:5ini}$_{\bar\w}$ has three uniformly separated
bounded solutions.
\item[\rm(d)] There exist three hyperbolic copies of the base for the restriction of
the family \eqref{eq:5ini} to the closure $\W_{\bar\w}$ of $\{\bwt\,|\;t\in\R\}$,
given by $\tml<\tmm<\tmuk$.
\end{itemize}
In this case, $t\mapsto \tilde l(t):=\tml(\bwt)$, $t\mapsto\tilde m(t):=\tmm(\bwt)$ and
$t\mapsto\tilde u(t):=\tmuk(\bwt)$ are the three unique uniformly separated solutions of
\eqref{eq:5ini}$_{\bar\w}$, they are hyperbolic, and there are no more hyperbolic solutions.
In addition, if $x_{\bar\w}(t,s,x)$ is the solution of \eqref{eq:5ini}$_{\bar\w}$ with
$x_{\bar\w}(s,s,x)=x$, then:
$\lim_{t\to\infty}(x_{\bar\w}(t,s,x)-\tilde u(t))=0$ if and only if $x>\tilde m(s)$;
$\lim_{t\to\infty}(x_{\bar\w}(t,s,x)-\tilde l(t))=0$ if and only if $x<\tilde m(s)$;
and $\lim_{t\to-\infty}(x_{\bar\w}(t,s,x)-\tilde m(t))=0$ if and only if $x\in(\tilde l(s),
\tilde u(s))$.
\end{teor}
\begin{proof}
The proof of this result follows the line of that Theorem \ref{th:3Cequiv}.
The assertions after the equivalences follow from (d) and Proposition \ref{prop:5Dcoer}(v)
and (vi).
We will check (b)$\,\Rightarrow\,$(c)$\,\Rightarrow\,$(d)$\,\Rightarrow\,$(a)$\,\Rightarrow\,$(b).
Recall that the hypotheses on $\mh$ are also valid for its restriction to $\W_{\bar\w}\times\R$:
see Remark \ref{rm:5DCtambien}.
\smallskip\par
(b)$\,\Rightarrow\,$(c)$\,\Rightarrow\,$(d). Obviously, (b) implies (c).
Now we assume (c) and observe that there is no restriction
in assuming that the three uniformly separated solutions are $l(t):=\ml(\bwt)$,
$u(t):=\muk(\bwt)$, and $m(t)$ with $l(t)<m(t)<u(t)$.
We will check that the closures $\mK_l$, $\mK_m$ and $\mK_u$ of the corresponding $\tau$-orbits
are three different ordered compact sets projecting on $\W_{\bar\w}$.
Reasoning as in Theorem \ref{th:3Cequiv}, we check that, for a $\delta>0$,
$x_0\le\muk(\w_0)-\delta$ and $x_0\ge\ml(\w_0)+\delta$
whenever $(\w_0,x_0)\in\mK_m$. Theorem \ref{th:5DC3medidas} allows us to
assert that all the Lyapunov exponents of $\mK_m$ are positive, and that
the upper and lower equilibria of $\mK_m$ coincide on a set $\W_0$ with
$m_0(\W_0)=1$ for all $m_0\in\mathfrak{M}_\mathrm{erg}(\W_{\bar\w},\sigma)$.
Theorem \ref{th:2copia} ensures that
$\mK_m$ is a repulsive hyperbolic copy of the base $\W_{\bar\w}$,
say $\mK_m=\{\tmm\}$. Now, it is easy to deduce from the previous separation properties that
$\mK_l<\mK_m<\mK_u$, as asserted. Theorem \ref{th:4DC3copias} applied to $\W_{\bar\w}\times\R$
shows that (d) holds.
\smallskip\par
(d)$\,\Rightarrow\,$(a)$\,\Rightarrow\,$(b). If (d) holds, then
$t\mapsto\tml(\bwt)$, $t\mapsto\tmm(\bwt)$ and
$t\mapsto\tmuk(\bwt)$ are three hyperbolic solutions of \eqref{eq:5ini}$_{\bar\w}$
(see Proposition \ref{prop:2extiende}),
so (a) holds. Let us assume (a), and let $\tilde x_1<\tilde x_2<\tilde x_3$
be the three hyperbolic solutions of \eqref{eq:5ini}$_{\bar\w}$.
Let us first eliminate the possibility that $\tilde x_2$ is
attractive, assuming it for contradiction. Let us call
$l(t):=\ml(\bwt)$ and $u(t):=\muk(\bwt)$, where $\ml$ and $\muk$
are the lower and upper $\tau$-equilibria: see \eqref{eq:5DdefA}.
Proposition \ref{prop:2haciaatras}(i)
yields $\delta>0$ such that $\inf_{t\le 0}(u(t)-\tilde x_2(t))>\delta$
and $\inf_{t\le 0}(\tilde x_2(t)-l(t))>\delta$. Let $\mM_l$,
$\mM_2$ and $\mM_u$ be the \upalfa-limit sets of $(\bar\w,l(0))$,
$(\bar\w,\tilde x_2(0))$ and $(\bar\w,u(0))$, which project on the \upalfa-limit
set $\W_-\subseteq\W_{\bar\w}$ of $\bar\w$. As in the previous paragraph, we check that
$\ml(\w_0)+\delta\le x_0\le\muk(\w_0)-\delta$ if $(\w_0,x_0)\in\mM_2$;
and deduce that $\mM_2$ is a repulsive copy of $\W_-$, which
contradicts Proposition \ref{prop:2extiende}.
\par
Hence, $\tilde x_2$ is repulsive. Proposition \ref{prop:2haciaatras}(i)
yields $\delta>0$ such that $\inf_{t\ge 0}(u(t)-\tilde x_2(t))>\delta$
and $\inf_{t\ge 0}(\tilde x_2(t)-l(t))>\delta$. Let
$\bar\mM_2$ be the \upomeg-limit set of $(\bar\w,\tilde x_2(0))$,
which projects on the \upomeg-limit set $\W_+\subseteq\W_{\bar\w}$ of $\bar\w$.
As in the proof of (c)$\,\Rightarrow\,$(d), we check that
$\ml(\w_0)+\delta\le x_0\le\muk(\w_0)-\delta$ whenever $(\w_0,x_0)\in\bar\mM_2$;
and we deduce that $\bar\mM_2$ is a repulsive copy of $\W_+$. Hence,
$\bar\mM_2$ does not intersect the \upomeg-limit sets $\bar\mM_1$ of
$(\bar\w,\tilde x_1(0))$ and $\bar\mM_3$ of $(\bar\w,\tilde x_3(0))$:
see Proposition \ref{prop:2haciaatras}(ii).
So, we have $\bar\mM_1<\bar\mM_2<\bar\mM_3$. Theorem \ref{th:4DC3copias} applied to $\W_+\times\R$
ensures that $\bar\mM_1$ and $\bar\mM_3$ are attractive hyperbolic copies of $\W_+$,
which according to Proposition \ref{prop:2extiende} is only possible
if $\tilde x_1$ and $\tilde x_3$ are attractive. Proposition \ref{prop:2haciaatras}(i)
ensures that the three solutions are uniformly separated: (b) holds.
\end{proof}
In Subsection \ref{subsec:61}, we will analyze a situation precluding the occurrence of
critical transitions. Some of the hypotheses will refer to the relative order of
the three hyperbolic copies of the base of two ``ordered" equations
which are in the situation of Theorem \ref{th:4DC3copias}. For the sake of completeness,
Theorem \ref{th:5DCorder} proves that there are just two relative positions in the case of a
minimal base under the hypotheses so far assumed. In Subsection \ref{subsec:ordenes},
we will check that the minimality of the base is indeed required for Theorem \ref{th:5DCorder},
whose proof is based on the next result.
\begin{prop}\label{prop:5DCorder}
Let $\mh_0,\mh_1\colon\WR\to\R$ satisfy
\hyperlink{d1}{\bf d1} and \hyperlink{d2}{\bf d2}, with $\mh_0(\w,x)\leq \mh_1(\w,x)$ for all $(\w,x)\in\WR$.
\begin{itemize}
\item[\rm(i)] Let $\ml_i$ (resp.~$\muk_i$) be the lower (resp.~upper) bounds of the global attractor
    \eqref{eq:5DdefA} of $x'=\mh_i(\wt,x)$ for $i=0,1$. Then, $\ml_0\leq \ml_1$ and $\muk_0\leq \muk_1$.
\item[\rm(ii)] Assume also that $\mh_0$ and $\mh_1$ satisfy \hyperlink{d3}{\bf d3} and \hyperlink{d4}{\bf d4},
    and that there exists $\bar\w\in\W$ such that $x'=\mh_i(\bwt,x)$ has three hyperbolic solutions
    $\tilde l_i<\tilde m_i<\tilde u_i$ for $i=0,1$. Then, $\inf_{t\in\R}(\tilde m_0(t)-\tilde l_1(t))>0$ if and only
    if $\inf_{t\in\R}(\tilde u_0(t)-\tilde m_1(t))>0$, in which case $\tilde l_0\leq \tilde l_1<\tilde m_1\leq
    \tilde m_0<\tilde u_0\leq \tilde u_1$. If, in addition, $\mh_0(\bar\w{\cdot}t,x)< \mh_1(\bar\w{\cdot}t,x)$ for all
    $(t,x)\in\RR$, then all the inequalities are strict.
\end{itemize}
\end{prop}
\begin{proof}
The proof repeats that of \cite[Proposition 4.2]{dno3}, which is
generalized by this result. We must just use
Theorem \ref{th:5DCthreecopies} instead of \cite[Theorem 3.3]{dno3}.
\end{proof}
\begin{teor}\label{th:5DCorder}
Assume that $(\W,\sigma)$ is minimal.
Let $\mh_0,\mh_1\colon\WR\to\R$ satisfy \hyperlink{d1}{\bf d1},
\hyperlink{d2}{\bf d2}, \hyperlink{d3}{\bf d3} and \hyperlink{d4}{\bf d4}, and with
$\mh_0(\w,x)<\mh_1(\w,x)$ for all $(\w,x)\in\WR$.
Assume that the family $x'=\mh_i(\wt,x)$ has three hyperbolic copies of the base $\tml_i<\tmm_i<\tmuk_i$
for $i=0,1$. Then, one of the following orders holds:
\begin{itemize}
\item[{\rm (1)}] $\tml_0<\tml_1<\tmm_1<\tmm_0<\tmuk_0<\tmuk_1$,
\item[{\rm (2)}] $\tml_0<\tmm_0<\tmuk_0<\tml_1<\tmm_1<\tmuk_1$.
\end{itemize}
\end{teor}
\begin{proof} Let us assume (2) does not hold. Hence, there exists
$\bar\w\in\W$ with $\tml_1(\bar\w)\le\tmuk_0(\bar\w)$. A comparison argument ensures that
$\tml_1(\bwt)<\tmuk_0(\bwt)$ for all $t<0$, and hence the minimality of $\W$ ensures that
$\tml_1\le\tmuk_0$. If, in addition, $\tml_1(\w_0)=\tmuk_0(\w_0)$ for an $\w_0\in\W$, then
a new comparison argument shows that $\tmuk_0(\w_0{\cdot}t)<\tml_1(\w_0{\cdot}t)$ for all $t>0$,
impossible. Therefore, $\tml_1<\tmuk_0$.
Now, for contradiction, we assume that neither (1) holds, and deduce
from Proposition~\ref{prop:5DCorder}(ii) and the minimality of $\W$ the existence of
$\bar\w\in\W$ with $\tmm_0(\bar\w)\leq\tml_1(\bar\w)$. Hence,
$\tmm_0(\bwt)<\tml_1(\bwt)$ for all $t>0$. We fix $t_0>0$ and call
$\w_0:=\bar\w{\cdot}t_0$. Proposition
\ref{prop:5Dcoer}(v) yields $\lim_{t\to\infty}(\tmuk_0(\w_0{\cdot}t)-v_0(t,\w_0,\tml_1(\w_0)))=0$,
where $v_0$ stands for the solutions of $x'=\mh_0(\wt,x)$.
Since $v_0(t,\w_0,\tml_1(\w_0))<\tml_1(\w_0{\cdot}t)$ for all $t>0$,
we deduce that $\limsup_{t\to\infty}(\tmuk_0(\w_0{\cdot}t)-\tml_1(\w_0{\cdot}t)))\le0$,
which combined with the minimality of $\W$ contradicts $\tml_1<\tmuk_0$.
\end{proof}
We point out that, under the hypotheses of Theorem \ref{th:5DCorder}, part of
the arguments of \cite[Section 5]{dno1} show that
the situation (1) is equivalent to the absence of a bifurcation value
$\lambda_0\in[0,1]$ for the family $x'=\mh_0(\wt,x)+\lambda(\mh_1(\wt,x)-\mh_0(\wt,x))$.
\section{Asymptotically d-concave transition equations}\label{sec:6}
Let $g\colon\RR\to\R$ be a $C^2$-admissible function.
As explained at the beginning of Section \ref{sec:4}, we can understand
\begin{equation}\label{eq:6Dgtran}
 x'=g(t,x)
\end{equation}
as a transition between the corresponding \upalfa-family and \upomeg-family.
Now, we will assume that these limit families satisfy conditions
\hyperlink{d1}{\bf d1}-\hyperlink{d4}{\bf d4}, as well as the existence of
three hyperbolic copies of the base for the \upalfa-family and the \upomeg-family.
These last conditions are those which provide a wider
range of dynamical possibilities for \eqref{eq:6Dgtran}
under conditions \hyperlink{d1}{\bf d1}-\hyperlink{d4}{\bf d4}:
the maximum number of uniformly separated solutions for each
equation of the \upalfa-family or the \upomeg-family is three: see
Theorem \ref{th:5DCthreecopies}; hence,
Proposition \ref{prop:2solhull} precludes the existence of more than
three uniformly separated solutions of
\eqref{eq:6Dgtran}; and, if there are three,
then Theorem \ref{th:5DCthreecopies} yields three hyperbolic copies of
the base for the \upalfa-family and the \upomeg-family.
\par
As in Section \ref{sec:4}, we will achieve the required properties
by assuming the existence of strictly
d-concave (in $x$) maps $g_-$ and $g_+$ forming asymptotic pairs
with $g$. More precisely, we fix $g$ and assume the existence of
$g_-$ and $g_+$ such that
\begin{list}{}{\leftmargin 23pt}
\item[\hypertarget{gd1}{\bf gd1}] $g,g_-,g_+\in C^{0,2}(\RR,\R)$.
\item[\hypertarget{gd2}{\bf gd2}] $\lim_{t\to\pm\infty}(g(t,x)-g_\pm(t,x))=0$
uniformly on each compact subset $\mJ\subset\R$.
\item[\hypertarget{gd3}{\bf gd3}] $\limsup_{x\to\pm\infty}(\pm h(t,x))<0$
uniformly on $\R$ for $h=g,g_-,g_+$.
\item[\hypertarget{gd4}{\bf gd4}] $\inf_{t\in\R}\big((g_\pm)_{xx}(t,x_1)-(g_\pm)_{xx}(t,x_2)\big)>0$
whenever $x_1<x_2$.
\item[\hypertarget{gd5}{\bf gd5}] Each one of the equations
\begin{equation}\label{eq:6Dglim}
 x'=g_-(t,x) \quad\text{and}\quad x'=g_+(t,x)
\end{equation}
has three hyperbolic solutions, $\tilde l_{g_-}<\tilde m_{g_-}<\tilde u_{g_-}$ and
$\tilde l_{g_+}<\tilde m_{g_+}<\tilde u_{g_+}$.
\end{list}
\par
As in Section \ref{sec:4} (see Remarks \ref{rm:4Chipos}), we will say that
``$g$ satisfies conditions \hyperlink{gd1}{\bf gd1}-\hyperlink{gd5}{\bf gd5}"
if there exist $g_-$ and $g_+$ such that all the listed conditions are satisfied, and
we will refer to the first and second equations in \eqref{eq:6Dglim} as
the {\em past\/} and {\em future equations} of the {\em transition equation\/}
\eqref{eq:6Dgtran}. Some of the results of this section extend part of those of \cite{dno3} to a much
more general setting: the setting and hypotheses of this section are
considerably less restrictive than those leading to the analogous
results in \cite{dno3}.
\par
Our initial purpose is to classify the dynamical scenarios for the
transition equation \eqref{eq:6Dgtran} when $g$ satisfies
\hyperlink{gd1}{\bf gd1}-\hyperlink{gd5}{\bf gd5},
which is achieved in Theorem~\ref{th:6Dcasos}. Its proof is based on some previous
results. The notation established before Lemma \ref{lema:4Cunion} is used in what follows.
The proof of the next result, which uses Lemma \ref{lema:4Cunion},
is almost identical to that of Lemma
\ref{lema:4Chiposc}. The only difference is that we must work with the second
derivative in the last step of the proof.
\begin{lema}\label{lema:6Dhiposd}
If $h\in C^{0,2}(\RR,\R)$ then $\mh$ satisfies \hyperlink{d1}{\bf d1} on $\W_h$.
If $h\in C^{0,2}(\RR,\R)$ and $\limsup_{x\to\pm\infty}(\pm h(t,x))<0$
uniformly on $\R$, then $\mh$ satisfies \hyperlink{d2}{\bf d2} on $\W_h$.
And, if \hyperlink{gd1}{\bf gd1}, \hyperlink{gd2}{\bf gd2} and \hyperlink{gd4}{\bf gd4} hold, then
$\mg$ and $\mg_\pm$ satisfy \hyperlink{d3}{\bf d3} and \hyperlink{d4}{\bf d4} on $\W_g$ and $\W_{g_\pm}$, respectively.
\end{lema}
Assume that $h\in C^{0,2}(\RR,\R)$ and $\limsup_{x\to\pm\infty}(\pm h(t,x))<0$.
Lemma \ref{lema:6Dhiposd} and Proposition \ref{prop:5Dcoer} ensure the existence of
the global attractor $\mA_h=\bigcup_{\w\in\W_h}\big(\{\w\}\times[\ml_h(\w),\muk_h(\w)]\big)$ of the
flow $\tau_h$ defined by $x'=\mh(\wt,x)$ on $\W_h\times\R$.
In particular, if $\w_0:=h$, then the maps $l_h(t):=\ml_h(\w_0{\cdot}t)$ and $u_h(t):=\muk_h(\w_0{\cdot}t)$
define the lower and upper bounded solutions of $x'=h(t,x)$.
In addition, the global pullback attractor of the induced process is $\{[l_h(s),u_h(s)]\,|\;s\in\R\}$
(see e.g.~\cite[Definition 1.12, Theorem 2.12 and Corollary 1.18]{carvalho1}, and the proof
of Proposition \ref{prop:5Dcoer}(ii)).
Recall that $x_h(t,s,x)$ satisfies $x'=h(t,x)$ and $x_h(s,s,x)=x$.
\par
Theorem \ref{th:6DCtressoluciones}, key in the proof of Theorem \ref{th:6Dcasos},
establishes the existence of three solutions which govern
the dynamics of \eqref{eq:6Dgtran} if \hyperlink{gd1}{\bf gd1}-\hyperlink{gd5}{\bf gd5} hold:
the two previously
described solutions $l_g$ and $u_g$, which are locally pullback attractive, and a
locally pullback repulsive one, $m_g$. Its proof requires the next previous result:
\begin{prop}\label{prop:6DCtreshiperbolicas}
Assume that \eqref{eq:6Dgtran} has three hyperbolic solutions $\tilde l_g<\tilde m_g<\tilde u_g$.
\begin{itemize}
\item[\rm(i)] If $g$ and $g_+$ satisfy all the conditions involving them in
    \hyperlink{gd1}{\bf gd1}-\hyperlink{gd5}{\bf gd5},
    then $\lim_{t\to\infty}(x_g(t,s,x)-\tilde u_{g_+}(t))=0$ for $x>\tilde m_g(s)$ and
    $\lim_{t\to\infty}(x_g(t,s,x)-\tilde l_{g_+}(t))=0$ for $x<\tilde m_g(s)$.
\item[\rm(ii)] If $g$ and $g_-$ satisfy all the conditions involving them in
    \hyperlink{gd1}{\bf gd1}-\hyperlink{gd5}{\bf gd5},
    then $t\mapsto x_g(t,s,x)$ is bounded from above (resp. from below) as time decreases
    if and only if $x\leq \tilde u_g(s)$ (resp. $x\geq \tilde l_g(s)$);
    and $\lim_{t\to-\infty}(x_g(t,s,x)-\tilde m_{g_-}(t))=0$ for $x\in(\tilde l_g(s),\tilde u_g(s))$.
\end{itemize}
\end{prop}
\begin{proof}
We proceed as in the proof of \cite[Proposition 3.5]{dno3}: Lemma \ref{lema:6Dhiposd}
allows us to use the information of Theorem~\ref{th:5DCthreecopies}
instead of that of \cite[Theorem 3.3]{dno3}, and Proposition \ref{prop:2haciaatras} provides the
necessary information on uniform separation of the \upalfa-limit and \upomeg-limit sets of
the points $(g,\tilde l_g(0))$, $(g,\tilde m_g(0))$ and $(g,\tilde u_g(0))$ for
the corresponding skew-product.
\end{proof}
\begin{teor}\label{th:6DCtressoluciones}
Assume that $g$ satisfies \hyperlink{gd1}{\bf gd1}-\hyperlink{gd5}{\bf gd5},
let $\tilde l_{g_\pm}<\tilde m_{g_\pm}<\tilde u_{g_\pm}$ be the hyperbolic solutions given by \hyperlink{gd5}{\bf gd5},
and let $l_g$ and $u_g$ be the lower and upper bounded solutions of
\eqref{eq:6Dgtran}. Then,
\begin{itemize}
\item[\rm(i)] $u_g$ and $l_g$ are the unique solutions of \eqref{eq:6Dgtran}
satisfying $\lim_{t\to-\infty}(u_g(t)-\tilde u_{g_-}(t))=0$
and $\lim_{t\to-\infty}(l_g(t)-\tilde l_{g_-}(t))=0$, and they are locally pullback attractive.
\item[\rm(ii)] There exists a unique solution
$m_g$ of \eqref{eq:6Dgtran} defined at least on a positive half-line and
satisfying $\lim_{t\to\infty}(m_g(t)-\tilde m_{g_+}(t))=0$,
and it is locally pullback repulsive.
\end{itemize}
Moreover, for $s\in\R$ in the interval of definition of $m_g$,
$\lim_{t\to\infty} (x_g(t,s,x)-\tilde u_{g_+}(t))=0$ if and only if $x>m_g(s)$,
and $\lim_{t\to\infty} (x_g(t,s,x)-\tilde l_{g_+}(t))=0$ if and only if $x<m_g(s)$.
In addition, for any $s\in\R$, $\lim_{t\to-\infty}(x_g(t,s,x)-\tilde m_{g_-}(t))=0$
if and only if $x\in(l_g(s),u_g(s))$.
\end{teor}
\begin{proof}
The assertions reproduce those of \cite[Theorem~3.7]{dno3}, formulated under more
restrictive hypotheses. The proof basically repeats step by step that one, using the information
provided by Proposition \ref{prop:6DCtreshiperbolicas} to check the last
assertions. The differences rely on the
first steps, which we detail. We take $m>0$ such that $\n{b}_\infty\le m$ for any
bounded solution of $x'=g(t,x)$ and $x'=g_\pm(t,x)$ (see Proposition \ref{prop:5Dcoer}(ii)).
Given $\ep>0$, Theorem \ref{th:2persistencia} provides $\delta_\pm>0$ such that,
if $\n{g_\pm-h_\pm}_{1,m}<\delta_\pm$, then each one of the equations $x'=h_\pm(t,x)$
has three hyperbolic solutions, at a uniform distance from those of
$x'=g_\pm(t,x)$ bounded by $\ep$. We choose $t^0=t^0(\ep)>0$ such that
$|g(t,x)-g_\pm(t,x)|<\delta_\pm/2$ and $|g_x(t,x)-(g_\pm)_x(t,x)|<\delta_\pm/2$ if
$\pm t\ge t^0$ and $|x|\le m$ (see Lemma \ref{lem:4Cderivadas}), and
define $f_\pm(t,x)$ as $g(t,x)$ if $\pm t>t^0$ and as $g_\pm(t,x)-g_\pm(\pm t^0,x)+g(\pm t^0,x)$
otherwise. The solutions of $x'=g(t,x)$ with values $\tilde l_{f_-}(-t^0)$, $\tilde u_{f_-}(-t^0)$ and $\tilde m_{f_+}(t^0)$ provide the solutions $l_g$, $u_g$ and $m_g$ of the statement,
as we can proof by repeating the remaining arguments of \cite{dno3}.
\end{proof}
We will denote $l_g$, $m_g$ and $u_g$ by $\tilde l_g$, $\tilde m_g$ and $\tilde u_g$ when they are hyperbolic.
Now we will formulate the announced result concerning the dynamical possibilities for \eqref{eq:6Dgtran}.
Recall that two uniformly separated
solutions are, by definition, bounded. Clearly, there exist (at least) two uniformly separated solutions
if and only if $l_g$ and $u_g$ satisfy this property.
\begin{teor}\label{th:6Dcasos}
Assume that $g$ satisfies \hyperlink{gd1}{\bf gd1}-\hyperlink{gd5}{\bf gd5},
let $\tilde l_{g_\pm}<\tilde m_{g_\pm}<\tilde u_{g_\pm}$ be the hyperbolic solutions given by \hyperlink{gd5}{\bf gd5},
and let $l_g$, $u_g$ and $m_g$ be the solutions of
\eqref{eq:6Dgtran} provided by Theorem {\rm\ref{th:6DCtressoluciones}}.
Then, the dynamics of the transition equation \eqref{eq:6Dgtran} fits in
one of the following dynamical scenarios:
\begin{itemize}[leftmargin=12pt]
\item \hypertarget{DCA}{{\sc Case A}}: there exist exactly three hyperbolic solutions, $\tilde l_g:=l_g$ and $\tilde u_g:=u_g$,
which are attractive, and $\tilde m_g:=m_g$, which is repulsive. In addition, the unique solution uniformly
separated from $\tilde l_g$ and $\tilde u_g$ is $\tilde m_g$.
In this case, $\tilde l_g<\tilde m_g<\tilde u_g$, $\lim_{t\to\pm\infty}(\tilde l_g(t)-\tilde l_{g_\pm}(t))=0$,
$\lim_{t\to\pm\infty}(\tilde m_g(t)-\tilde m_{g_\pm}(t))=0$ and $\lim_{t\to\pm\infty}(\tilde u_g(t)-\tilde u_{g_\pm}(t))=0$.
\item \hypertarget{DCB}{{\sc Case B}}: there exists exactly one hyperbolic solution, which is attractive, and uniformly separated only from
another solution, which is locally pullback attractive and repulsive. There are two possibilities:
\begin{itemize}[leftmargin=12pt]
\item {\sc Case B1}: $\tilde u_g=u_g$ is hyperbolic attractive and uniformly separated of $l_g=m_g$.
In this case, $\lim_{t\to\infty}(\tilde u_g(t)-\tilde u_{g_+}(t))=0$
and $\lim_{t\to\infty}(l_g(t)-\tilde m_{g_+}(t))=0$.
\item {\sc Case B2}: $\tilde l_g=l_g$ is hyperbolic attractive and uniformly separated of $m_g=u_g$.
In this case, $\lim_{t\to\infty}(\tilde l_g(t)-\tilde l_{g_+}(t))=0$
and $\lim_{t\to\infty}(u_g(t)-\tilde m_{g_+}(t))=0$.
\end{itemize}
\item \hypertarget{DCC}{{\sc Case C}}: there are no uniformly separated solutions. In this case,
$\tilde l_g=l_g$ and $\tilde u_g=u_g$ are the unique hyperbolic solutions, they are attractive,
and the locally pullback repulsive solution $m_g$ is unbounded.
There are two possibilities:
\begin{itemize}[leftmargin=12pt]
\item {\sc Case C1}: $m_g<\tilde l_g$ in its domain of definition.
In this case, $\lim_{t\to\infty}(\tilde l_g(t)-\tilde u_{g_+}(t))=\lim_{t\to\infty}(\tilde u_g(t)-\tilde u_{g_+}(t))=0$.
\item {\sc Case C2}: $m_g>\tilde u_g$ in its domain of definition.
In this case, $\lim_{t\to\infty}(\tilde l_g(t)-\tilde l_{g_+}(t))=0=\lim_{t\to\infty}(\tilde u_g(t)-\tilde l_{g_+}(t))=0$.
\end{itemize}
\end{itemize}
\end{teor}
\begin{proof}
Theorem \ref{th:6DCtressoluciones} allows us to repeat the proofs of
\cite[Theorems 3.9 and 3.10 and Corollary 3.11]{dno3}
under conditions \hyperlink{gd1}{\bf gd1}-\hyperlink{gd5}{\bf gd5}.
The statement of this theorem follows from those ones and Theorem \ref{th:6DCtressoluciones}.
\end{proof}
Figures 2, 3 and 4 of \cite{dno3} depict these five dynamical possibilities
in the case of a map $f$  which is asymptotically periodic with respect to $t$.
In addition, they can be characterized in terms of the
forward attraction properties of the global pullback attractor for \eqref{eq:6Dgtran}:
the proof of \cite[Proposition 3.12]{dno3} can be repeated in our framework.
\par
\begin{figure}
\centering
\includegraphics[width=\textwidth]{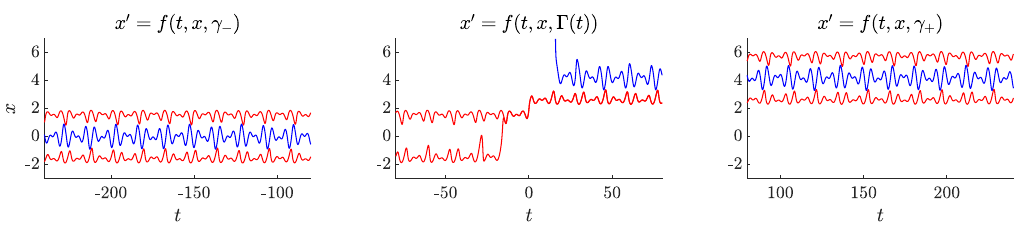}
\caption{
We define $\G(t)=\arctan(5t)/\pi+1/2\in(0,1)$,
$g(t,x)=-x^3+\sin(t)+\sin(\sqrt{2}\, t)+(5/2)x+\G(t)\,a(3x^2-3ax+a^2-5/2)$,
and $g_-(t,x)$ and $g_+(t,x)$ by replacing $\G(t)$ by $0$ and $1$ in $g(t,x)$, respectively.
Then, \protect\hyperlink{gd1}{\bf gd1}-\protect\hyperlink{gd5}{\bf gd5} hold.
The central panel shows that the dynamics of $x'=g(t,x)$ for $a=4.2$
corresponds to \protect\hyperlink{DCC}{\sc Case C2}: we
depict in red the two attractive hyperbolic solutions $\tilde l_g$ and $\tilde u_g$,
and in blue the unbounded locally pullback repulsive solution $m_g$.
It is easy to check that $g_-(t,x)\le g(t,x)$ for this choice of $a$.
But, as checked below, any bounded solution of $x'=g_-(t,x)$ (which are bound by the red
curves in the left panel) is below $m_{g_+}$ (depicted in blue in
the right panel), and hence neither the hypotheses nor the thesis
of Proposition~\ref{prop:6DCaseA}(i) are fulfilled.
It is also easy to check that $g_-(t,x-a)=g_+(t,x)$ and that $\pm g_-(t,r)<0$
for $\pm r>2$. Hence, $-2\le\tilde l_{g_-}(t)<\tilde u_{g_-}(t)\le 2$, and
$2.2\le\tilde l_{g_+}(t)<\tilde u_{g_+}(t)\le 6.2$, which implies the assertion.
In addition, since $\liminf_{t\to-\infty}(\tilde l_{g_+}(t)-u_g(t))>
\lim_{t\to-\infty}(\tilde u_{g_-}(t)-u_g(t))+0.2=0.2$, we get $u_g(t)<\tilde l_{g_+}(t)$
for $t\le t_0$; and $u_g'(t)<\tilde l'_{g_+}(t)$ if $u_g(t)=\tilde l_{g_+}(t)$,
from where we deduce that $u_g(t)<\tilde l_{g_+}(t)$ for all $t\in\R$. This is
only possible in \protect\hyperlink{DCC}{\sc Case C2}.}
\label{fig:badlyordered}
\end{figure}

In the line of Proposition \ref{prop:4CCaseA}, the next result establishes two conditions
precluding some of the five cases, and which, together, guarantee
\hyperlink{DCA}{\sc Case A}. The example depicted in Figure~\ref{fig:badlyordered}
shows that the hypotheses concerning the relative order of $\tilde m_{g_+}$
and the bounded solution $b_i$ are not superfluous. So, the conditions
are more exigent than those of the analogous result in the concave case,
Proposition~\ref{prop:4CCaseA}.
\begin{prop}\label{prop:6DCaseA}
Assume that $g$ satisfies \hyperlink{gd1}{\bf gd1}-\hyperlink{gd5}{\bf gd5}. Then,
\begin{itemize}
\item[\rm(i)] if there exists $h_1\colon\RR\to\R$ such that
    $h_1(t,x)\leq g(t,x)$ for all $(t,x)\in\RR$, and $x'=h_1(t,x)$
    has a bounded solution $b_1$ such that $\liminf_{t\to\infty}(b_1(t)-\tilde m_{g_+}(t))>0$,
    then $x'=g(t,x)$ is in \hyperlink{DCA}{\sc Case A},
    \hyperlink{DCB}{\sc B1} or \hyperlink{DCC}{\sc C1}.
\item[\rm(ii)] If there exists $h_2\colon\RR\to\R$ such that $h_2(t,x)\ge g(t,x)$
    for all $(t,x)\in\RR$, and $x'=h_2(t,x)$ has a bounded solution $b_2$ such that
    $\liminf_{t\to\infty}(\tilde m_{g_+}(t)-b_2(t))>0$, then $x'=g(t,x)$
    is in \hyperlink{DCA}{\sc Case A}, \hyperlink{DCB}{\sc B2} or \hyperlink{DCC}{\sc C2}.
\end{itemize}
\end{prop}
\begin{proof} Let us check (i): the second proof is analogous.
Let $u_g$ be the upper bounded solution for $x'=g(t,x)$ (see
Proposition \ref{prop:5Dcoer}(ii)). Then, $b_1\le u_g$
(see Proposition \ref{prop:5Dcoer}(iv)).
Hence, $\liminf_{t\to\infty}(u_g(t)-\tilde m_{g_+}(t))\ge
\liminf_{t\to\infty}(b_1(t)-\tilde m_{g_+}(t))>0$,
which according to Theorem \ref{th:6Dcasos} precludes
\hyperlink{DCB}{\sc Cases B2} and \hyperlink{DCC}{\sc C2}.
\end{proof}

\par
As in the concave case, we will focus on critical transitions associated
to one-parametric families of equations which occur when the dynamics
moves from \hyperlink{DCA}{\sc Case A} to one of the \hyperlink{CCC}{\sc Cases C}
of Theorem \ref{th:6Dcasos} as the parameter crosses a {\em critical value}.
Theorem \ref{th:6Dsaddle-node} shows the persistence of \hyperlink{DCA}{\sc Cases A},
\hyperlink{DCC}{C1} and \hyperlink{DCC}{C2} under small suitable parametric variations,
as well as the occurrence of a saddle-node bifurcation phenomenon when
\hyperlink{DCA}{\sc Case A} transits to one of the \hyperlink{DCB}{\sc Cases B}
as the parameter varies.
\begin{teor} \label{th:6Dsaddle-node}
Let $\mC\subseteq\R$ be an open interval, and let $\bar g\colon\R\times\R\times\mC\to\R$
be a map such that $g^c(t,x):=\bar g(t,x,c)$
satisfies \hyperlink{gd1}{\bf gd1}-\hyperlink{gd5}{\bf gd5} for all $c\in\mC$.
Let $\bar g_x$ be the partial derivative with respect to the second variable, and
assume that $\bar g$ and $\bar g_x$ are admissible on $\R\times\R\times\mC$.
Assume also that $\limsup_{x\to\pm\infty}(\pm\bar g(t,x,c))<0$ uniformly on $\R\times\mJ$
for any compact interval $\mJ\subset\mC$.
\begin{itemize}
\item[\rm(i)] Assume that there exist $c_1,c_2$ in $\mC$ with
$c_1<c_2$ such that the dynamics of $x'=g^c(t,x)$ is in \hyperlink{DCA}{\sc Case A} for $c=c_1$
and not for $c=c_2$.
If $c_0:=\inf\{c>c_1\,|\;\text{\hyperlink{DCA}{\sc Case A} does not hold\/}\}$,
then $c_0>c_1$. Let $\tilde l_{g^c}<\tilde m_{g^c}<\tilde u_{g^c}$ be the three hyperbolic solutions
of $x'=g^c(t,x)$ for $c\in[c_1,c_0)$. Then, the dynamics of $x'=g^{c_0}(t,x)$ is either
in \hyperlink{DCB}{\sc Case B1},
with $\lim_{c\to c_0^-}(\tilde m_{g^c}(t)-\tilde l_{g^c}(t))=0$
for all $t\in\R$, or in \hyperlink{DCB}{\sc Case B2},
with $\lim_{c\to c_0^-}(\tilde u_{g^c}(t)-\tilde m_{g^c}(t))=0$ for all $t\in\R$.
The results are analogous if $c_1>c_2$.
\item[\rm(ii)] Assume that there exist $c_3,c_4$ in $\mC$ with
$c_3<c_4$ such that the dynamics of $x'=g^c(t,x)$ is in \hyperlink{DCC}{\sc Case C1} for $c=c_3$
and not for $c=c_4$. If $c_0:=\inf\{c>c_3\,|\;\text{\hyperlink{DCC}{\sc Case C1} does not hold}\}$,
then $c_0>c_3$, and the dynamics of $x'=g^{c_0}(t,x)$ is in
\hyperlink{DCB}{\sc Case B1}. The results are analogous by replacing
\hyperlink{DCC}{\sc C1} and \hyperlink{DCB}{\sc B1} by \hyperlink{DCC}{\sc C2}
and \hyperlink{DCB}{\sc B2}, and also if $c_3>c_4$.
\end{itemize}
\end{teor}
\begin{proof}
As in the proof of Theorem \ref{th:4Csaddle-node}, the admissibility hypotheses
combined with Theorems \ref{th:6Dcasos} and \ref{th:2persistencia} guarantee
the persistence of \hyperlink{DCA}{\sc Case A} under small variations of $c$.
Let us check that also \hyperlink{DCC}{\sc Case C1} is persistent, assuming for
contradiction that $x'=g^{c_3}(t,x)$ is in this case and the existence of a sequence
$(c_n)$ with limit $c_3$ such that $x'=g^{c_n}(t,x)$ is not \hyperlink{DCC}{\sc Case C1}
for all $n\in\N$. (The same argument works for \hyperlink{DCC}{\sc Case C2}.)
Theorem \ref{th:2persistencia} shows that $x'=g^{c_n}(t,x)$ has two different
attractive hyperbolic solutions for large enough $n$, which must be
$\tilde l_{g^{c_n}}$ and $\tilde u_{g^{c_n}}$ (see Theorem \ref{th:6Dcasos}),
and which satisfy $\lim_{n\to\infty}\big\|\tilde l_{g^{c_n}}-\tilde l_{g^{c_3}}\,\big\|_\infty=
\lim_{n\to\infty}\big\|\tilde u_{g^{c_n}}-\tilde u_{g^{c_3}}\,\big\|_\infty=0$.
This precludes \hyperlink{DCB}{\sc Cases B}.
Let $\rho$ be the radio of the common domains of attraction also provided
by Theorem \ref{th:2persistencia}, and let us take $n_0$ such that
$\big\|\tilde l_{g^{c_n}}-\tilde l_{g^{c_3}}\,\big\|_\infty<\rho/3$ and
$\big\|\tilde u_{g^{c_n}}-\tilde u_{g^{c_3}}\,\big\|_\infty<\rho/3$ for all
$n\ge n_0$. If $n\ge n_0$, we deduce from $\lim_{t\to\infty}(\tilde u_{g^{c_3}}(t)-\tilde l_{g^{c_3}}(t))=0$ the existence of $t_0$ such that
$|\tilde u_{g^{c_n}}(t_0)-\tilde l_{g^{c_n}}(t_0)|<\rho$,
and hence that $\lim_{t\to\infty}(\tilde u_{g^{c_n}}(t)-\tilde l_{g^{c_n}}(t))=0$,
which precludes \hyperlink{DCA}{\sc Case A}. That is, $x'=g^{c_n}(t,x)$
is in \hyperlink{DCC}{\sc Case C2} for all $n\ge n_0$.
\par
Let $k$ be a common bound for the $\n{{\cdot}}_\infty$-norm of the bounded solutions of
$x'=g^{c_3}(t,x)$ and $x'=g_+^{c_3}(t,x)$, and
let $\ep>0$ be smaller than $\inf_{t\in\R}(\tilde u_{g_+^{c_3}}(t)-\tilde m_{g_+^{c_3}}(t))$
and $\inf_{t\in\R}(\tilde m_{g_+^{c_3}}(t)-\tilde l_{g_+^{c_3}}(t))$.
Theorem \ref{th:2persistencia} applied to $\ep/4$ provides $\delta>0$
such that, if $f$ is $C^1$-admissible and $\n{f-g_+^{c_3}}_{1,k}<\delta$,
then $x'=f(t,x)$ has three hyperbolic solutions at a $\|{\cdot}\|_\infty$-distance
of those of $x'=g^{c_3}_+(t,x)$ less that $\ep/4$, and hence with a separation
between of at least $\ep/2$. The admissibility of $\bar g$ and condition
\hyperlink{gd2}{\bf gd2} applied to $g^{c_3}$ and $g^{c_3}_+$
allow us to choose $t_0$ and $n_0$ large enough to get
\[
\sup_{(t,x)\in[t_0,\infty)\times[-k,k]} |g^{c_n}(t,x)-g_+^{c_3}(t,x)|+\!\!
\sup_{(t,x)\in[t_0,\infty)\times[-k,k]} |(g^{c_n})_x-(g_+^{c_3})_x(t,x)|<\delta
\]
for all $n\ge n_0$: we just write $|g^{c_n}-g_+^{c_3}|\le|g^{c_n}-g^{c_3}|+|g^{c_3}-g_+^{c_3}|$,
do the same with the derivatives, and apply Lemma \ref{lem:4Cderivadas}.
Let us define $f^{c_n}_+(t,x)$ by truncating $g^{c_n}$ at $t_0$,
as in the proof of Theorem \ref{th:6DCtressoluciones}. Since
$\n{f^{c_n}_+-g^{c_3}_+}_{1,k}<\delta$, $x'=g^{c_n}(t,x)$
has three (possibly locally defined) solutions,
$b_1^{c_n}<b_2^{c_n}<b_3^{c_n}$, with $|b_i^{c_n}(t)|\le k+\ep/4$ and
$b_{i+1}^{c_n}(t)-b^{c_n}_i(t)\ge\ep/2$ for all $t\ge t_0$ and $n\ge n_0$.
We define $\bar b^{c_3}_i(t):=\lim_{n\to\infty} b_i^{c_n}(t)$ for $i=1,2,3$,
and get three solutions of $x'=g^{c_3}(t,x)$ defined and uniformly separated
by $\ep/2$ on $[t_0,\infty)$.
Since we are in \hyperlink{DCC}{\sc Case C2},
we have $b_2^{c_n}(t)\ge \tilde u_{g^{c_n}}(t)$ for all $t\in[t_0,\infty)$:
there cannot be two different solutions separated on $[t_0,\infty)$
strictly below $\tilde u_{g^{c_n}}$. Hence,
$\bar b_3^{c_3}(t)\ge\bar b_2^{c_3}(t)+\ep/2\ge\tilde u_{g^{c_3}}(t)+\ep/2$
for all $t\in[t_0,\infty)$,
which is not possible in \hyperlink{DCC}{\sc Case C1}. This is the sought-for contradiction.
\par
Let us complete the proof of (i) with $c_1<c_2$.
The persistence of \hyperlink{DCA}{\sc Cases A} and \hyperlink{DCC}{\sc C}
ensures that $c_0>c_1$ and that $x'=g^{c_0}(t,x)$ is in one of the \hyperlink{DCB}{\sc Cases B},
say \hyperlink{DCB}{\sc B1}. Let us prove that
$\lim_{c\to c_0^-}(\tilde m_{g^c}(t)-\tilde l_{g^c}(t))=0$
by checking that, given $(c_n)\uparrow c_0$,
$\lim_{n\to\infty}\tilde m_{g^{c_n}}(t)=
\lim_{n\to\infty}\tilde l_{g^{c_n}}(t)=m_{g^{c_0}}(t)$ for all $t\in\R$.
The existence of these limits follows from the existence of a common bound
for all the bounded solutions if $n$ is large enough. A new application
of last assertion of Theorem \ref{th:2persistencia} applied to
$\tilde u_{g^{c_0}}$ and its approximants $\tilde u_{g_{c_n}}$
shows that $\tilde u_{g_{c_n}}(t)-\tilde l_{g_{c_n}}(t)>
\tilde u_{g_{c_n}}(t)-\tilde m_{g_{c_n}}(t)\ge\rho$
if $n$ is large enough, with a common $\rho>0$. And hence both
limits are $m_{g^{c_0}}$, which is the unique bounded solution of $x'=g^{c_0}(t,x)$
uniformly separated from $\tilde u_{g^{c_0}}=\lim_{n\to\infty}\tilde u_{g^{c_n}}$
(see Theorem \ref{th:6Dcasos}). The remaining situations are proved
with similar arguments.
\par
To complete the proof of (ii) if $c_3<c_4$ and with $x'=g^{c_3}(t,x)$ in \hyperlink{DCC}{\sc Case C2},
we deduce from the proved persistence that $c_0>c_3$ and that the dynamics of $x'=g^{c_0}(t,x)$ is in one of the \hyperlink{DCB}{\sc Cases B}. Let us assume for
contradiction that it is in \hyperlink{DCB}{\sc Case B1}, so that $\tilde u_{g^{c_0}}$ is
hyperbolic. We take $(c_n)\uparrow c_0$, with $x'=g^{c_n}(t,x)$ in \hyperlink{DCC}{\sc Case C2},
and get the sought-for contradiction by repeating the last paragraph of the proof of
the persistence of \hyperlink{DCC}{\sc Case C1}: just replace $c_3$ by $c_0$.
The remaining cases are proved with similar arguments.
\end{proof}
We complete this part with an analogue of Corollary \ref{coro:4Cunpunto} for the d-concave case.
\begin{teor} \label{th:6Ddospuntos}
Let $\mC\subseteq\R$ be an open interval, and let
$\{g^c\,|\;c\in\mC\}$ be a family of functions satisfying \hyperlink{gd1}{\bf gd1}-\hyperlink{gd5}{\bf gd5} and
such that, if $\bar g(t,x,c):=g^c(t,x)$, then $\bar g$ and $\bar g_x$ are admissible on $\R\times\R\times\mC$.
Assume that there exists $\bar c\in\mC$ such that the dynamics of $x'=g^{\bar c}(t,x)$
is in \hyperlink{DCB}{\sc Case B1} (resp. \hyperlink{DCB}{\sc Case B2}), and such that,
for all $c_-,c_+\in\mC$ with $c_-<\bar c<c_+$: $g^{c_-}(t,x)\le g^{\bar c}(t,x)\le g^{c_+}(t,x)$
for all $(t,x)\in\RR$; and there exist $t_{c_-}$ and $t_{c_+}$ such that the first
and second inequalities are strict for $t=t_{c_-}$ and $t=t_{c_+}$ (respectively) and all $x\in\R$.
Then, there exists $\rho>0$ such that $x'=g^c(t,x)$ is in \hyperlink{DCA}{\sc Case A}
(resp. \hyperlink{DCC}{\sc Case C2}) for $c\in(\bar c-\rho,\bar c)$
and in \hyperlink{DCC}{\sc Case C1} (resp.~\hyperlink{DCA}{\sc Case A}) for $c\in(\bar c,\bar c+\rho)$.
\end{teor}
\begin{proof}
Let $l_c$, $m_c$ and $u_c$ be the three solutions of $x'=g^c(t,x)$ given by Theorem \ref{th:6DCtressoluciones}.
Let $g^c_+$ be the globally bounded and $C^2$-admissible function
associated to $g^c$ by \hyperlink{gd2}{\bf gd2} at $+\infty$, and let
$\tilde l_{g_+^c}<\tilde m_{g_+^c}<\tilde u_{g_+^c}$ be the three hyperbolic solutions of
$x'=g_+^c(t,x)$ provided by \hyperlink{gd5}{\bf gd5}.
Let $k$ be a common bound for the $\n{\cdot}_\infty$-norm of these three solutions.
We take $2\ep>0$ smaller than $\inf_{t\in\R}(\tilde u_{g_+^{\bar c}}(t)-\tilde m_{g_+^{\bar c}}(t))$ and
$\inf_{t\in\R}(\tilde m_{g_+^{\bar c}}(t)-\tilde l_{g_+^{\bar c}}(t))$. Then,
Theorem~\ref{th:2persistencia} provides
$\delta>0$ such that if $\|g_+^{\bar c}-g\|_{1,k}<\delta$ for a $C^1$-admissible map $g$,
then $x'=g(t,x)$ has three hyperbolic solutions at $\n{{\cdot}}_\infty$-distance of
those of $x'=g_+^{\bar c}(t,x)$ less than $\ep$.
\par
The admissibility of $\bar g$ and $(\bar g)_x$ provide
$\rho>0$ such that $\|g^{\bar c}-g^c\|_{1,k}<\delta/3$ for all $c\in(\bar c-\rho,\bar c+\rho)$.
Our hypotheses provide $t_c\ge 0$ with
$\sup_{(t,x)\in[t_c,\infty)\times[-k,k]}|g^c(t,x)-g_+^c(t,x)|+
\sup_{(t,x)\in[t_c,\infty)\times[-k,k]}|g^c_x(t,x)-(g_+^c)_x(t,x)|<\delta/3$: see
\hyperlink{gd2}{\bf gd2} and Lemma \ref{lem:4Cderivadas}.
We take $c^*\in (\bar c-\rho,\bar c+\rho)$ and $t^0:=\max(t_{\bar c},t_{c^*})$. For
$h(t,x)$, we denote by $\hat h(t,x)$ the map given by
$h(t,x)$ for $t\ge t^0$ and by $g^{\bar c}_+(t,x)-g^{\bar c}_+(t^0,x)+h(t^0,x)$
for $t<t^0$. In this way, we construct $\hat g^{c^*}_+$, $\hat g^{\bar c}$ and
$\hat g^{c^*}$ from $g^{c^*}_+$, $g^{\bar c}$ and $g^{c^*}$, and note that they are
$C^1$-admissible.
Then: $\|\hat g^{\bar c}-\hat g^{c^*}\|_{1,k}<\delta/3$,
since the difference is $\hat g^{\bar c}(t_0,x)-\hat g^{c^*}(t_0,x)$
for $t<t^0$ and $g^{\bar c}(t,x)-g^{c^*}(t,x)$ for $t\ge t^0$; and
$\|g^{c^*}-g^{c^*}_+\|_{1,k}<\delta/3$ and
$\|g_+^{\bar c}-\hat g^{\bar c}\|_{1,k}<\delta/3$ for analogous reasons.
So, $\|g_+^{\bar c}-\hat g^{c^*}_+\|_{1,k}<\delta$,
and hence $x'=\hat g^{c^*}_+(t,x)$ has three hyperbolic solutions at a distance
less than $\ep$ of those of $x'=g^{\bar c}_+(t,x)$.
In addition, since they solve $x'=g^{c^*}_+(t,x)$ on $[t^0,\infty)$,
the middle one coincides with $\tilde m_{g^{c^*}_+}$ on $[t^0,\infty)$: $\tilde m_{g^{c^*}_+}$
is the unique solution of $x'=g^{c^*}_+(t,x)$ uniformly separated from
two other solutions as $t$ increases. Hence,
$\tilde u_{g^{\bar c}_+}(t)\ge\tilde m_{g^{c^*}_+}(t)+\ep$
and $\tilde l_{g^{\bar c}_+}(t)\le\tilde m_{g^{c^*}_+}(t)-\ep$ for $t\ge t^0$.
\par
Let us assume that $x'=g^c(t,x)$ is in \hyperlink{DCB}{\sc Case B1}
for $c=\bar c$, associate $\rho$ to $\bar c$ as above, and check that
$x'=g^c(t,x)$ is in \hyperlink{DCC}{\sc Case C1}
for any $c^*\in(\bar c,\bar c+\rho)$, which we fix.
Since $g^{c^*}(t,l_{c^*}(t))\ge g^{\bar c}(t,l_{c^*}(t))$ for all $t\in\R$,
Proposition \ref{prop:5Dcoer}(iv) shows that $l_{\bar c}\le l_{c^*}$.
These inequalities combined with $g^{c^*}(t_0,l_{c^*}(t_0))>g^{\bar c}(t_0,l_{c^*}(t_0))$
yield $l_{\bar c}(t)<l_{c^*}(t)$ for all $t>t_0$, and hence
$\lim_{t\to\infty}(x_{\bar c}(t,t_0+1,l_{c^*}(t_0+1))-\tilde u_{g_+^{\bar c}}(t))=0$:
see Theorems \ref{th:6DCtressoluciones} and \ref{th:6Dcasos}.
A standard comparison argument shows that $x_{\bar c}(t,t_0+1,l_{c^*}(t_0+1))\le l_{c^*}(t)$
for $t\ge t_0+1$, and hence $\liminf_{t\rightarrow\infty}(l_{c^*}(t)-\tilde u_{g_+^{\bar c}}(t))\geq 0$.
Thus, $\liminf_{t\rightarrow\infty}(l_{c^*}(t)-\tilde m_{g_+^{c^*}}(t))\ge\ep$,
which means \hyperlink{DCC}{\sc Case C1} for $c^*$: see Theorem \ref{th:6Dcasos}.
Similar comparison arguments show that \hyperlink{DCA}{\sc Case A} holds
$c^*\in(\bar c-\rho,\bar c)$, as well as the stated properties if
$x'=g^{\bar c}(t,x)$ is in \hyperlink{DCB}{\sc Case B2}.
\end{proof}
\subsection{Some scenarios of critical transitions in the d-concave case}\label{subsec:61}
Let $\mI\subseteq\R$ be an open interval, and let $f\colon\RR\times\mI\to\R$ and
$\G,\G_-,\G_+\colon\RR\to\mI$ satisfy
\begin{list}{}{\leftmargin 20pt}
\item[\hypertarget{fd1}{\bf fd1}] there exist the derivatives $f_x$, $f_{xx}$, $f_\gamma$
and $f_{\gamma\gamma}$, and $f$, $f_x$, $f_{xx}$, $f_\gamma$ and
$f_{\gamma\gamma}$ are admissible on $\R\times\R\times\mI$.
\item[\hypertarget{fd2}{\bf fd2}] $\G,\G_-$ and $\G_+$ take values in $[a,b]\subset\mI$, are $C^2$-admissible,
and $\lim_{t\to\pm\infty}(\G(t,x)-\G_\pm(t,x))=0$ uniformly on each compact subset $\mJ\subset\R$.
\item[\hypertarget{fd3}{\bf fd3}] $\limsup_{x\to\pm\infty}(\pm f(t,x,\gamma))<0$
uniformly in $(t,\gamma)\in\R\times\mJ$ for all compact interval $\mJ\subset\mI$.
\item[\hypertarget{fd4}{\bf fd4}]
$\inf_{t\in\R}\big((\partial^2/\partial x^2)f(t,x,\G_\pm(t,x))|_{x=x_1}-
(\partial^2/\partial x^2)f(t,x,\G_\pm(t,x))|_{x=x_2}\big)>0$ whenever $x_1<x_2$.
\item[\hypertarget{fd5}{\bf fd5}] Each equation $x'=f(t,x,\G_\pm(t,x))$ has
three hyperbolic solutions $\tilde l_{\G_{\pm}}<\tilde m_{\G_{\pm}}<\tilde u_{\G_{\pm}}$.
\end{list}
\par
With the same abuse of language as in Subsection \ref{subsec:41},
we will say that $(f,\G)$ satisfies \hyperlink{fd1}{\bf fd1}-\hyperlink{fd5}{\bf fd5}
if there exist maps $\G_-$ and $\G_+$ such that the previous conditions are satisfied,
and refer to the equations
\begin{equation}\label{eq:6Dlim}
 x'=f(t,x,\G_-(t,x)) \quad\text{and}\quad x'=f(t,x,\G_+(t,x))
\end{equation}
as the ``past" and ``future" of
\begin{equation}\label{eq:6Dtran}
 x'=f(t,x,\G(t,x))\,.
\end{equation}
We can easily prove the next result, analogous to Proposition \ref{prop:4Chiposcomp}:
\begin{prop}\label{prop:6Dhiposcomp}
Assume that $(f,\G)$ satisfies \hyperlink{fd1}{\bf fd1}-\hyperlink{fd5}{\bf fd5}.
Then, the maps $g,g_-,g_+$ given by $g(t,x):=f(t,x,\G(t,x))$, $g_-(t,x):=f(t,x,\G_-(t,x))$
and $g_+(t,x):=f(t,x,\G_+(t,x))$ satisfy the conditions
\hyperlink{gd1}{\bf gd1}-\hyperlink{gd5}{\bf gd5}. Therefore, the dynamical
possibilities for \eqref{eq:6Dtran} are those described in Theorem {\rm\ref{th:6Dcasos}}.
\end{prop}
\begin{nota}\label{rm:6Dpar}
The conditions on $(f,\G)$ can be weakened, as in Remark \ref{rm:4Cpar}.1.
\end{nota}
In the line of Theorem \ref{th:4Cnotran}, Theorem \ref{th:6DCnotipping},
based on Proposition \ref{prop:6DCaseA}, establishes conditions
providing a {\em safety interval\/} $[\gamma_1,\gamma_2]$: if
$\G(\RR)\subseteq[\gamma_1,\gamma_2]$, then neither rate-induced tipping nor
phase-induced tipping takes place. As seen in its statement, this safety interval depends
on the $C^2$-admissible function $\G_+$ determining the future equation,
which is an important difference with respect to the concave analogue,
Theorem \ref{th:4Cnotran}.
And Theorems \ref{th:6Dsize} and \ref{th:6Dsize2},
based on Theorem \ref{th:6Ddospuntos},
provide the d-concave analogues of Theorems \ref{th:4Csize} and \ref{th:4Csize2}:
under hypotheses precluding the transition map $\G^d$
to take values in any fixed interval for all the values of the
parameter, they show either the absence of critical
transition or the occurrence of exactly two tipping points.
Looking for clarity in their statements, we just analyze the situation precluding
$\G^d$ to be always bounded from below.
\begin{teor}\label{th:6DCnotipping}
Assume that $(f,\G)$ satisfies \hyperlink{fd1}{\bf fd1}-\hyperlink{fd5}{\bf fd5}.
Assume also that
$\gamma\mapsto f(t,x,\gamma)$ is nondecreasing for all $(t,x)\in\RR$,
and that there exist $\gamma_1\le \gamma_2$ such that: $\G(\RR)\subseteq[\gamma_1,\gamma_2]$,
$x'=f(t,x,\gamma_1)$ has a bounded solution $b_1$ with
$\liminf_{t\to\infty}(b_1(t)-\tilde m_{\Gamma_+}(t))>0$;
and $x'=f(t,x,\gamma_2)$ has a bounded solution $b_2$ with
$\liminf_{t\to\infty}(\tilde m_{\Gamma_+}(t)-b_2(t))>0$.
Then, \eqref{eq:6Dtran} is in \hyperlink{DCA}{\sc Case A}.
\par
If, in addition, we assume that $\G_\pm$
do not depend on $t$, then the equations $x'=f(t,x,\G(ct,x))$ and
$x'=f(t,x,\G(t+c,x))$ are in \hyperlink{DCA}{\sc Case A} for all $c>0$ and $c\in\R$,
respectively: there is neither rate-induced tipping nor phase-induced tipping.
\end{teor}
\begin{proof} Take $h_i(t,x)=f(t,x,\gamma_i)$ for $i=1,2$ and apply Proposition~\ref{prop:6DCaseA}.
\end{proof}
\begin{teor}\label{th:6Dsize}
Assume that $\mI=\R$. Let $\G\colon\RR\to\R$ and $\G_0\colon\RR\to[0,\infty)$
be globally bounded and $C^2$-admissible, and such
that the pair $(f,\G+d\,\G_0)$ satisfies \hyperlink{fd1}{\bf fd1}-\hyperlink{fd5}{\bf fd5}
for all $d\in\R$. Assume that $\G_0(t_0,x)>0$ for all $x\in\R$ and a $t_0\in\R$. Assume also that
$\gamma\mapsto f(t,x,\gamma)$ is strictly increasing on $\R$ for all $(t,x)\in\RR$,
with $\lim_{\gamma\to\pm\infty} f(t,x,\gamma)=\pm\infty$ uniformly on compact sets of $\RR$.
Then,
\begin{equation}\label{eq:6uniparametric}
 x'=f(t,x,\G(t,x)+d\,\G_0(t,x))
\end{equation}
is either in \hyperlink{DCC}{\sc Case C1} for all $d\in\R$,
in \hyperlink{DCC}{\sc Case C2} for all $d\in\R$, or there exist
$d_-<d_+$ such that it is in \hyperlink{DCC}{\sc Case C2} for $d<d_-$,
in \hyperlink{DCB}{\sc Case B2} for $d=d_-$, in \hyperlink{DCA}{\sc Case A} for
$d\in(d_-,d_+)$, in \hyperlink{DCB}{\sc Case B1} for $d=d_+$,
and in \hyperlink{DCC}{\sc Case C1} for $d>d_+$.
\end{teor}
\begin{proof}
It is easy to check that the family of maps $g^d(t,x):=f(t,x,\G(t,x)+d\,\G_0(t,x))$ satisfies all the
hypotheses of Theorem \ref{th:6Ddospuntos}, with $\mC=\R$ and $\bar c$ equal to any $\bar d\in\R$.
We assume that \eqref{eq:6uniparametric}$_d$ is in
\hyperlink{DCA}{\sc Case A} for $d=\bar d$ and, for contradiction that
$d_+:=\inf\{d>\bar d\,|\;\text{\hyperlink{DCA}{\sc Case A} does not hold}\}=\infty$.
Theorem~\ref{th:2persistencia} ensures that $\bar d<d_+$.
Let $\tilde l_d<\tilde m_d<\tilde u_d$ be the three hyperbolic solutions of \eqref{eq:6uniparametric}$_d$
for $d\in[\bar d,d_+)$. Proposition
\ref{prop:5Dcoer}(iv) yields $\tilde l_{\bar d}\le \tilde l_d$
for all $d\in[\bar d,d_+)$. Let us prove that $\tilde m_d\le\tilde m_{\bar d}$ for all
$d\in[\bar d,d_+)$. Clearly, it suffices to check it for $d\in[\bar d,\bar d+\rho)$
for a small $\rho>0$, which we choose (applying Theorem \ref{th:2persistencia})
to ensure $\tilde u_{\bar d}>\tilde m_d+\ep$ for an
$\ep>0$ common for all $d\in[\bar d,\bar d+\rho)$. For contradiction, we take
$x_0\in(\tilde m_{\bar d}(0),\tilde m_d(0))$. Theorem \ref{th:6Dcasos}
yields $\lim_{t\to\infty}(x_d(t,0,x_0)-\tilde l_d(t))=0$ and
$\lim_{t\to\infty}(x_{\bar d}(t,0,x_0)-\tilde u_{\bar d}(t))=0$, and hence
$x_d(t,0,x_0)\ge x_{\bar d}(t,0,x_0))$ for $t\ge 0$ yields
$\limsup_{t\to\infty}(x_d(t,0,x_0)-\tilde u_{\bar d}(t))\ge 0$. Therefore,
$\limsup_{t\to\infty}(x_d(t,0,x_0)-\tilde m_d(t))\ge \ep$, impossible.
\par
So, we can take $m_1<m_2$ such that
$m_1\leq \tilde l_{\bar d}\le\tilde l_d<\tilde m_d\le\tilde m_{\bar d}\le m_2$.
An argument similar to that involving $k_d$ in the proof of Theorem~\ref{th:4Csize}
provides the sought-for contradiction.
Similarly, $d_-:=\sup\{d<\bar d\colon\,|\;\text{\hyperlink{DCA}{\sc Case A} does not hold}\}<\bar d$
is finite. It easy to deduce from Theorem \ref{th:6Ddospuntos}
that the variation is the stated one: \hyperlink{DCC}{\sc C2} for all $d<d_-$,
\hyperlink{DCB}{\sc B2} at $d_-$, \hyperlink{DCB}{\sc B1} at $d_+$, and \hyperlink{DCC}{\sc C1} for all $d>d_+$.
\par
Since \hyperlink{DCA}{\sc Case A} cannot occur for all $d$, there are equations either in
\hyperlink{DCC}{\sc Case C1} or in \hyperlink{DCC}{\sc Case C2}. Assume
that \eqref{eq:6uniparametric}$_d$ is in \hyperlink{DCC}{\sc Case C2} for $d=\bar d$,
but not for all $d$. Theorem \ref{th:6Dsaddle-node} and \ref{th:6Ddospuntos} ensure the existence of
$d_->\bar d$ for which the dynamics is in \hyperlink{DCB}{\sc Case B2} and that
\hyperlink{DCA}{\sc Case A} holds for close values of $d>d_-$. So, we are in the situation
of the previous paragraphs. The argument
is analogous if \eqref{eq:6uniparametric}$_d$ is in \hyperlink{DCC}{\sc Case C1} for
$d=\bar d$, but not for all $d$, and the proof is complete.
\end{proof}
By reviewing the previous proof, we observe that we have proved the next result:
\begin{teor}\label{th:6Dsize2}
Assume that $\mI=\R$. Let $\G\colon\RR\to\R$ and $\G_0\colon\RR\to[0,\infty)$
be globally bounded and $C^2$-admissible, and such
that the pair $(f,\G+d\,\G_0)$ satisfies \hyperlink{fd1}{\bf fd1}-\hyperlink{fd5}{\bf fd5}
for all $d\in\R$. Assume that there exists $\bar d\in\R$ such that
\begin{equation}\label{eq:6DGG_02}
 x'=f(t,x,\G(t,x)+d\,\G_0(t,x))
\end{equation}
is in \hyperlink{DCA}{\sc Case A} for $d=\bar d$, and let
$\tilde l_{\bar d}<\tilde m_{\bar d}<\tilde u_{\bar d}$ be its three hyperbolic solutions.
Let $m_1<m_2$ and $m_3<m_4$ be such that
$m_1\leq\tilde l_{\bar d}(t)<\tilde m_{\bar d}(t)\leq m_2$ for all $t\in\R$ and
$m_3\leq\tilde m_{\bar d}(t)<\tilde u_{\bar d}(t)\leq m_4$ for all $t\in\R$.
\begin{enumerate}[label={\rm(\arabic*)},leftmargin=*]
\item Assume that there exists $t_0$
such that $\G_0(t_0,x)>0$ for all $x\in[m_1,m_2]$,
that $\gamma\mapsto f(t,x,\gamma)$ is nondecreasing for all $(t,x)\in\RR$ and strictly increasing
for $(t,x)\in\R\times [m_1,m_2]$, with
$\lim_{\gamma\to\infty} f(t,x,\gamma)=\infty$ uniformly on compact sets of $\R\times [m_1,m_2]$.
Then, there exists $d_+>\bar d$ such that \eqref{eq:6DGG_02}$_d$
is in \hyperlink{DCC}{\sc Case C1} for $d>d_+$,
in \hyperlink{DCB}{\sc Case B1} for $d=d_+$, in \hyperlink{DCA}{\sc Case A} for
$d\in[\bar d,d_+)$.
\item Assume that there exists $t_0$
such that $\G_0(t_0,x)>0$ for all $x\in[m_3,m_4]$,
that $\gamma\mapsto f(t,x,\gamma)$ is nondecreasing for all $(t,x)\in\RR$ and strictly increasing
for $(t,x)\in\R\times [m_3,m_4]$, with
$\lim_{\gamma\to-\infty} f(t,x,\gamma)=-\infty$ uniformly on compact sets of $\R\times [m_3,m_4]$.
Then, there exists $d_-<\bar d$ such that \eqref{eq:6DGG_02}$_d$
is in \hyperlink{DCC}{\sc Case C2} for $d<d_-$,
in \hyperlink{DCB}{\sc Case B2} for $d=d_-$, in \hyperlink{DCA}{\sc Case A} for
$d\in(d_-,\bar d]$.
\end{enumerate}
\end{teor}
An analysis similar to that of Remark \ref{rm:4Ccasoparticular} applies to this
d-concave case.
\subsection{Numerical simulations in asymptotically d-concave models}\label{subsec:62}
In this section, we consider two different single-species population models whose
internal dynamics are driven by nonautonomous cubic equations and which include
predation phenomena. The intrinsic cubic dynamics is indebted to the Allee effect
(see \cite{courchamp1,dolo1,dno3}), e.g., due to some breeding cooperative mechanism or to an
easier mate finding. In both cases, the evolution of the population is modeled by
\begin{equation}\label{eq:6populationsimulationgeneral}
 x'=r(t)\,x\,\left(1-\frac{x}{K(t)}\right)\frac{x-S(t)}{K(t)}+\Delta(t,x)\,,
\end{equation}
where we assume $r$, $K$ and $S$ to be quasiperiodic functions with $r$ and $K$
positively bounded from below, and $\Delta$ to be $C^2$-admissible.
So, analogously to Subsection~\ref{subsec:42},
the map $h(t,x,\delta):=r(t)\,x\,(1-x/K(t))(x-S(t))/K(t)+\delta$
satisfies \hyperlink{fd1}{\bf fd1} and \hyperlink{fd3}{\bf fd3}.
In addition, we will assume that $(h,\Delta)$ satisfies \hyperlink{fd2}{\bf fd2},
\hyperlink{fd4}{\bf fd4} and \hyperlink{fd5}{\bf fd5} for some
maps $\Delta_\pm$. The meaning of $r$ and $K$ is the same as in Subsection~\ref{subsec:42},
while $S$ (on which we assume $K(t)+S(t)\geq 0$ for all $t\in\R$) stands for the force of
the Allee effect, and $\Delta$ models the contribution of a single external effect: predation.
\par
Each one of the two examples tries to emphasize some of the novel aspects of the theory
presented in this paper: the possibility of non asymptotically constant transition functions
in the first one, and the possibility of intrinsically $x$-dependent transition functions in the second one.
As in Subsection~\ref{subsec:42}, we find
\hyperlink{CCA}{\sc Cases A}, \hyperlink{CCB}{B} and \hyperlink{CCC}{C}
for different values of certain parameters, and we point to
certain parametric variations as possible causes of tipping.
Throughout the section,
\hyperlink{DCA}{\sc Case A} means the survival of the species, while \hyperlink{DCC}{\sc Case C2}
means its extinction (and \hyperlink{DCB}{\sc Case B2} is the highly unstable situation which separates the other two).
\begin{exa} We begin by assuming that, in \eqref{eq:6populationsimulationgeneral}, the predation
is modeled by a H\"{o}lling type III functional response term $-\gamma\, x^2/(b(t)+x^2)$,
where $\gamma$ and $b>0$ have the same meaning as in Example~\ref{ex:4Cejemplo}:
\begin{equation}\label{eq:6Ddeff}
 x'=r(t)\,x\,\left(1-\frac{x}{K(t)}\right)\frac{x-S(t)}{K(t)}-\gamma\,\frac{x^2}{b(t)+x^2}\,.
\end{equation}
Next, we assume that the population is attacked by a predator species which behaves as follows:
the habitat is initially free of predators; at a certain time a group of predators arrives
at the ecosystem, which they leave after some time; and this behavior repeats yearly.
Such a pattern may correspond to the colonization of a new patch by a migratory species of predators,
due to the reproductive, nutritional, breeding or wintering interest of the habitat.
(See, e.g., \cite{aabpf} for a study on the evolution of some migration patterns
of common swift, an insectivorous bird.)

\begin{figure}
     \centering
         \includegraphics[width=\textwidth]{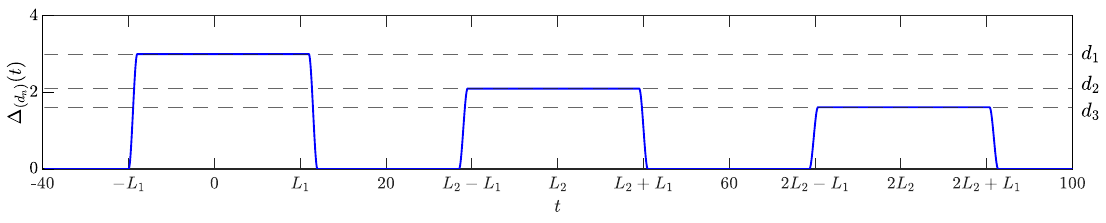}
         \caption{The transition map $\Delta_{(d_n)}$ defined in \eqref{eq:series} for
         $\rho=1$, $L_1=10$, $L_2=40$, $d_+=0.5$, $d=2.5$, $d_n=d_++ d/((n-1)/4+1)^2$,
         and $p_n=(n-1)L_2+(-1)^{(n-1)}/n$ for all $n\in\N$.}
        \label{fig:seriesimpulses}
\end{figure}
\begin{figure}
     \centering
         \includegraphics[width=\textwidth]{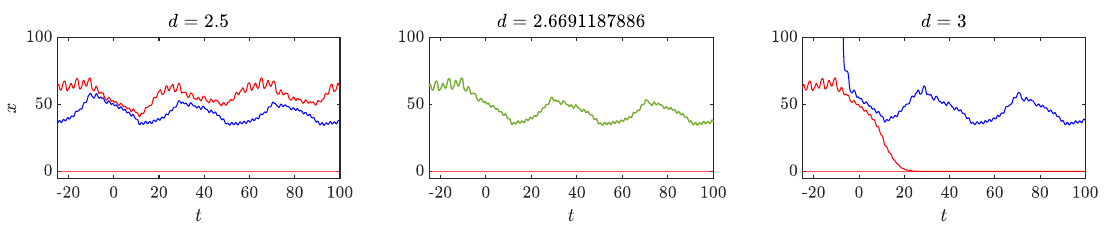}
         \caption{Numerical depiction of the existence of a unique size-tipping
         point for \eqref{eq:62transitiontoperiodic}$_d$.
         The central panel shows the dynamics for an accurate approximation to the tipping point $d_0$:
         the two upper hyperbolic solutions are so close within the representation window that are a good
         approximation (green) to the nonhyperbolic solution of \protect\hyperlink{DCB}{\sc Case B2}.
         The left panel depicts \protect\hyperlink{DCA}{\sc Case A}, which is the dynamics for
         any $d\in[0,d_0)$ and means survival: the attractive hyperbolic solutions in red, and the repulsive one in blue.
         The right panel depicts \protect\hyperlink{DCC}{\sc Case C2}, which is the dynamics for any $d>d_0$ and means
         extinction: the hyperbolic solutions in red, and the locally pullback repulsive solution in blue.}
        \label{fig:variablend}
\end{figure}
Let $L_2$ be the length of the year.
We assume that the $n$-th predation season occurs during the time $[p_n-L_1-\rho,p_n+L_1+\rho]$,
with maximum number of predators during $[p_n-L_1,p_n+L_1]$: $\rho>0$ is the (short) time needed
to reach and leave the patch. We assume $L_2>2(L_1+\rho)$, $p_{n+1}-p_n>2(L_1+\rho)$ for all $n\in\N$,
and $p_n-(n-1)\,L_2\to 0$ as $n\to\infty$. The size of the $n$-th group of predators is
determined by the constant $d_n\ge 0$, and we assume that the sequence $(d_n)$ is bounded with limit $d_+$.
The possible differences between $p_n$ and $(n-1)L_2$ capture
variations in the start date of the yearly predation season, and
the hypothesis $p_{n+1}-p_n\to L_2$ is made for the sake of simplicity:
combined with the existence of $d_+$, it describes an asymptotically periodic
phenomenon, which means that
the behavior of the predators becomes as regular as possible over time.
Other more complicated types of recurrence in the future equation may fit in the model.
(See \cite{gosalo} for a study on the variation of arrival dates of common swift and barn
swallow to the Iberian Peninsula.)
The phenomenon of lack of predators in some occasional years can be described through null
elements in the sequence $(d_n)$.
We use the map $\Gamma_{\rho,L}$ of Example \ref{ex:4Cejemplo} (see Figure~\ref{fig:Gammafunction})
to model this behavior: the amount of predators at the ecosystem at time $t$ is
\begin{equation}\label{eq:series}
 \Delta_{(d_n)}(t):=\sum_{n=1}^\infty d_n\, \Gamma_{\rho,L_1}(t-p_n)\,,
\end{equation}
which is a bounded continuous function due to the boundedness of $(d_n)$ and to the
disjointness of the intervals of predation.
Figure~\ref{fig:seriesimpulses} depicts $\Delta_{(d_n)}$ for $\rho=1$, $L_1=10$, $L_2=40$
and certain sequences $(d_n)$ and $(p_n)$.
\par
So, we study the transition equation
\begin{equation}\label{eq:62transitiontoperiodic}
 x'=r(t)\,x\,\left(1-\frac{x}{K(t)}\right)\frac{x-S(t)}{K(t)}-\Delta_{(d_n)}(t)\,\frac{x^2}{b(t)+x^2}\,,
\end{equation}
which represents the dynamics of the single-species population through the repeated passage
(which tends to be periodic) of groups of predators starting at certain time $p_1-L_1-\rho$.
We define $\Delta_-:=0$ and
\[
 \Delta_+(t):=\sum_{n=-\infty}^\infty d_+\Gamma_{\rho,L_1}\big(t-(n-1)L_2\big)\,,
\]
which is bounded, continuous, and $L_2$-periodic in time.
Then, $\lim_{t\to-\infty}(\Delta_{(d_n)}(t)-\Delta_-(t))=0$, since
$\Delta_{(d_n)}(t)=0$ for all $t\leq p_1-L_1-\rho$, and
$\lim_{t\to\infty}(\Delta_{(d_n)}(t)-\Delta_+(t))=0$, since the uniform
continuity of $\G_{\rho,L_1}$ on compact sets ensures that
$\n{\G_{\rho,L_1}(t-p_n)-\G_{\rho,L_1}(t-(n-1)L_2)}_\infty\to 0$ as $n\rightarrow\infty$,
and the separation of the supports of the terms of the series guaranteed by the conditions
$L_2>2(L_1+\rho)$ and $p_n-(n-1)L_2\rightarrow0$ ensures that we can compare the series term-by-term.
That is, \eqref{eq:series} corresponds to a transition between these two limit functions, and
\hyperlink{fd2}{\bf fd2} is fulfilled. It can be checked that the right-hand side
of equation \eqref{eq:62transitiontoperiodic} is not d-concave if
$\max_{t\in\R}\Delta_{(d_n)}(t)=\max_{n\in\N}d_n$ is large enough,
while $r(t)\,x\,(1-x/K(t))(x-S(t))/K(t)-\Delta_+(t)\,x^2/(b(t)+x^2)$ is d-concave
if $d_+$ is not too large, in which case also \hyperlink{fd4}{\bf fd4} is fulfilled:
see \cite[Subsection 5.2]{dno3}.
\par
Let us choose: $r(t):=0.7+0.3\sin^2(t)$, $K(t):=70+20\cos(\sqrt{5}\,t)$ and $S(t)=20+30\cos^2(\sqrt{3}\,t)$
for the internal dynamics of the species, $b(t):=200$ for the influence of the predation,
and $L_1=10$, $L_2=40$, $d_+=0.3$, $d_n=d_++ d/((n-1)/20+1)^2$ and $p_n=(n-1)L_2+(-1)^{(n-1)}/n$ (for all $n\in\N$)
for the shape of the transition function.
The particular expression of $d_n$ implies that the yearly density of predators $d_n$ decreases to $d_+$.
The decreasing attractiveness of the habitat can be indebted to different causes: learning of defensive
mechanisms, overpopulation in the previous season, insufficient nesting or breeding space, etc.
The constant $d$ of the definition of $d_n$ is a size bifurcation parameter
in terms of which we will study the dynamical cases of \eqref{eq:62transitiontoperiodic}.
The choice of $d_+$ (below $0.32$) guarantees \hyperlink{fd4}{\bf fd4}.
We numerically check \hyperlink{fd5}{\bf fd5}, and hence \hyperlink{fd1}{\bf fd1}-\hyperlink{fd5}{\bf fd5} hold
for all $d\ge 0$.
In addition, the size of $d_n$ for small $n$
provides a not d-concave equation \eqref{eq:62transitiontoperiodic} if $d$ is large enough (above $0.96$).
\par
Clearly, $\Delta_{(d_n)}=\tilde\Delta_++d\,\Delta_0$
for $\tilde\Delta_+(t):=\sum_{n=1}^\infty d_+\Gamma_{\rho,L_1}\big(t-p_n\big)$ and
the continuous nonnegative map $\Delta_0(t):=\sum_{n=1}^\infty(1/((n-1)/20+1)^2)\G_{\rho,L_1}(t-p_n)$,
whose limits as $t\to\pm\infty$ are $0$.
We define $f(t,x,\gamma):=r(t)\,x\,(1-x/K(t))(x-S(t))/K(t)-\tilde\Delta_+(t)\,x^2/(b(t)+x^2)-\gamma\,x^2/(b(t)+x^2)$ and
$g(t,x,\gamma):=f(t,x,-\gamma)$, and check that the pairs $(g,d\,\Delta_0)$ satisfy the hypotheses
of Theorem \ref{th:6Dsize2}(ii) (with $\G(t):=0$, $\G_0(t):=\Delta_0(t)$, and $\bar d=0$). To this end,
we numerically check that $x'=g(t,x,0)$ has three hyperbolic copies of the base and that the lower one, attractive, is
$0$ (and hence $\tilde m_0$ is positively bounded from below).
Hence, Theorem~\ref{th:6Dsize} ensures the existence of a unique size-induced
tipping point $d_0>0$ for $x'=f(t,x,d)$ (i.e., for \eqref{eq:62transitiontoperiodic}):
\hyperlink{DCA}{\sc Case A} holds for $0\le d<d_0$, and \hyperlink{DCC}{\sc Case C2} holds for $d>d_0$.
That is, an excessive increase in the number of predators visiting the habitat leads to
the extinction of the species. The existence of this critical transition is depicted
in Figure~\ref{fig:variablend}.
\end{exa}
\begin{exa}
Now, we consider that a flock of $x$ animals described by
\eqref{eq:6populationsimulationgeneral} grazes in a patch which is initially free of predators.
We assume that at time $t=0$ a group of predators, which we suppose that have constant density $d$
(due to the time scale in which we work) and whose predation mechanism is assumed to be suitably
modeled by a H\"{o}lling type III functional response term $-d\,x^2\,/(b(t)+x^2)$
reaches the patch. (See Example~\ref{ex:4Cejemplo} for the meaning of $d$ and $b$.)
The function $b$ is assumed to be quasiperiodic and positively bounded from below.
At time $L>0$, the threat is identified by the flock owner and $s$ shepherds per unit
of time are hired to protect the flock: there are $s\,(t-L)$ shepherds at time $t\geq L$,
and each shepherd is assumed to be able to protect $h$ heads of livestock.
As soon as there are enough shepherds to protect the whole herd, i.e, when $x\le h\,s\,(t-L)$,
predators are not able to attack the flock. That is, predation occurs while $0\le t\le L\,(cx+1)$,
where $c=1/hsL$.
So, for $x\geq0$, we can model the evolution of the flock by the equation
\begin{equation}\label{eq:6Dant}
 x'=r(t)\, x\,\left(1-\frac{x}{K(t)}\right)\frac{x-S(t)}{K(t)}-d\,\G_L\left(\frac{2\,t}{cx+1}-L\right)\,\frac{x^2}{b(t)+x^2}\,,
\end{equation}
where we take $\G_L:=\G_{\rho,L}$ for some small fixed $\rho>0$, with $\G_{\rho,L}$
defined in Example \ref{ex:4Cejemplo} (see Figure \ref{fig:Gammafunction}). So, the predation
term practically vanishes when $t$ is outside the interval $[0,L\,(cx+1)]$.
By multiplying the H\"{o}lling type III functional response term by $\G_L$, it is implicitly
assumed that the search for prey mechanism, i.e. the H\"{o}lling type III interaction, is not affected
by the presence of shepherds as long as there are not enough of them to protect the whole herd.
(See the biological meaning of H\"{o}lling functional response in \cite{holling1}.)
This assumption, made for the sake of simplicity, can be understood as
follows: if a shepherd has more sheep in his care than he can protect,
then a predator, once it has located its prey, can wait a negligible amount of time on the timescale
we are working with until the shepherd moves on to other sheep, far enough away to allow the predator to hunt
the prey.
\par
Since $\R\times[0,\infty)$ is an invariant set for the process given by \eqref{eq:6Dant}
and only nonnegative solutions have biological meaning, we can replace
the predation term by a globally defined one. To this end, we take a
a globally defined $C^2$-map $k(x)$ which coincides with $1/(cx+1)$ on $[0,\infty)$,
and consider the equation
\begin{equation}\label{eq:livestock-k}
 x'=r(t)\, x\,\left(1-\frac{x}{K(t)}\right)\frac{x-S(t)}{K(t)}-d\,\G_L(2\,t\,k(x)-L)\:\frac{x^2}{b(t)+x^2}\,,
\end{equation}
Let $\Lambda_{L,c}(t,x):=\Gamma_L(2\,t\,k(x)-L)$.
Then, for any choices of $d\geq0$, $L>0$ and $c>0$, $d\,\Lambda_{L,c}$ is globally bounded,
$C^2$-admissible on $\R\times\R$, and with $\lim_{t\to\pm\infty}d\,\Lambda_{L,c}(t,x)=0$
uniformly on each compact set $\mJ\subset\R$. That is,
$d\,\Lambda_{L,c}$ globally satisfies \hyperlink{fd2}{\bf fd2}, with $\Lambda_\pm=0$.
In addition, if $f$ is the right-hand term of \eqref{eq:6Ddeff}, then
it is not difficult to check that \hyperlink{fd1}{\bf fd1}, \hyperlink{fd3}{\bf fd3}
and \hyperlink{fd4}{\bf fd4} hold.
\par
We choose $r(t):=0.7+0.3\sin^2(t)$, $K(t):=70+20\cos(\sqrt{5}\,t)$, $S(t):=20+30\cos^2(\sqrt{3}\,t)$,
and $b(t):=20+\cos(t)$ to construct Table~\ref{table:4} and Figures~\ref{fig:dsingletippingpointdconcave},
\ref{fig:dsingletippingpointdconcaveL} and \ref{fig:dsingletippingpointdconcaveC}, and
numerically check that \hyperlink{fd5}{\bf fd5} holds for these choices, being $0$ the lower bounded solution of $x'=f(t,x,0)$. That is $(f,d\,\Lambda_{L,c})$
satisfies \hyperlink{fd1}{\bf fd1}-\hyperlink{fd5}{\bf fd5} for all $d\geq0$,
$L>0$ and $c>0$, and hence the dynamics of \eqref{eq:livestock-k} fits in one of the cases
described by Theorem~\ref{th:6Dcasos}.
Moreover, since $0$ is the lowest bounded solution for the past and future equations,
\hyperlink{DCB}{\sc Cases B1} and \hyperlink{DCC}{\sc C1} are preluded.
In addition, if $g(t,x,\gamma):=f(t,x,-\gamma)$, then
the pairs $(g,d\,\Lambda_{L,c})$ satisfy all the hypotheses of Theorem~\ref{th:6Dsize2}(2)
(with $\G(t):=0$, $\G_0(t):=\Lambda_{L,c}(t)$ and $\bar d=0$).
This result shows the existence of a unique tipping value $d(L,c)>0$:
\eqref{eq:livestock-k}$_d$ is in \hyperlink{DCA}{\sc Case A} for all $d\in[0,d(L,c))$,
in \hyperlink{DCB}{\sc Case B2} for $d=d(L,c)$ and in \hyperlink{DCC}{\sc Case C2} for all $d>d(L,c)$.
Figure~\ref{fig:dsingletippingpointdconcave} depicts the upper locally pullback attractive and
the locally pullback repulsive solutions of the transition equation \eqref{eq:livestock-k}$_d$
for $d$ close to the bifurcation point, for some fixed $L$ and $c$, and
Table~\ref{table:4} shows numerical approximations to $d(L,c)$ for different $L,c>0$.
\begin{figure}
     \centering
         \includegraphics[width=\textwidth]{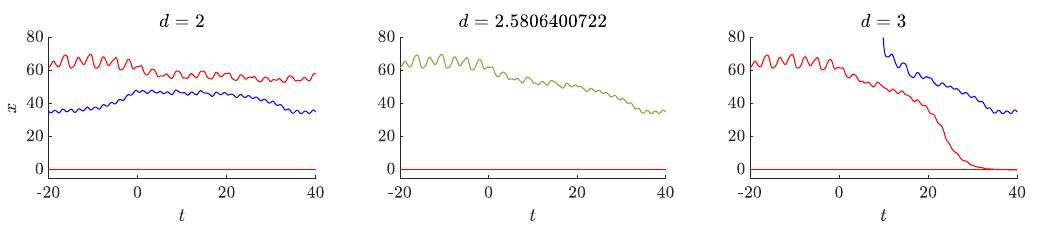}
         \caption{Numerical depiction of the existence of a unique size-tipping point for \eqref{eq:livestock-k}$_d$ when $L$ and $c$ are fixed.
         In this example, $L=20$ and $c=0.02$. The central panel corresponds
         to the approximation for $d(20,0.02)$ of Table ~\ref{table:4}: \protect\hyperlink{DCB}{\sc Case B2}.
         To its left and right, we find \protect\hyperlink{DCA}{\sc Cases A} and
         \protect\hyperlink{DCC}{\sc C2}. See Figure \ref{fig:variablend} to understand the color code.}
        \label{fig:dsingletippingpointdconcave}
\end{figure}
\begin{table}[h!]
\centering
\begin{tabular}{|c | c | c|c|}
 \hline &&&\\[-2.2ex]
 $d(L,c)$  & $c=0.01$ & $c=0.02$ & $c=0.03$\\[0.4ex]
\hline &&&\\[-2.2ex]
$L=\mskip4.5mu 2\mskip4.5mu$  &9.5918417988 &7.8400146619 &6.6406325271\\[0.4ex]
$L=10$  &3.5156887400 &3.1640725896 &2.9522195572\\[0.4ex]
$L=20$ &2.7559336044 &2.5806400722 &2.4622290038\\[0.4ex]
$L=30$ &2.4757094854&2.3677420953&2.3132184604\\[0.4ex]
$L=40$ & 2.3543746813& 2.2850546293&2.2459305139\\[0.4ex]
 \hline
\end{tabular}
\caption{Numerical approximations up to ten places to the bifurcation points $d(L,c)$ of \eqref{eq:livestock-k}$_d$.
The displayed number is a value of $d$ for which \eqref{eq:livestock-k}$_d$ is in
\protect\hyperlink{DCA}{\sc Case A} and such that \eqref{eq:livestock-k}$_{d+1e-10}$ is
in \protect\hyperlink{DCC}{\sc Case C}.
The numerical integration has been done using Matlab2023a \texttt{ode45} algorithm with
\texttt{AbsTol} and \texttt{RelTol} equal to \texttt{1e-12}.
The final integration has been carried out over the interval $[-\texttt{1e4},\texttt{1e4}]$.}
\label{table:4}
\end{table}

Let the other parameters vary. The monotonicity of $L\mapsto \Gamma_L(2t/(cx+1)-L)$ for any
$(t,x)\in\R\times[0,\infty)$ yields the uniqueness of a possible tipping point $L_0$ for
\eqref{eq:livestock-k}$_L$ for $d$ and $c$ fixed.
In fact,
if $\tilde u_L$ and $\tilde m_L$ are the upper and middle hyperbolic solutions if
\eqref{eq:livestock-k}$_L$ is in \hyperlink{DCA}{\sc Case A} and $L_1<L_2$ provide this case,
then Proposition \ref{prop:5Dcoer}(iv) shows that $\tilde u_{L_1}>\tilde u_{L_2}$, and
a new comparison argument shows that $\tilde m_{L_1}\le \tilde m_{L_2}$. So, if
\hyperlink{DCB}{\sc Case B2} (the unique possible one) occurs as $L\uparrow L_0$, then they collide,
and \hyperlink{DCA}{\sc Case A} cannot occur for $L>L_0$.
\par
Analogously, the monotonicity of $c\mapsto \Gamma_L(2t/(cx+1)-L)$ for any
$(t,x)\in\R\times[0,\infty)$ ensures the uniqueness of the bifurcation
for \eqref{eq:livestock-k}$_c$ for $d$ and $L$ fixed in the case of existence.
The biological sense of the problem makes reasonable expecting no more than
one critical transition as $L$ or $c$ varies: the decrease in $L$ means an
earlier detection of the problem and therefore the extinction of the hinders;
and the decrease in $c$ means an increase in the rate of recruitment of shepherds, i.e.,
a faster response to the problem that facilitates survival.
\par
Figures \ref{fig:dsingletippingpointdconcaveL} and \ref{fig:dsingletippingpointdconcaveC}
represent the behaviour of the locally pullback attractive and locally pullback repulsive solutions
of the transition equation \eqref{eq:livestock-k}$_L$ for fixed $d$ and $c$, and \eqref{eq:livestock-k}$_c$
for fixed $d$ and $L$, respectively. As in the case of Figure \ref{fig:dsingletippingpointdconcave},
the left-hand panel corresponds to the survival of the species (\hyperlink{CCA}{\sc Case A}),
the right-hand panel corresponds to extinction (\hyperlink{DCC}{\sc Case C2}), and the middle panel is
an approximation to the intermediate unstable situation between them (\hyperlink{CCB}{\sc Case B2}).

\begin{figure}
     \centering
         \includegraphics[width=\textwidth]{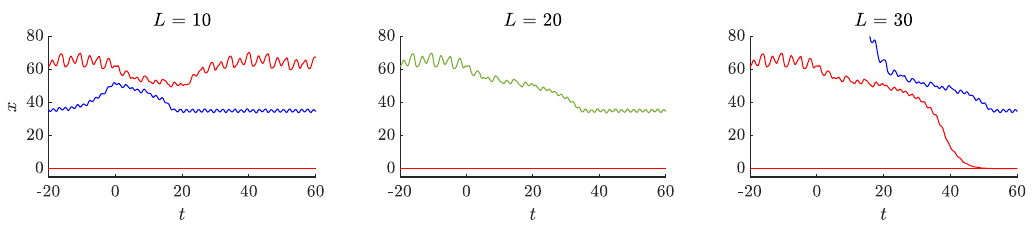}
         \caption{Numerical depiction of an $L$-induced tipping point:
         for $c=0.02$ and $d=2.5806400722$, we find \protect\hyperlink{DCA}{\sc Case A} for $L=10$,
         \protect\hyperlink{DCB}{\sc Case B2} for $L=20$ and \protect\hyperlink{DCC}{\sc Case C2} for $L=30$.}
        \label{fig:dsingletippingpointdconcaveL}
\end{figure}

\begin{figure}
     \centering
         \includegraphics[width=\textwidth]{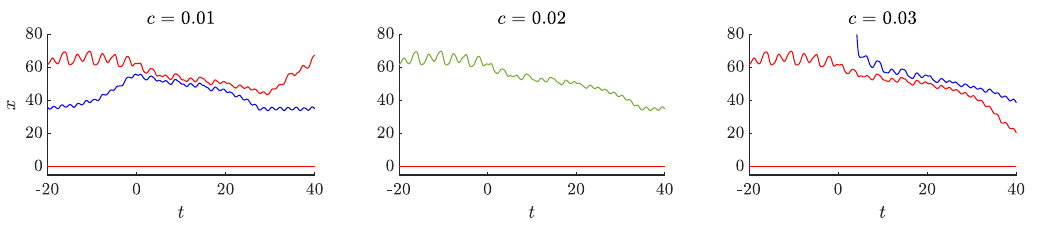}
         \caption{Numerical depiction of a $c$-induced tipping point:
         for $L=20$ and $d=2.5806400722$, we find \protect\hyperlink{DCA}{\sc Case A} for $c=0.01$, \protect\hyperlink{DCB}{\sc Case B2}
         for $c=0.02$ and \protect\hyperlink{DCC}{\sc Case C2} for $c=0.03$.}
        \label{fig:dsingletippingpointdconcaveC}
\end{figure}

\end{exa}
\subsection{Example on the necessity of minimality in Theorem \ref{th:5DCorder}}\label{subsec:ordenes}
In this section, we present an example which shows that the minimality of $(\W,\sigma)$
is indeed required in Theorem~\ref{th:5DCorder}. We will construct a non minimal set $\W$
and a pair of functions $\mh_1,\,\mh_2\colon\WR\to\R$ satisfying
\hyperlink{d1}{\bf d1}-\hyperlink{d4}{\bf d4} with $\mh_1(\w,x)>\mh_2(\w,x)$
for all $(\w,x)\in\WR$, such that $x'=\mh_i(\w{\cdot}t,x)$ has three hyperbolic
copies of the base $\ml_i<\mm_i<\muk_i$ for $i=1,2$ which do not satisfy none
of the two possible orders described in Theorem~\ref{th:5DCorder}.
We make use of the transition framework of Section~\ref{sec:6} to construct the example:
$x'=\mh_i(\w{\cdot}t,x)$ for $i=1,2$ will be transition equations, with $\W$ composed by
a heteroclinic orbit connecting two singletons, which are minimal.
The cornerstone of the example is the fact that we construct three hyperbolic
copies of the base $\W$ contained in $\W\times\R$ projecting onto each one of
the two minimal subsets of $\W$, and the three copies of the base have the two
distinct orders allowed by Theorem~\ref{th:5DCorder} over each minimal.
\par
Let $\G\colon\R\to(0,1)$ be a continuous map with
$\lim_{t\to\infty}\G(t)=\G_+:=1$ and $\lim_{t\to-\infty}\G(t)=\G_-:=0$
(as $\G(t):=\arctan(t)/\pi+1/2$). We take $a\ge\sqrt{10}$ and
\[
 h_b(x,\alpha):=-x^3+x+\alpha\,(3x^2a-3xa^2+a^3-a)+\alpha\,(1-\alpha)\,b\,,
\]
for some $b\ge 0$ which will be properly fixed later on.
Note that: $h_b(x,\alpha)=h_0(x,\alpha)+\alpha(1-\alpha)\,b$;
$h_b(x,0)=-x(x-1)(x+1)$; $h_b(x,1)=-(x-a)(x-a-1)(x-a+1)$; and
$3x^2a-3xa^2+a^3-a\geq 0$ for all $x\in\R$ by the choice of $a$,
so $\alpha\mapsto h_0(x,\alpha)$ is nondecreasing for all $x\in\R$.
For each $b\ge 0$, we consider the equation
\begin{equation}\label{eq:discusionejemplo}
 x'=h_b(x,\G(t))\,.
\end{equation}
It is easy to check that $(h_b,\G)$ satisfies
\hyperlink{fd1}{\bf fd1}-\hyperlink{fd5}{\bf fd5} for any $b\ge0$:
the past equation $x'=h_b(x,0)$ has three hyperbolic critical points $-1$,
$0$ and $1$, and the future equation $x'=h_b(x,1)$, which is a shift of the past one,
has three hyperbolic critical points $a-1$, $a$ and $a+1$.
So, the dynamics of \eqref{eq:discusionejemplo}$_b$ fits in one of the
dynamical cases of Theorem~\ref{th:6Dcasos}.
\par
We will check later the existence of $b_0>0$ such that
\eqref{eq:discusionejemplo}$_b$ is in \hyperlink{DCA}{\sc Case A} for $b=b_0$.
Let $\W$ be the hull of $(t,x)\mapsto h_{b_0}(x,\G(t))$ (see Subsection~\ref{subsec:hull}),
and let $\mh_1\colon\WR\to\R$ be given by $\mh_1(\w,x):=\w(0,x)$ for $(\w,x)\in\WR$,
that is, the extension of $h_{b_0}$ to $\W$. Then, $\mh_1(\w,x)$ is a cubic polynomial
with $-1$ as leading coefficient for all $\w\in\W$, and hence $\mh_1$ satisfies
\hyperlink{d1}{\bf d1}-\hyperlink{d4}{\bf d4}. Note that $\W$ is the union of the
(heteroclinic) $\sigma$-orbit $\{h_{b_0}(x,\G(t+s))\,|\;s\in\R\}$ and its \upalfa-limit and
\upomeg-limit sets, $\{h_{b_0}(x,0)\}$ and $\{h_{b_0}(x,1)\}$:
see Lemma \ref{lema:2hull}. Theorem~\ref{th:5DCthreecopies} ensures that
$x'=\mh_1(\wt,x)$ has three hyperbolic copies of the base $\ml_1<\mm_1<\muk_1$.
In particular,
the restrictions of these three copies to the \upalfa-limit set $\{h_{b_0}(x,0)\}$ are
$-1,\,0$ and $1$, and to the \upomeg-limit set $\{h_{b_0}(x,1)\}$ are $a-1,\,a$ and $a+1$.
\par
Next, we define $\mh_2(\w,x):=-x^3+x-\ep$ for $\ep\in(0,2/(3\sqrt{3}\,))$, which
clearly satisfies \hyperlink{d1}{\bf d1}-\hyperlink{d4}{\bf d4} and
$\mh_1(\w,x)>\mh_2(\w,x)$ for all $(\w,x)\in\WR$. It can be checked that
$x'=\mh_2(\w{\cdot}t,x)$ has three copies of the base: three constant equilibria
$\ml_2,\mm_2,\muk_2$ satisfying $\ml_2<-1<0<\mm_2<\muk_2<1$.
So, the order of $\ml_1$, $\mm_1$, $\muk_1$ and $\ml_2$, $\mm_2$, $\muk_2$
is $\ml_2<-1<0<\mm_2<\muk_2<1$ (like in Theorem~\ref{th:5DCorder}(i))
over the minimal set $\{h_{b_0}(x,0)\}\subset\W$, and
$\ml_2<\mm_2<\muk_2<a-1<a<a+1$ (like in Theorem~\ref{th:5DCorder}(ii))
over the minimal set $\{h_{b_0}(x,1)\}\subset\W$.
Hence, the continuity of the copies of the base preclude
any of the two possibilities of Theorem~\ref{th:5DCorder} to hold over the whole $\W$.
\par
It remains to check the existence of $b_0>0$ such that \eqref{eq:discusionejemplo}$_{b_0}$
is in \hyperlink{DCA}{\sc Case A}, for which it suffices to check \eqref{eq:discusionejemplo}$_0$
is in \hyperlink{DCC}{\sc Case C2} and that there exists $b_1>0$ such that
\eqref{eq:discusionejemplo}$_{b_1}$ is in \hyperlink{DCC}{\sc Case C1}:
Theorem \ref{th:6Dsaddle-node} precludes moving from \hyperlink{DCC}{\sc Case C2} to \hyperlink{DCC}{\sc Case C1}
as $b$ varies without crossing \hyperlink{DCA}{\sc A}.
\par
We denote by $l_b$ and $u_b$
(resp.~$m_b$) the locally pullback attractive (resp.~repulsive) solutions of
\eqref{eq:discusionejemplo}$_b$ provided by Theorem \ref{th:6DCtressoluciones},
and recall that $\lim_{t\to-\infty} u_b(t)=1$, $\lim_{t\to-\infty} l_b(t)=-1$,
and $\lim_{t\to\infty} m_b(t)=a$.
Since $\G(t)\le 1$ for all $t\in\R$ and $\alpha\mapsto h_0(x,\alpha)$
is nondecreasing for all $x\in\R$, we have $h_0(a-1,\G(t))\leq h_0(a-1,1)=0$
for all $t\in\R$, so $\R\times(-\infty,a-1]$ is positively invariant
for \eqref{eq:discusionejemplo}$_0$. Since $\lim_{t\to-\infty} u_0(t)=1<a-1$,
we have $u_0(t)\in (-\infty,a-1]$ for all $t\in\R$, and hence
$\lim_{t\to\infty}u_0(t)=a-1$: the other possible future limits $a$ and $a+1$
are uniformly separated from $u_0$.
That is, \eqref{eq:discusionejemplo}$_0$ is in \hyperlink{DCC}{\sc Case C2}.
To look for $b_1$, we first check that all the bounded solutions of
\eqref{eq:discusionejemplo}$_b$ take values in $[-1,\infty)$, since
$h_b(x,0)<h_b(x,\G(t))$ for all $(t,x)\in\RR$, and hence
any $m_1<-1$ satisfies the initial hypothesis of Proposition \ref{prop:5Dcoer}.
Next, we take $t_0>0$ in the domain of definition of $m_0$ with $m_0(t)<a+1/2$ for all $t\ge t_0$
and assume for contradiction that $l_b(t)\le a+1/2$ for all $b>0$ and $t\in[t_0,t_0+1]$.
Let $\gamma$ be a lower bound for $\G(t)(1-\G(t))$ for $t\in[t_0,t_0+1]$.
Then $l_b(t_0+1)\ge-2+\int_{t_0}^{t_0+1}(-(a+1/2)^3+\gamma\,b)\,ds$ for all $b>0$, which
is impossible. We take $b_0$ and $t_1$ with $l_{b_0}(t_1)>a+1/2>m_0(t_1)$. Theorem~\ref{th:6DCtressoluciones}
ensures that $\lim_{t\to\infty}(x_0(t,t_1,l_{b_0}(t_1))-(a+1))=0$,
and a comparison argument yields $l_{b_0}(t)=x_{b_0}(t,t_1,l_{b_0}(t_1))\ge x_0(t,t_1,l_{b_0}(t_1))$
for $t\ge t_1$. That is, $\liminf_{t\to\infty}(l_{b_0}(t)-(a+1))\ge 0$, which may only happen
in \hyperlink{DCC}{\sc Case C1} (see Theorem~\ref{th:6Dcasos}). This completes the proof.

\end{document}